\documentclass[twoside,a4paper,reqno,11pt]{amsart} 

\usepackage{enumerate}

\usepackage[margin=2cm]{geometry}
\usepackage[
backend=biber,
style=numeric,
sorting=ynt
]{biblatex}

\numberwithin{equation}{section}

\usepackage{amsfonts, amsmath, amssymb, mathrsfs, bm, latexsym, stmaryrd, array, hyperref, mathtools, bold-extra, xcolor, comment, parskip, dirtytalk, lipsum}
\usepackage{booktabs}

\counterwithin{table}{section}

\theoremstyle{plain}
\newtheorem{theorem}{Theorem}[section]

\newtheorem{claim}[theorem]{Claim}
\newtheorem{lemma}[theorem]{Lemma}
\newtheorem{definition}[theorem]{Definition}
\newtheorem{corollary}[theorem]{Corollary}
\newtheorem{prop}[theorem]{Proposition}
\newtheorem{example}[theorem]{Example}

\newtheorem{remark}[theorem]{Remark}
\newtheorem*{theorem*}{Theorem}
\newtheorem*{definition*}{Definition}
\newtheorem*{prop*}{Proposition}
\newtheorem*{lemma*}{Lemma}
\newtheorem*{corollary*}{Corollary}
\newtheorem*{example*}{Example}
\newtheorem{hypothesis}[theorem]{Hypothesis}
\newtheorem*{hypothesis*}{Hypothesis}

\begin{document}

\title{On the Orbital Diameter of Primitive Affine Groups}
\author{Kamilla Rekv\'enyi }
\address{Department of Mathematics,
University of Manchester,
Manchester,
M13 9PL \newline Heilbronn Institute for Mathematical Research, Bristol, UK}
\email{kamilla.rekvenyi@manchester.ac.uk}
\date{September 2024}

\begin{abstract}
     The orbital diameter of a primitive permutation group is the maximal diameter of its orbital graphs. There has been a lot of interest in bounds for the orbital diameter. In this paper we provide explicit bounds on the diameters of groups of primitive affine groups with almost quasisimple point stabilizer. 
    As a consequence we obtain a partial classification of primitive affine groups with orbital diameter less than or equal to 2. 
\end{abstract}

\maketitle

\section{Introduction}
Let $G\leq Sym(\Omega)$ be a transitive permutation group on a finite set $\Omega.$ We can define a componentwise action of $G$ on $\Omega \times \Omega.$ An orbital is an orbit of $G$ on $\Omega\times\Omega$. There is a unique diagonal orbital $\Delta = \{(\alpha,\alpha):\alpha \in \Omega\},$ all others are called non-diagonal orbitals. For a non-diagonal orbital $\Gamma$ we define the orbital graph $\Gamma$ to be an undirected graph with vertex set $\Omega$, and edge set the pairs in the orbital $\Gamma$. A theorem of Donald Higman states that the non-diagonal orbital graphs are all connected if and only if the action of $G$ is primitive. Hence we can define the \textit{orbital diameter} of a primitive permutation group to be the supremum of the diameters of its orbital graphs. We will denote this by $orbdiam(G,\Omega).$ \par
Now we describe primitive affine groups. 
Let $V=V_n(q_0)$ be an $n-$dimensional vector space over $\mathbb{F}_{q_0}$. Then $A\Gamma L(V)$ is the set of permutations of $V$ of the form $x \to Ax+b,$ where $A\in \Gamma L(V)$ and $b\in V.$ 
A primitive group of affine type has the form $G = VG_0\leq A\Gamma L(V),$  where $G_0 \leq \Gamma L(V ),$ the stabilizer of $0,$ acts irreducibly on $V$ and $V$ acts by translations. \par 
\par In \cite{orbdiam}, Martin Liebeck, Dugald MacPherson and Katrin Tent classified infinite families of primitive permutation groups such that there is an upper bound on the orbital diameter of all groups in the family. Their motivation and methods of proof were model theoretical and they provided no explicit bounds on the orbital diameter. Hence two natural goals in the study of the orbital diameters are to find explicit bounds and to classify groups with small orbital diameter. 
In this paper we fulfil these goals for primitive affine groups. We provide explicit lower bounds for the orbital diameter, and using these and further study we provide an overview of primitive affine groups of orbital diameter at most 2. This complements the results in \cite{attilaskresanov} and  \cite{skresanov}, which provide some explicit upper bounds on the orbital diameter of primitive affine groups. In fact, in \cite{skresanov} the author proves that for a finite primitive affine permutation group, $G\leq AGL_n(p)$ such that $p$ divides $\vert G_0\vert,$ the orbital diameter is bounded above by $ 9 n^3.$

 Let $C$ be an infinite class of finite affine primitive groups and suppose $C$ is a bounded class, i.e. there is some $d$ such that $orbdiam(G,V)\leq d$ for all $G\in C$.  Theorem 1.1 in \cite{orbdiam} states that for such a class $C$,  all $G\in C$ are of $t$-bounded classical
type, defined as follows, for some $t$ bounded by a function of $d$. We will denote a quasisimple classical group that has natural module $V_n(q_0)$ by $Cl_n(q_0).$
\begin{definition}\textup{\cite{orbdiam}}
An affine primitive group $G=VG_0$ where $V=V_n(q_0),$ an $n$-dimensional vector space over $\mathbb{F}_{q_0},$ and $G_0\leq \Gamma L_n(q_0)$ is of $t$-bounded classical type if both of the following hold. \begin{itemize}

       \item $G_0$ stabilizes a direct sum decomposition $ V_1\bigoplus\dots\bigoplus V_k$ of $V$ and acts transitively on the set $\{V_1,\dots,V_k\},$ where $k\leq t.$
        \item  There is a tensor decomposition  $V_1=V_m(q_0)\otimes_{\mathbb{F}_{q_0}}Y$ where $\dim Y\leq t.$ The group $G_1$ induced by $G_0$ on $V_1$ contains $Cl_m(q_0')\otimes 1_Y\unlhd G_1$ acting naturally on $V_1,$ where $\vert \mathbb{F}_{q_0}:\mathbb{F}_{q_0'}\vert \leq t.$
    
\end{itemize}
\end{definition}


A theorem of Aschbacher  \cite{aschbacher} characterises the subgroup structure of classical groups. It says that any subgroup of $\Gamma L(V)$ either lies in one of $8$ well-understood geometric classes, denoted $C_1-C_8$ or is almost quasisimple, absolutely irreducible on $V$ and not realizable over a subfield of $\mathbb{F}_{q_0}.$ We focus our attention on this last class, the so-called $\mathcal{S}$ or $C_9$ class. 
An integral part of the proof of \cite[Thm 1.1 (1)]{orbdiam} deals with the case when $G_0$ is almost quasisimple. Note that $G_0^\infty$ denotes the last term of the derived series of $G_0.$

\begin{prop}\textup{\cite[Prop 3.6]{orbdiam}}\label{f(d)tetel}
Fix $d\in \mathbb{N}$. Let $G=V_n(q_0) G_0\leq A\Gamma L_n(q_0)$ and suppose $G_0^\infty$ is quasisimple and acts absolutely irreducibly on $V_n(q_0).$ Also suppose that $orbdiam(G,V)\leq d.$ Then there exists a function $f\colon \mathbb{N} \to \mathbb{N}$ such that one of the following holds; \begin{enumerate}
    \item $n\leq f(d)$
    \item $G_0^\infty = Cl_n(q_0')$, where $\vert \mathbb{F}_{q_0}:\mathbb{F}_{q_0'}\vert \leq d.$
\end{enumerate}
\end{prop}

Note that here the second case is an example of a primitive affine group of $t$-bounded classical type. The first case restricts the dimension of the vector space by bounding it by a function of the diameter. In \cite{orbdiam} the function $f(d)$ was not explicitly determined. In this paper we provide explicit lower bounds for $d$ for  groups as in the following hypothesis.\par

 
\begin{hypothesis}\label{hypothesisgen}
  Let $G=VG_0$ be a primitive affine group such that $G_s:=\frac{G_0^\infty}{Z(G_0^\infty)}$ is a non-abelian finite simple group. Suppose that $V$ is an absolutely irreducible $\mathbb{F}_{q_0}G_0^\infty-$module in characteristic $p.$ Also let $n$ be the dimension of $V$  and assume that $V$ cannot be realised over a proper subfield of $\mathbb{F}_{q_0}.$ 
\end{hypothesis}
 When $G_s$ is of Lie type, we found lower bounds for $d,$ which are expressed as a function of the Lie rank of $G_s$ and a function of the degree for $G_s=A_n.$ We use these lower bounds to determine the affine groups of orbital diameter $2.$


\subsection{Lie type Stabilizer in Defining Characteristic}
Here we give a lower bound on the orbital diameter for the case when $G$ is as in Hypothesis \ref{hypothesisgen} and $G_0^\infty=X_l(q),$ where $X_l(q)$ is a finite quasisimple group of Lie type of Lie rank $l$ over a field $\mathbb{F}_q$ in characteristic $p.$ (Here we have $l$ as the rank of the ambient algebraic group.)
The case when $G_0$ is a classical group and $V$ is its natural module is covered in Lemma \ref{naturalmod}.

\begin{theorem}\label{lhalf}
Let $G$ be as in Hypothesis \ref{hypothesisgen} such that  $G_0^\infty=X_l(q),$ where $X_l(q)$ is a finite quasisimple group of Lie type in characteristic $p.$  Assume that if $G_0^\infty$ is a classical group then $V$ is not a natural module for $G_0^\infty.$ \par 
 Then $$orbdiam(G,V)\geq \lfloor\frac{l}{2}\rfloor.$$ Moreover, for $n>(2l+1)^2,$  $$orbdiam(G,V)\geq \frac{l^2}{18}.$$ 
\end{theorem}

\underline{Note}  Some classical groups are isomorphic to others and hence have several \say{natural modules.} For example $PSL_2(q)\cong \Omega_3(q)$ has natural modules of dimension $2$ and $3.$ In Theorem \ref{lhalf}, for all such classical groups, all their natural modules are excluded. A complete list of such isomorphisms can be found in \cite[p. 96]{wilsonbook}. 

Using Theorem \ref{lhalf}, we achieve a classification of such groups with orbital diameter at most 2. Our results are described in Table \ref{tab:defchar}, where we use the notation $V=V(\lambda)$ for the highest weight module with highest weight $\lambda,$ as defined in Section 2.

\begin{theorem} \label{diam2}
 Let $G$ be as in Hypothesis \ref{hypothesisgen} such that  $G_0^\infty=X_l(q),$ where $X_l(q)$ is a finite quasisimple group of Lie type in characteristic $p.$ Suppose $orbdiam(G,V)\leq 2.$ Then one of the following holds. 
 \begin{enumerate}
     \item $G_0$ is as in Table \ref{tab:defchar}. Moreover, under the assumption that $G_0$ contains the group $\mathbb{F}_{q_0}^*$ of scalars, the permutation rank $r$ and the orbital diameter $d$ are as in Table \ref{tab:defchar}.
     \item $G_0^\infty\cong \,^2D_5(q),$ $(\lambda,\dim V(\lambda))=(\omega_5,16)$ with $q_0=q^2,$ or $G_0^\infty\cong ^3D_4(q),$ $(\lambda,\dim V(\lambda))=(\omega_1,8)$ with $q_0=q^3.$
 \end{enumerate}

\end{theorem}
 \begin{table}
     \centering
 \[
 \begin{array}{c|c|c|c|c|c}
 G_0\,\, \triangleright &\lambda&\dim(V(\lambda))&\text{extra conditions}&r&d\\
 \hline
 \text{classical}&\omega_1&\dim(V(\omega_1))&&2 \,\,\text{or}\,\, 3&1 \,\,\text{or}\,\, 2\\
 \hline
 A_4(q)&\omega_2&10&& 3&2\\\hline
  G_2(q)&\omega_1&6&q \,\,\text{even}& 2&1\\\hline
   G_2(q)&\omega_1&7&q \,\,\text{odd}& 4&2\\\hline
   D_5(q)&\omega_5&16&& 3&2\\\hline 
    B_4(q)&\omega_4&16&& 4&2\\\hline
     B_3(q)&\omega_3&8&& 3&2\\\hline
      ^2B_2(q)&\omega_1&4&& 3&2\\\hline 
 \hline
 \end{array}
 \]
     \caption{Small orbital diameter cases in defining characteristic}
     \label{tab:defchar}
 \end{table}

\begin{remark} \rm
In part (2.) we have not been able to determine whether the orbital diameter is 2.  We conjecture that it is at least $3.$
\end{remark}

\subsection{Alternating Stabilizer}
Here we provide a lower bound on the orbital diameter for the case when $G$ is as in Hypothesis \ref{hypothesisgen} and $G_s\cong A_r.$ We start with a definition.

\begin{definition}\label{fullydeleted}
Let $G=A_r$ or $S_r$ and let $\mathbb{F}_{q_0}^r$ be the permutation module of $G$ over $\mathbb{F}_{q_0}.$ Define submodules $W:=\{(a_1,\dots,a_r):\sum a_i=0\}\leq \mathbb{F}_{q_0}^r$ and $T:=Span(1,\dots,1).$ The fully deleted permutation  module is ${W}/{W\cap T}$. The fully deleted permutation  module has dimension $r-1$ if $p \nmid r$  and $r-2$ if $p\vert r.$
\end{definition}

Now we provide an asymptotic upper bound for $n$ of the form $f(d)$ as in Proposition \ref{f(d)tetel}.
\begin{theorem}\label{fdtheoremalt}
 Fix $\epsilon>0$. Let $G$ be as in Hypothesis \ref{hypothesisgen} such that $G_s\cong A_r$ and $d=orbdiam(G,V).$  Assume $V$ is not the fully deleted permutation module. 
Then there exists an $R\in \mathbb{N}$ such that for all $r\geq R$ we have $n\leq d^{2+\epsilon},$ where $n$ is the dimension of $V.$
 
 \end{theorem}

 The following result concerning the case when $V_n(q_0)$ is the fully deleted permutation  module 
 gives an explicit linear function $f(d)$ as in Proposition \ref{f(d)tetel} as an upper bound for $n.$

\begin{prop}\label{permalt}
Let $G=VG_0$ be as in Hypothesis \ref{hypothesisgen} such that $G_s\cong A_r$ and $V$ be the fully deleted permutation module.
Let $d$ denote the orbital diameter of $G.$ Then 
\begin{enumerate}[(i)]
    \item if $p\nmid r,$ then $orbdiam(G,V)\geq \frac{r-1}{2}$
    \item if $p\nmid r,$ then $orbdiam(G,V)\geq \frac{r-2}{4}.$
\end{enumerate}
\end{prop} 

Using this, we provide the following classification.

\begin{corollary}\label{classification}
 Let $G=V_n(q_0)G_0$ and suppose $G_0^\infty$ is quasisimple with $G_s\cong\frac{G_0^\infty}{Z(G_0^\infty)}\cong A_r$ and $V_n(q_0)$ is the fully deleted permutation  module. Then
 \begin{enumerate}
     \item $orbdiam(G,V)=1$ if and only if $r=6$ and $q_0=2.$
     \item $orbdiam(G,V)=2$ if and only if one of the following holds: \begin{enumerate}
         \item $q_0=2$ and $r=5,8$ or $10.$ 
         \item $q_0=3$ and $r=6.$
         \item $q_0=5,$ $r=5$ and $4\times A_5\leq G_0.$
     \end{enumerate}
 \end{enumerate}
\end{corollary}

 We also provide an explicit lower bound for the orbital diameter for the cases when $V$ is not the fully deleted permutation module, which we will use in our classification of groups with orbital diameter 2.
 
 \begin{theorem}\label{lowerdiamalt}
 Let $G$ be as in Hypothesis \ref{hypothesisgen} such that $G_s\cong A_r$ and $d=orbdiam(G,V).$ Assume that $V$ is not the fully deleted permutation module. Then one of the following holds.
\begin{itemize}
         \item  $r\geq 15$ and $d\geq \frac{r^{2}-5r-2}{2r\log_2(r)}\geq \frac{r-6}{2\log_2(r)}.$
         \item  $r\leq 14$ and a lower bound is as follows: \begin{tabular}{c|c|c|c|c|c|c|c|c|c|c|c}
             $r$ &5&6&7&8&9&10&11&12&13&14  \\\hline
              $d\geq$&1&1&1&1&2&2&2&2&2&2
         \end{tabular}
     \end{itemize}  
\end{theorem}

We also explore the possibilities for cases with small orbital diameter when $V$ is not the fully deleted permutation module and get the following result. 

\begin{theorem}\label{if2}
Let $G$ be as in Hypothesis \ref{hypothesisgen} such that $G_s\cong A_r$. Assume $V$ is not the fully deleted permutation module.  If $orbdiam(G,V)\leq 2$ then one of the following holds: 
\begin{enumerate}
\item $r\leq 7$ and $n\leq 9.$
    \item
    $(r,n,q_0)$ is as in the following table: \[\begin{array}{c|c|c}
    r & n&q_0 \\ \hline
     12&16&4\\ \hline
     11&16&4\\ \hline
     11&16&5\\ \hline
     9&8&2\\ \hline
          8&4&2\\ \hline
\end{array}\]
    
\end{enumerate}

\end{theorem}

We provide some examples as in parts $(1.)$ and $(2.)$ for groups with small orbital diameters in Example \ref{partialconverse}. \par 

\subsection{Lie type Stabilizer in Cross Characteristic}

Here we list our results for the case when $G$ is as in Hypothesis \ref{hypothesisgen} and $G_0^\infty = X_l(r),$ a quasisimple group of Lie type such that $(r,p)=1.$ \par

\underline{Remark} In this section, for $PSp_{2l}(r)$ we assume that $l\geq 2$ and for $P\Omega^\epsilon_s(q_0)$ we assume that $s\geq 7$ as for the smaller values of $l$ and $s$ they are isomorphic to other classical groups that we cover.  \par

We start with an asymptotic upper bound for $n$ of the form $f(d)$ as in Proposition \ref{f(d)tetel}.
 
\begin{theorem}\label{fdtheoremlie}
Fix $\epsilon>0$. Let $G$ be as in Hypothesis \ref{hypothesisgen} such that $G_0^\infty= X_l(r),$ a quasisimple group of Lie type such that $(r,p)=1$ and let $d=orbdiam(G,V).$  There is an $R\in \mathbb{N}$ such that if $\vert X_l(r)\vert\geq R,$  then $n\leq d^{1+\epsilon}.$

 \end{theorem} 
 
 We also provide a lower bound for the orbital diameter. 
 
 \begin{theorem}\label{lowerdiamlie}
 Let $G$ be as in Hypothesis \ref{hypothesisgen} such that $G_s\cong X_l(r),$ a quasisimple group of Lie type such that $(r,p)=1$ and let $orbdiam(G,V)= d.$ 
 \begin{enumerate}
     \item If $ X_l(r)$ is either an untwisted exceptional group of Lie rank $l$ or $X_l(r)\cong \,^2E_6(r)$ or $^2\!F_4(r)$, then one of the following holds:\begin{itemize} 
         \item $d\geq \frac{r^l}{l\log_2(r)}$
         \item $X_l(r)\cong $  $^2\!F_4(2)',$  $G_2(3),$ $G_2(4)$ or $F_4(2)$ and $d\geq 2,$ or $X_l(r)\cong G_2(5)$ and $d\geq 4,$ or $X_l(r)\cong G_2(7)$ and $d\geq 8.$
     \end{itemize}
     \item If $ X_l(r)\cong\, ^2B_2(r),$ $^2G_2(r)$ or $^3D_4(r),$ then one of the following holds:
     \begin{itemize} 
         \item $d\geq \frac{r^{l-1}}{{(l-1)}\log_2(r)}$
         \item $X_l(r)\cong $   $^2B_2(8),$  or $^3D_4(2),$ and $d\geq 2,$ or $X_l(r)\cong ^2B_2(32),$ and $d\geq 5.$ 
     \end{itemize}
     \item If $X_l(r)$ is a classical group of Lie rank $l,$ then 
         $$d\geq \frac{r^{l}-3}{(l+1)^3\log_2(r)}.$$

\end{enumerate}
 \end{theorem}

These can be used to study such groups with orbital diameter at most 2, but many possibilities are unresolved - for more details see \cite{kalathesis}.
\subsection{Sporadic Stabilizer}
 Here we list our results for the case when $G$ is as in Hypothesis \ref{hypothesisgen} and $G_s$ is a  sporadic simple group. 
 \begin{theorem}\label{sporadic}
     Let $G$ be as in Hypothesis \ref{hypothesisgen} and $G_s$ a  sporadic simple group. \begin{enumerate}[(i)]
         \item If $G_s$ is $M,$ $BM,$ $Ly,$ $HN,$ $Th,$ $O'N,$ $Fi_{24}'$ or $Fi_{23},$ then $orbdiam(G,V)\geq 3.$ 
         \item If $n> N,$ where $N$ is as in the following table, then $orbdiam(G,V)\geq 3.$ 
         \[ {\begin{array}{c|c|c|c|c|c|c|c|c|c|c|c|c|c|c|c|c|c|c}
           G_s   & M_{11}&M_{12}&M_{22}&M_{23}&M_{24}&J_1&J_2&J_3&J_4&HS&McL&He&Ru&Suz&Co1&Co2&Co3&Fi22 \\\hline
           N&11&15&34&44&44&20&36&18&112&22&22&51&28&12&24&23&23 &78 
         \end{array}}\]
         \item If $(G_s,n,q_0)$ are as in the following table and $G_0$ contains the scalars in $GL_n(q_0),$ then $orbdiam(G,V)=2.$ 
        \[\begin{array}{c|c|c|c|c|c}
           G_s   & M_{11}&M_{24}&Suz&J_2&J_2 \\ \hline
            n  &  5&11&12&6&6\\ \hline
            q_0& 3&2&3&4&5 \\ \hline
           G_s\leq & PSL_5(3)& PSL_{11}(2)& PSp_{12}(3)& G_2(4)\leq Sp_6(4)&PSp_6(5)
         \end{array}\]
     \end{enumerate}
 \end{theorem}

\section{Preliminary Results}

In the proofs of our results we extensively use facts from the representation theory of groups of Lie type. We include those and some preliminary lemmas about the orbital diameter of primitive affine groups.

\subsection{Representations of the finite groups of Lie type in defining characteristic}

Let $\overline{L}$ be a simple, simply connected algebraic group over $\overline{\mathbb{F}_p}$ ($p$ prime), and let $L=\overline{L}^F,$ where $F$ is a Frobenius morphism, so let $L=X_l(q)$ be a finite group of Lie type over a finite field $\mathbb{F}_q$. Let $\Pi=\{\alpha_1,\dots,\alpha_l\}$ be a system of fundamental roots for $\overline{L}$ and $\omega_l,\dots,\omega_l$ be the corresponding fundamental dominant weights. Let $X_q=\{\sum_i c_i\omega_i\vert 0\leq c_i\leq q-1\}$ unless  $L=^2B_2(q),$ $^2G_2(q)$ or $^2F_4(q).$ In the latter cases, 
let $X_q=\{\sum_i c_i\omega_i\vert 0\leq c_i\leq q(\alpha_i)-1\},$ where $q=q^{2a+1}$ with $p=2,$ $3,$ $2,$ respectively, and $q(\alpha)=p^a$ if $\alpha$ is a long root and $q(\alpha)=p^{a+1}$ if $\alpha$ is a short root, and  otherwise. Let $V(\lambda)$ be the irreducible $\overline{L}$-module of highest weight $\lambda.$ For $\overline{L}$ of type $A_l(k),$ $D_l(k),$ $D_4(k)$ and $E_6(k),$ let $\tau_0$ denote a graph automorphism of $\overline{L}.$

\begin{theorem}\textup{\cite[Thm 5.4.1 and Remark after Thm 5.4.1]{kleidmanliebeck}}\label{thm541}
A above, let $L$ be a simply connected group of Lie type over $\mathbb{F}_q.$  Then  for $\lambda\in X_q$ the modules $V(\lambda)$ remain irreducible and inequivalent upon restriction to $L$ and exhaust the irreducible $kL$-modules.  
\end{theorem}

For $J\subseteq \Pi,$  let $P_J$ be the parabolic subgroup of $\overline{L}$ corresponding to deleting the nodes in $J$ from the Dynkin diagram of $\overline{L}.$
\begin{definition}\label{parabolicremark}
   For a dominant weight $\lambda=\sum a_i\omega_i,$ we define $P_\lambda$ to be the parabolic $P_J,$ where $J=\{i:a_i\neq 0\}.$

\end{definition}

\begin{lemma}\label{para}\textup{\cite[2.3]{liebeckaffine}}
Let $V(\lambda)$ the irreducible highest weight module with highest weight $\lambda=\sum_{i=1}^r a_i\omega_i$ and $v^+$ be a maximal vector fixed by a Borel subgroup of $L=X_l(q)$. Then one of the following holds. 
\begin{enumerate}[(i)]
    \item For $X_l(q)$ untwisted, the stabilizer in $L$ of $\langle v^+\rangle$ is the parabolic subgroup $P_\lambda^F.$
    \item For $X_l(q)$ twisted, the stabilizer in $L$ of $\langle v^+\rangle$ is the parabolic subgroup ${(P_{\lambda+\tau_0(\lambda)})}^F$ except for $X_l(q)=\,^3D_4(q)$ in which case it is ${(P_{\lambda+\tau_0(\lambda)+\tau_0^2(\lambda)})}^F.$
\end{enumerate}
\end{lemma}

The following lemma characterises the fields over which the absolutely irreducible representations for groups of Lie type in defining characteristic are defined.  \par 

\begin{lemma}\label{q0lemma}
    Let $X_l(q)$ be a quasisimple group of Lie type and $V=V_n(q_0)$ an absolutely irreducible module for $X_l(q)$ in defining characteristic that cannot be realized over a proper subfield of $\mathbb{F}_{q_0}$. 
    \begin{enumerate}[(i)]
        \item \label{qnode1} If  $X_l(q)$ is untwisted, or $X_l(q)$ is $\,^2\!A_l(q),$ $\,^2\!D_l(q),$ $\,^2\!E_6(q)$ or $\,^3\!D_4(q)$ with $V\cong V^{\tau_0},$ then one of the following holds:
  
    \begin{enumerate}
        \item $q=q_0.$
        \item \label{subfield}$q=q_0^k,$ $k\geq 2$ and $V=V(\lambda)\otimes V(\lambda)^{q_0}\otimes \dots\otimes V(\lambda)^{q_0^{k-1}}$ realised over $\mathbb{F}_q$ for some $\lambda.$ 
     
    \end{enumerate}
        \item \label{twistedqhatvany} If  $X_l(q)$ is $\,^2\!A_l(q),$ $\,^2\!D_l(q),$ $\,^2\!E_6(q)$ or $\,^3\!D_4(q)$ with $V\not \cong V^{\tau_0}$ then one of the following holds: \begin{enumerate}
            \item $^2\!A_l(q)\leq A_l(q^2)\leq GL(V)$ and $V$ satisfies (\ref{qnode1}) for $A_l(q^2).$
            \item $^2\!D_l(q)\leq D_l(q^2)\leq GL(V)$ and $V$ satisfies (\ref{qnode1}) for $D_l(q^2).$
            \item $^2\!E_6(q)\leq E_6(q^2)\leq GL(V)$ and $V$ satisfies (\ref{qnode1}) for $E_6(q^2).$
            \item $^3\!D_4(q)\leq D_4(q^3)\leq GL(V)$ and $V$ satisfies (\ref{qnode1}) for $D_4(q^3).$
        \end{enumerate}
        \item \label{twistedexceptional} If $X_l(q)$ is $\,^2\!B_2(q),$ $\,^2\!G_2(q)$ or $\,^2\!F_4(q),$ then one of the following holds: \begin{enumerate}
            \item  $^2\!B_2(q)\leq B_2(q)\leq GL(V)$ and $V$ satisfies (\ref{qnode1}) for $B_2(q).$
            \item  $^2\!G_2(q)\leq G_2(q)\leq GL(V)$ and $V$ satisfies (\ref{qnode1}) for $G_2(q).$
            \item  $^2\!F_4(q)\leq F_4(q)\leq GL(V)$ and $V$ satisfies (\ref{qnode1}) for $F_4(q).$
        \end{enumerate}
      \end{enumerate}
\end{lemma} 
\begin{proof}

(i). Assume $X_l(q)$ is untwisted. Write $q=p^e$ and $q_0=p^f.$ Then \cite[Proposition 5.4.6.(i)]{kleidmanliebeck} gives $f\vert e$ and if $k=\frac{e}{f},$  $V=V(\lambda)\otimes V(\lambda)^{q_0}\otimes \dots\otimes V(\lambda)^{q_0^{k-1}}$ as required. \par 
    For $X_l(q)$ twisted and $V\cong V^{\tau_0},$ the same reasoning proves the result using \cite[Proposition 5.4.6.(ii)(a)]{kleidmanliebeck}. \par 
(ii). In each case we want to prove the inclusion $$^s\!X_l(q)\leq X_l(q^s)\leq GL_n(q_0),$$ where $s=2$ for $\,^2\!A_l(q),$ $\,^2\!D_l(q),$ $\,^2\!E_6(q)$ and $s=3$ for $\,^3\!D_4(q).$ \par
    Part (ii) of \cite[Proposition 5.4.6]{kleidmanliebeck} and \cite[5.4.7(a)]{kleidmanliebeck} gives $q=p^e,$ $q_0=p^f$ and $f\vert se.$ Then if $k=\frac{se}{f},$ we have $V=V(\lambda)\otimes V(\lambda)^{q_0}\otimes \dots\otimes V(\lambda)^{q_0^{k-1}}$ for some $\lambda.$  Part (\ref{qnode1}) says that as an $X_l(q^s)$-module, $V$ is also realized over $\mathbb{F}_{q_0}.$  \par 

 (iii). In each case we want to prove the inclusion $$^2\!X_l(q)\leq X_l(q)\leq GL_n(q_0).$$ Remark \cite[5.4.7(b)]{kleidmanliebeck} gives $q=p^e,$ $q_0=p^f$ and $f\vert e.$  Then if $k=\frac{e}{f},$ $V=V(\lambda)\otimes V(\lambda)^{q_0}\otimes \dots\otimes V(\lambda)^{q_0^{k-1}}$ for some $\lambda.$ Part (\ref{qnode1}) says that as an $X_l(q)$-module, $V$ is also realized over $\mathbb{F}_{q_0}.$ 

\end{proof}

We will repeatedly use the following results from \cite{alvaro} and \cite{lubeck}. 

\begin{theorem}\textup{\cite[Thm 1.1 and 1.2]{alvaro}}\label{alvarobounds}
   Let $L=X_l(q)$ be a finite quasisimple group of Lie type, $\lambda$ a $p$-restricted weight and $V(\lambda)$ an irreducible module for $L$.  Assume $l\geq K$ and $\dim V(\lambda)<N$ where $K$ and $N$ are in Table \ref{tab:alvarodims}. Let $\epsilon_p(k)$ be 1 if $p\vert k$ and zero otherwise. Then $\lambda$ and  $\dim (V(\lambda))$ are as in Table \ref{tab:alvarodims}.
   \end{theorem}
   \begin{table}[]
       \centering
       \[\begin{array}{|c|c|c|c|c|}
       \hline
      L& A_l(q) & B_l(q)&C_l(q)&D_l(q) \\\hline
      K  & 9&14&14&16\\\hline
      N&\binom{l+1}{4}& 16\binom{l}{4} &16\binom{l}{4} &16\binom{l}{4} \\\hline
   \end{array}\] 
   \[\begin{array}{|c|c|c|}
   \hline 
      L& \lambda &  dim V(\lambda)\\ \hline
    A_l(q)&   \omega_1 & l+1\\
      &  \omega_2 & {l+1 \choose 2} \\
     &  2\omega_1 & {l+2 \choose 2}\\
     &  \omega_1+\omega_l & (l+1)^2-1-\epsilon_p(l+1)\\
     &  \omega_3 & {l+1 \choose 3}\\
     &  3\omega_1 & {l+3 \choose 3}\\
     &  \omega_1+\omega_2 & 2{l+2 \choose 3}-\epsilon_p(3){l+1 \choose 3}\\
      & \omega_1+\omega_{l-1} & 3{l+2 \choose 3}-{l+2 \choose 2}-\epsilon_p(l)(l+1)\\
      &  2\omega_1+\omega_{l} & 3{l+2 \choose 3}+{l+1 \choose 2}-\epsilon_p(l+2)(l+1)\\ \hline
          C_l(q)& \omega_1 & 2l\\
      & \omega_2 & 2l^2-l-1-\epsilon_p(l)\\
      & 2\omega_1 & {2l+1 \choose 2}\\
      & \omega_3 & {2l \choose 3}-2l-\epsilon_p(l-1)(2l)\\
      & 3\omega_1 & {2l+2 \choose 3}\\
     &  \omega_1+\omega_2 & 16{l+1 \choose 3}-\epsilon_p(2l+1)(1-\epsilon_p(3)(2l))-\epsilon_p(3)({2l \choose 3}-2l)\\ \hline 
        B_l(q)&\omega_1 & 2l+1\\
      & \omega_2 & {2l+1 \choose 2}\\
      & 2\omega_1 & {2l+2 \choose 2}-\epsilon_p(2l+1)\\
      & \omega_3 & {2l+1 \choose 3}\\
      & 3\omega_1 & {2l+3 \choose 3}-(2l+1)-\epsilon_p(2l+3)(2l+1)\\
     &  \omega_1+\omega_2 & 16{l+\frac{3}{2} \choose 3}-\epsilon_p(l)(2l+1)-\epsilon_p(3)({2l+1 \choose 3})\\ \hline 
      D_l(q)& \omega_1 & 2l\\
     &  \omega_2 & {2l \choose 2}-\epsilon_p(2)(1+\epsilon_p(l))\\
      & 2\omega_1 & {2l+1 \choose 2}-1-\epsilon_p(l)\\
      & \omega_3 & {2l \choose 3}-\epsilon_p(l+1)(2l)\\
     &  3\omega_1 & {2l+2 \choose 3}-2l-\epsilon_p(l+1)(2l)\\
     &  \omega_1+\omega_2 & 16{l+1 \choose 3}-\epsilon_p(2l-1)(2l)-\epsilon_p(3)({2l \choose 3})\\\hline
   \end{array}\]
    \caption{Nonzero $p$-restricted dominant weights $\lambda$ such that $dim V(\lambda) \leq N$ and $l\geq K.$}
    \label{tab:alvarodims}
\end{table}

The following is an analogous result for small values of $l$ using the results in \cite{lubeck}.

\begin{theorem}\label{lowerlprestrict}
\begin{enumerate}[(i)]
    \item Let $L=A_l(q)$ or $^2A_l(q)$ and $V=V(\lambda)$ where $\lambda$ is a $p$-restricted weight of $L.$ The smallest three possible dimensions of $V$ are for $\lambda=\omega_1,$ $\omega_2$ and $2\omega_1.$  Furthermore, for $l\leq 8$ one of the following holds.
    \begin{enumerate}
        \item $\lambda$ is as in Table \ref{tab:alvarodims}.
        \item $\lambda=\omega_4$ and and $(l,\dim(V(\lambda)))\in \{(7,70), (8,126)\}.$ 
        \item $\dim (V(\lambda))\geq N_A,$ where $N_A$ is as in the following table.
        \[\begin{array}{c|c|c|c|c|c|c|c}
         l & 2  &3&4&5&6&7&8  \\\hline
          N_A   & 14&19&45&90&147&112&156\\\hline
        \end{array}\]
    \end{enumerate}
    Moreover, if $\tau_0(\lambda)=\lambda,$ one of the following holds. 
    \begin{enumerate}
        \item[(d)] $(\lambda,\dim V(\lambda))=(\omega_1+\omega_l,(l+1)^2-\epsilon_p(l+1))$ 
        \item[(e)] $(l,\lambda, \dim(V(\lambda))) \in$ $\{(3, 2\omega_2, 19), (5, \omega_3, 20), (7, \omega_4, 70) \}$ 
       \item[(f)] $\dim (V(\lambda))\geq N_{A'},$ where $N_{A'}$ is as in the following table.
        \[\begin{array}{c|c|c|c|c|c|c|c}
         l & 2  &3&4&5&6&7&8  \\\hline
          N_{A'}   & 19&44&74&154&344&657&1135\\\hline
        \end{array}\]
    \end{enumerate}
    \item Let $L=B_l(q)$ with $q$ odd, and $V=V(\lambda)$ where $\lambda$ is a $p$-restricted weight of $L.$  For $3\leq l\leq 13$ one of the following holds.
    \begin{enumerate}
        \item $\lambda$ is as in Table \ref{tab:alvarodims}.
        \item  $(\lambda, \dim(V(\lambda)))=(\omega_l, 2^l)$ 
        \item $\dim (V(\lambda))\geq N_B,$ where $N_B$ is as in the following table. 
        \begin{center}
\begin{tabular}{c|c|c|c|c|c|c|c|c|c|c|c}
    $l$& $3$&$4$&$5$&$6$&$7$&$8$&$9$&$10$&$11$&$12$&$13$\\\hline
       $N_B$  & $27$&$64$&$100$&$208$&$128$&$256$&$512$&$1000$&$1331$&$1728$&$2197$ \\\hline

    \end{tabular}
    \end{center}
    \end{enumerate}
   \item Let $L=C_l(q)$ and $V=V(\lambda)$ where $\lambda$ is a $p$-restricted weight of $L.$ For $2\leq l\leq 13$ one of the following holds.
    \begin{enumerate}
        \item $\lambda$ is as in Table \ref{tab:alvarodims}.
        \item $l,$ $\lambda$ and $\dim(V(\lambda))$ are as in the following table. \[\begin{array}{c|c|c|c}
         l & \lambda  & \dim(V(\lambda)) &\text{extra conditions} \\\hline
          \text{all}   & \omega_l&2^l&q\,\,\text{even}\\\hline
          4   & \omega_3&\geq 40&\\\hline
          3   & \omega_3&\geq 13&\\\hline
          2   & 2\omega_2& 10&\\\hline
        \end{array}\]
        \item $\dim (V(\lambda))\geq N_C,$ where $N_C$ is as in the following table. 
        \begin{center}
            
        \begin{tabular}{c|c|c|c|c|c|c|c|c|c|c|c|c}
    $l$& $2$&$3$&$4$&$5$&$6$&$7$&$8$&$9$&$10$&$11$&$12$&$13$\\ \hline
       $N_C$  & $11$&$25$&$64$&$100$&$208$&$128$&$256$&$512$&$1000$&$1331$&$1728$&$2197$ \\\hline
      
    \end{tabular}
    
        \end{center}
    \end{enumerate} 

\item Let $L=D_l(q)$ or $^2D_l(q)$ and $V=V(\lambda)$ where $\lambda$ is a $p$-restricted weight of $L.$ For $4\leq l\leq 15$ one of the following holds.
    \begin{enumerate}
        \item $\lambda$ is as in Table \ref{tab:alvarodims}.
        \item $l,$ $\lambda$ and $\dim(V(\lambda))$ are as in the following table. \[\begin{array}{c|c|c}
         l & \lambda  & \dim(V(\lambda))  \\\hline
          \text{all}   & \omega_l&2^l\\\hline
          4   & \omega_1+\omega_3&48\\\hline
          5   & \omega_3&100\\\hline
          6   & \omega_3&208\\\hline
          7   & \omega_3&336\\\hline
        \end{array}\]
        \item $\dim (V(\lambda))\geq l^3.$
        \end{enumerate}
        \item Let $L=\,\,^\epsilon X_l(q)$ be an exceptional quasisimple group of Lie type and $V=V(\lambda)$ where $\lambda$ is a $p$-restricted weight of $L.$ Then one of the following holds. 
        \begin{enumerate}
            \item $L,$ $\lambda$ and $\dim(V(\lambda))$ are as in the following table. \[\begin{array}{c|c|c}
         L & \lambda  & \dim(V(\lambda))  \\\hline
        E_8(q)   & \omega_8&248\\\hline
      E_7(q)   & \omega_7&56\\\hline
       E_7(q)   & \omega_1&133-\epsilon_p(2)\\\hline
        ^\epsilon E_6(q)   & \omega_6&27\\\hline
        ^\epsilon E_6(q)   & \omega_3&78-\epsilon_p(3)\\\hline
         ^\epsilon F_4(q)   & \omega_4&26-\epsilon_p(3)\\\hline
         ^\epsilon F_4(q)   & \omega_1&52\\\hline
         ^\epsilon G_2(q)   & \omega_2&7-\epsilon_p(2)\\\hline
          ^\epsilon G_2(q)   & \omega_1&14\\\hline
          ^3D_4(q)&\omega_1&8\\\hline
          ^3D_4(q)   & \omega_2&28-2\epsilon_p(2)\\\hline
        \end{array}\]
         \item $\dim (V(\lambda))\geq N_E,$ where $N_E$ is as in the following table. 
         \[\begin{array}{c|c|c|c|c|c|c|c}
             G & E_8(q)&E_7(q)& E_6(q)& ^2E_6(q)& ^\epsilon  F_4(q)&^\epsilon G_2(q)&^3D_4(q)\\\hline
              N_E&3626&856&324&572&196&26&195 \\\hline
         \end{array}\]
        \end{enumerate}
\end{enumerate}
\end{theorem}
\begin{proof}
These are all from tables in \cite{lubeck}.
\end{proof}

\subsection{Primitive Affine Groups}
In this section we include some preliminary results that we will use in proving our theorems. Recall that a primitive affine group is of the form $G = VG_0\leq A\Gamma L(V),$  where $V=V_n(q_0)$ is a finite dimensional vector space over a finite field $\mathbb{F}_{q_0},$ the stabilizer of $0$ is $G_0 \leq \Gamma L(V ),$  acting irreducibly on $V$ and $V$ acts by translations. \par
Firstly, let us give an expression for the diameter of an orbital graph of a primitive affine group. Let us denote the distance between two vertices in the graph $a,b\in V$ by $d(a,b).$ 
\begin{lemma}\label{diamfact}
Let $G=VG_0$ be a primitive group of affine type, $a\in V\setminus{0}$ and $O=\{0,a\}^G$ be an orbital. Then in the corresponding orbital graph the following holds for all $b\in V:$  $$d(0,b)=\min (k: b \,\,\mbox{can be expressed as a sum of}\,\,k\,\,\mbox{elements in}\,\,\pm a^{G_0}).$$ 
\end{lemma}
\begin{proof}
We can show this by induction on the distance from $0$. The base case holds as by definition $0$ is joined to $b$ if and only if $b\in \pm a^{G_0}.$ \par
The elements of distance $m-1$ can be expressed as a sum of minimum $m-1$ elements in $\pm a^{G_0}.$ 
As $0$ is adjacent to every element in $\pm a^{G_0}$, the neighbours of $l\in V$ are of the form $l\pm a^{G_0}.$ If $d(0,l)=m-1$ and $l'$ is a neighbour of $l$ such that $d(0,l')\geq m-1,$ then $d(0,l')=m$ and $l'$ can can be expressed as a sum of minimum $m$ elements in $\pm a^{G_0}.$  
\end{proof}

Using this result we obtain the following upper bound for the orbital diameter. 

\begin{lemma}\textup{\cite[Lemma 3.1]{orbdiam}}\label{upperbound}
      Let $G = V_n(q_0).G_0 $ and assume that
$G_0$ contains the scalar matrices of $GL_n(q_0).$ Then $orbdiam(G,V)\leq n.$ 
  \end{lemma}  
  \begin{proof}
    Let $\{0,u\}^{G}$ be an orbital. Now as $G_0$ acts irreducibly, $u^{G_0}$ contains a basis $u_1,\dots,u_n$ of $V_n(q_0).$ Also $ku\in u^{G_0}$ for all $k\in \mathbb{F}_{q_0}^\star$ by assumption, so we have a path of length $n,$ $$0\frac{\hspace{0.2cm}\hspace{0.2cm}}{}k_1u_1\frac{\hspace{0.2cm}\hspace{0.2cm}}{}\cdots\frac{\hspace{0.2cm}\hspace{0.2cm}}{} k_1u_1+\dots+k_nu_n$$ where the $k_i$ are arbitrary scalars.
  \end{proof}

The next two results are clear.

\begin{lemma}\label{subgroup}
  Let $H_0\leq G_0\leq \Gamma L(V).$ Then $orbdiam(VG_0,V)\leq orbdiam(VH_0,V).$  
\end{lemma}

\begin{lemma}\label{rankbound}
 Let $G$ be a primitive group acting on a set $\Omega$ with permutation rank $r.$ Then $orbdiam(G,\Omega)\leq r-1.$   
\end{lemma}

Next we include a complete classification of primitive affine groups with orbital diameter 1. 
Clearly $G=V_n(q)G_0$ has orbital diameter 1 if and only if $G$ is 2-homogeneous. The 2-homogeneous affine permutation groups that are not 2-transitive have been classified in \cite{kantor} and the 2-transitive affine groups were classified in \cite{hering} and \cite{liebeckaffine}, as described in the following theorem.

\begin{theorem}\textup{\cite[Appendix 1]{liebeckaffine}\cite[Thm 1]{kantor}}\label{diam1grps} 
Let $G=VG_0$ with $V\cong (\mathbb{F}_p)^d$ be an affine permutation group with orbital diameter 1. Then one of the following holds. 
\begin{enumerate}
    \item [{\rm (i)}] $(G,V)$ is $2$-transitive, listed in \cite{liebeckaffine}[Appendix 1]
    \item [{\rm (ii)}] $G \leq A \Gamma L_1(q)$  with $q\equiv 3 \,\,(\!\!\!\!\mod 4)$ and $G$ is 2-homogeneous but not 2-transitive.
\end{enumerate}
\end{theorem}

We will also need the following result from \cite{orbdiam}.

\begin{lemma}\label{lemma2.1}
Let $G = V_n(q_0).G_0 $ be a primitive affine group and $V=V_n(q)$. Let $orbdiam(G,V)=d$ and let $\mathcal{O}$ be an orbit of $G_0$ on $V \setminus \{0\}$. Then 
\begin{enumerate}[(i)]
    \item \label{part1} The following inequality holds:$${q_0}^n\leq 1+\sum_{i=1}^d((2,q_0-1)\vert \mathcal{O}\vert)^i.$$
    \item \label{boundformula1}If $\mathcal{O}=-\mathcal{O}$, then $${q_0}^n\leq 1+\sum_{i=1}^d(\vert \mathcal{O}\vert)^i\leq 2(\vert \mathcal{O}\vert)^d. $$
    \item \label{boundeq1} If $\mathcal{O}=-\mathcal{O}$ and $\vert \mathcal{O}\vert \leq q_0^{r},$ then $$d \geq \frac{n-\log_{q_0}(2)}{r}.$$
    \item\label{boundeq2} If $\mathcal{O}=-\mathcal{O}$ and $\vert \mathcal{O}\vert \leq \frac{q_0^{r}}{2},$ then $$d \geq \frac{n}{r}.$$
\end{enumerate}

\end{lemma}

\begin{proof}
    Part $(i)$ is \cite{orbdiam}[Lemma 2.1]. Since we are considering undirected graphs, $0$ is adjacent to $\pm \mathcal{O}$ and so the number of vertices at distance $k$ is at most $1+\sum_{i=1}^k(2\vert \mathcal{O}\vert)^i.$ If $2\vert q,$ then $\pm \mathcal{O}=\mathcal{O}$ so we can omit the multiplication by 2. In part $(ii)$ the orbitals corresponding to $Y$ are self-paired so again we can omit the multiplication by $(2,q-1).$ We obtain parts $(iii)$ and $(iv)$ by substituting the bounds on the orbit sizes into $(ii).$
\end{proof}
Recall Hypothesis \ref{hypothesisgen} from the Introduction.

\begin{lemma}\label{autsorbit}
  Let $G$  be as in Hypothesis \ref{hypothesisgen}. Let $\mathcal{O}$ be an orbit of $G_0$ on $V.$ Then $$\vert \mathcal{O} \vert \leq (q_0-1)\vert Aut(G_s)\vert.$$   
\end{lemma}

\begin{proof}
    Let $Z=\mathbb{F}_{q_0}^\star I_n.$ 
   Since $C_{PGL_n(q_0)}(G_s)=1$ by \cite{kleidmanliebeck}[Lemma 4.0.5] 
    we have that $G_s \cong \frac{G_0^\infty Z}{Z}$ and 
    $\frac{G_0 Z}{Z}\leq Aut(G_s), $ so the bound follows.
\end{proof}

\begin{lemma}\label{eq3}
Let $G$ be as in Hypothesis \ref{hypothesisgen} and let $d=orbdiam(G,V)$. Then $$ n \leq 1+d\log_2(\vert Aut(G_s)\vert).$$ 
\end{lemma}

\begin{proof}
Call $k=\vert Aut(G_s)\vert .$ Then $\vert G_0\vert \leq (q_0-1)k$ by Lemma \ref{autsorbit}. Hence Lemma \ref{lemma2.1}(\ref{boundformula1}) tells us that $1 + (q_0-1)k + \dots + ((q_0-1)k)^d \geq q_0^n.$ \par Using the fact that $1+x+x^2+\dots+x^d\leq 2x^d$ for $x\geq 2,$ if follows that 
$q_0^n\leq 2((q_0-1)k )^d$ which is equivalent to 
$$n\leq \log_{q_0}(2((q_0-1)k )^d).$$ We will now show that this is bounded above by $1+d\log_2(k) $ as required.\par  We have two claims to prove. \par
\textit{\underline{Claim 1}} For $q_0\geq 3$ the value of $(\log_{q_0}(2(2(q_0-1)k )^d)$ is maximal when $q_0=3.$\par
We want to show that $$\log_{q_0}(2(2(q_0-1)k )^d)\leq \log_3(2(4k )^d)$$ for $q_0\geq 3,$
which is equivalent to 
$$\log_{q_0}(2)+d\log_{q_0}(2(q_0-1))+d\log_{q_0}(k)\leq\log_3(2)+d\log_3(4)+d\log_3(k).$$ 
\par Since for $q_0\geq3$ it is clear that $\log_{q_0}(2)\leq \log_{3}(2),$ it suffices to show that $$\log_{q_0}(2(q_0-1)k)\leq \log_3(4k).$$

To solve inequalities with one unknown, we use Wolfram-Alpha \cite{Wolfram|Alpha}.
As $k=\vert Aut(G_s)\vert\geq 60,$ Wolfram-Alpha shows that this inequality holds for any $4\leq q_0\leq 9.$  \par 
For $q_0\geq 10,$ by Wolfram-Alpha, $\log_{q_0}(2(q_0-1))\leq\log_3(4),$ so as  $\log_{q_0}(k)\leq \log_{3}(k),$ Claim 1 holds.

\par 
\textit{\underline{Claim 2}} $\log_3(2(4k)^d)\leq 1+d\log_2(k)$ \par 
This is equivalent to $$\log_3(2)+d\log_3(4k)\leq1+d\log_2(k). $$ Since $k\geq 60,$ $\log_3(4k)\leq\log_2(k)$ and $\log_3(2) <1,$ Claim 2 follows. \par
In Claim 1 we showed that for $q_0\geq 3,$ the value $\log_{q_0}(2(2(q_0-1)k )^d),$ which is an upper bound for $\log_{q_0}(2((q_0-1)k )^d),$ is maximal for $q_0=3.$ In Claim 2 we showed that $\log_{q_0}(2(2(q_0-1)k )^d)$ with $q_0=3$ is less than $1+d\log_2(k).$ 
In particular we showed
$$n\leq \log_{q_0}(2((q_0-1)k )^d)\leq \log_{q_0}(2(2(q_0-1)k )^d) \leq \log_3(2(4k)^d)\leq 1+d\log_2(k)$$ where the third inequality is Claim 1 and the last in Claim 2. So the result follows. \par 
\end{proof}
A similar proof to Claims 1 and 2 shows that the relation below also holds. 

\begin{lemma}\label{log2.2}
Let $q_0$ be a prime power and $k\geq 60.$
Then $\log_{q_0}(1+(q_0-1)k+((q_0-1)k)^2)$ is maximal for $q_0=2.$
\end{lemma}

\subsection{Use of Computation}
We use computation to compute or bound the orbital diameter for affine groups satisfying Hypothesis \ref{hypothesisgen} for various specific simple $G_s$ in specific representations. Matrix generators for such groups can be constructed using GAP \cite{GAP4}, Magma \cite{magma}, the AtlasRep \cite{atlasrep} Package or the online ATLAS, http://groupatlas.org/Atlas/v3/index.html. Exact orbital diameters can be calculated in many cases using the Grape \cite{grape} Package. 
\section{ Lie Type Stabilizer In Defining Characteristic}

In this section we prove Theorems \ref{lhalf} and \ref{diam2}. The groups considered here all satisfy the following hypothesis.

\begin{hypothesis}\label{hypothesis}
  Let $G=VG_0$ be a primitive affine group such that $G_0^\infty/Z(G_0\infty)=X_l(q),$ where $X_l(q)$ is a finite simple group of Lie type in characteristic $p.$ Suppose that $V$ is an absolutely irreducible $\mathbb{F}_{q_0}G_0^\infty-$module in characteristic $p$ of dimension $n.$ Also assume that $V$ cannot be realised over a proper subfield of $\mathbb{F}_{q_0}.$ Let $orbdiam(G,V)=d.$
\end{hypothesis}

We begin with a natural way to estimate the size of the orbit of a maximal vector under $G_0$ as defined in \cite[Def 15.11]{malletesterman}. We will denote this orbit in our proofs by $\mathcal{O}.$\par 
Recall that for a dominant weight $\lambda,$ denote the parabolic subgroup stabilizing a maximal $1$-space in $V(\lambda)$ by $P_\lambda^F$ as defined above Lemma \ref{para}. For simplicity we will abuse notation and denote this by $P_\lambda.$ The next result is clear.
 
\begin{lemma}\label{orbit-stabparaborel}
   Let $G_0\leq \Gamma L(V)$ be as in Hypothesis \ref{hypothesis} and $V=V(\lambda)$. Let $P_\lambda$ be the parabolic fixing a maximal $1$-space, $\langle v^+ \rangle$, $B\leq P_\lambda$ a Borel and $\mathcal{O}={v^+}^{G_0}$. Then   $$\vert \mathcal{O} \vert \leq (q_0-1)\vert G_0\colon P_\lambda\vert \leq(q_0-1)\vert G_0\colon B\vert.$$ 
\end{lemma}

We give an example of finding such an upper bound. We will repeatedly use this method in our proofs.

\begin{example}\label{paralambdaexample}\rm 
Consider the case when $\frac{G_0^\infty}{Z(G_0^\infty)}=A_l(q)$ and $\lambda=\omega_1+\omega_l.$ Lemma \ref{para} tells us that the parabolic $P_{\lambda}$ fixes a $1$-space in $V$. By definition, $P_{\lambda}=P_{1,l}.$ Then $\vert G_0\colon P_{1,l}\vert = \frac{(q^l-1)(q^{l+1}-1)}{(q-1)^2}$ and Lemma  \ref{orbit-stabparaborel} gives $\vert \mathcal{O} \vert \leq (q-1)\frac{(q^l-1)(q^{l+1}-1)}{(q-1)^2}\leq q^{2l+1}-1.$ 
\end{example}

We now provide a lemma concerning the examples with orbital diameter at most $2.$

\begin{lemma} \label{smalldiamlemma}
     Let $G$ be as in Hypothesis \ref{hypothesis}. Assume that if $G_0^\infty$ is a classical group then $V$ is not a natural module for $G_0^\infty.$ If $G$ is as in Table \ref{tab:defchar} and contains the scalars $\mathbb{F}_q^\star$, then the orbital diameters and ranks are as in Table \ref{tab:defchar}.
\end{lemma}
\begin{proof}
   First recall from Lemma \ref{rankbound} that the orbital diameter is bounded above by $r-1$, where $r$ is the permutation rank. It follows from the proof of Lemma \ref{rankbound} that $orbdiam(G,V)=r-1$ if and only if $G$ has a distance-transitive orbital graph.  \par 
    Consider $G_0 \triangleright B_4(q)$ with $\lambda=\omega_4.$ By \cite[Lemmas 2.9, 2.11]{glms}, $G$ is a rank 4 group, and by \cite[Thm 1.1]{disttransclassical} it has no distance-transitive orbital graphs. Hence it has orbital diameter 2. \par 
    Consider $G_0 \triangleright G_2(q)$ with $\lambda=\omega_1$ with $q$ odd.  By \cite[page 498]{liebeckaffine},  $G$ is a rank 4 group with no distance-transitive orbital graphs by \cite[Thm 1.1]{distancetransitive} so it has orbital diameter 2. \par 
     For the remaining cases in Table \ref{tab:defchar}, by \cite{liebeckaffine}, $G$ has rank $2$ or $3,$ so the orbital diameter is $1$ or $2,$ respectively.  \par 
\end{proof}
We will also use the following two lemmas in our proofs, for which we thank Aluna Rizzoli.

\begin{lemma}\label{fullrank}
Let $\overline{G}$ be a simple algebraic group over an algebraically closed field. Let $P\leq \overline{G}$ be a parabolic subgroup. Then for all $g \in \overline{G},$ $P\cap P^g$ contains a maximal torus.
\end{lemma}
\begin{proof}
Let $T\leq P$ be a maximal torus, $W\cong \frac{N(T)}{T}$ be the Weyl group of $\overline{G}$ and $B$ be a Borel subgroup such that $T\leq B\leq P.$ Then by the Bruhat decomposition, $\overline{G}=\bigcup_{w\in W}Bn_wB,$ so $g=b_1n_wb_2$ with $b_1,b_2\in B$ and $n_w$ a preimage of $w$ in $N_{\overline{G}}(T).$ Then \begin{align*}
    P\cap P^g&=P\cap P^{b_1n_wb_2}\\
    &=P\cap P^{n_wb_2}\\
    &=(P^{b_2^{-1}}\cap P^{n_w})^{b_2}\\
    &=(P\cap P^{n_w})^{b_2}.\\
    \end{align*}
As $T\leq P\cap P^{n_w}, $  $T^{b_2}\leq (P\cap P^{n_w})^{b_2}$ so the intersection of two conjugates of $P$ contains a maximal torus.   
\end{proof}
The next result specifies some possible intersections of parabolics. We will use the notation $Q_k$ for a connected unipotent group of dimension $k$ and $T_i$ for a torus of rank $i$. 

\begin{lemma} \label{alunalemma}
\begin{enumerate}
    \item Let $k=\overline{\mathbb{F}_3}$ and $\overline{G}=C_3(k).$  The possible intersections of two conjugates of the parabolic $P_2$ of $\overline{G}$ are $$P_2,\,\,Q_8T_3,\,\,Q_5A_1T_1,\,\,Q_5T_3,\,\,A_1A_1T_1.$$
    \item Let $k=\overline{\mathbb{F}_p}$ and $\overline{G}=E_7(k).$ The possible intersections of two conjugates of the parabolic $P_7$ of $\overline{G}$ are $$P_7,\,\,Q_{42}D_5T_2,\,\,Q_{33}D_5T_2,\,\,E_6T_1.$$
\end{enumerate}
\end{lemma}
\begin{proof}
Let $\overline{G}$ be a simple algebraic group over an algebraically closed field $\overline{\mathbb{F}_p}$. Let $P=P_J\leq \overline{G}$ be a parabolic subgroup, $T\leq P$ a maximal torus, and $W\cong \frac{N(T)}{T}$ the Weyl group of $\overline{G}$ with respect to $T$. Then for $g \in \overline{G},$ the different conjugacy classes for $P \cap P^g$ correspond to distinct $(P,P)$-double cosets in $\overline{G},$ as described in \cite[Section 2.8]{carter2}.
By the Bruhat decomposition we have $G =  \sqcup (Pn_{w_i}P)$ with $w_i\in W$ and $n_{w_i}$ a preimage of $w_i$ in $N_{\overline{G}}(T).$ which are the representatives of the $(W_J,W_J)$ double cosets in $W.$ Therefore all the distinct intersections $P\cap P^g$ are given by $P\cap P^{n_{w_i}}$. For (1) and (2), we can obtain these $w_i$'s by performing computations involving intersections of double cosets in the Weyl groups $W(C_3)=2^3.S_3$ and $W(E_7)=2\times Sp_6(2)$ using GAP \cite{GAP4} or Magma\cite{magma}. By \cite[Thm 2.8.7]{carter2}, the intersections of two conjugates of $P$ are generated by the maximal torus together with all root subgroups contained in the intersection. By considering the action of the $w_i$s on the root system, using computations in GAP and Magma we get the lists of possible intersections as in the statement of the lemma. 
\end{proof}
\subsection{Classical Stabilizers} In this section we prove Theorems \ref{lhalf}  and \ref{diam2} for the case when $G_0$ is a classical group.

We begin with a result on the diameter when $G=VG_0$ and $V$ is the natural module of $G_0.$

\begin{lemma}\label{naturalmod}
    Let $G=V_n(q)G_0$ with $G_0\leq \Gamma L_n(q)$ a classical group and $V_n(q)$ the natural module of $G_0.$ \begin{enumerate}
        \item If $G_0 \triangleright SL_n(q)$ then $orbdiam(G,V)=1.$
        \item If $G_0 \triangleright Sp_n(q)$ then $orbdiam(G,V)=1.$
        \item If $G_0 \triangleright SU_n(q^{1/2})$ with  $n\geq 3,$ then $orbdiam(G,V)=2.$
        \item If $G_0 \triangleright \Omega^\epsilon_n(q)$ with  $n\geq 4,$ or $n=3$ and $q\not \equiv 1 \mod 4,$ then $orbdiam(G,V)=2.$
        \item If $G_0 \triangleright \mathbb{F}_q^\star. \Omega^\epsilon_3(q)$ with $q\equiv 1 \mod 4,$ then $orbdiam(G,V)=2.$
        \end{enumerate}
\end{lemma}
\begin{proof}
We prove the statements in turn. \par 
1. \& 2. Here $G_0$ acts transitively on $V_n(q)$ and the result is clear. \par 
3. The orbits of $SU_n(q^{1/2})$ on $V_n(q)$ are of the form $O_{\lambda}=\{v\in V\setminus{0}: B(v,v)=\lambda\}$ where $ \lambda\in \mathbb{F}_{q^{1/2}}$ and $B$ is the associated Hermitian form \cite[Lemma 2.10.5]{kleidmanliebeck}. Therefore, by Lemma \ref{diamfact} to show that $orbdiam(G,V_n(q))=2$ we have to prove that we can express all vectors as a sum of two vectors of a given norm. Let $v,w\in O_{\lambda}.$ Then $B(v+w,v+w)=2\lambda+B(v,w)+\overline{B(v,w)},$ and we want to show that this can be arbitrary.
It is sufficient to prove that for all $\sigma \in \mathbb{F}_{q}$ there is $v,w\in O_\lambda$ such that $B(v,w)=\sigma,$ because $2\lambda+\sigma+\overline{\sigma}$ is arbitrary in $\mathbb{F}_{q^{1/2}}.$ \par
    Recall the standard basis of $V_n(q),$ $\{e_1,\dots,e_k,f_1,\dots,f_k,x\}$ for $n=2k+1$ and $\{e_1,\dots,e_k,f_1,\dots,f_k\}$ for $n=2k,$ where for all $i,j$ we have $B(e_i,e_j)=B(f_i,f_j)=0,$ $B(e_i,f_j)=\delta_{i,j}$ and $B(e_i,x)=B(f_i,x)=0,$ $B(x,x)=1$\cite[Prop 2.2.2]{giudici}.\par First assume that $\lambda=0.$ Then choose $v=e_1$ and $w=\sigma f_1$ which gives $B(v,w)=\sigma$ and we are done. \par 
    Now assume $\lambda\neq 0.$ Since the trace map $\mathbb{F}_q\to \mathbb{F}_{q^{1/2}} $ sending  $\alpha\to \alpha+\overline{\alpha}$ is surjective, there is $\mu\in \mathbb{F}_{q}$ such that $\mu+\overline{\mu}=\lambda.$ For $n\geq 4$ put $v=e_1+\mu f_1.$ Let $w=\sigma f_1+e_2+\mu f_2$ and so $B(v,w)=\sigma$ and we are done. \par
    Since the map $\mathbb{F}_q\to \mathbb{F}_{q^{1/2}} $ sending  $\alpha\to \alpha\overline{\alpha}$ is surjective, so
   there is $\chi\in \mathbb{F}_{q} $ such that $\chi \overline{\chi}=\lambda.$ For $n=3$ the pair $v=e_1+\mu f_1$ and $w=\sigma f_1+\chi x$ works.  \par 
4 and 5. The orbits of $\Omega^\epsilon_n(q)$ on $V_n(q)$ for $n\geq 4$ are of the form $O_{\lambda}=\{v\in V\setminus{0}: Q(v)=\lambda\}$ where $ \lambda\in \mathbb{F}_q,$ $Q$ is the associated quadratic and $B$ is the associated bilinear form \cite[Lemma 2.10.5]{kleidmanliebeck}. For $n=3$ the orbit $O_0$ splits into two orbits of size $\frac{q^2-1}{2}$ \cite[Lemma 2.10.5(iv)]{kleidmanliebeck}. For $q\equiv 3\mod 4,$ these are negatives of each other, so produce one undirected orbital graph, so we can regard them as one. For $q\equiv 1\mod 4,$ assuming all scalars are present, these also produce one undirected orbital graph.\par 
    Since $Q(v+w)=Q(v)+Q(w)+B(v,w)$ it is sufficient to show that for given $\sigma\in \mathbb{F}_q$ there is $v,w\in O_\lambda$ such that $B(v,w)=\sigma.$ This is achieved in a similar fashion to part (3).\par 

\end{proof}

\subsubsection{$G_0\triangleright A_l(q)$}
We continue with the case of the proof of Theorems \ref{lhalf} and \ref{diam2} when $G=VG_0$ and $\frac{G_0^\infty}{Z(G_0^\infty)}=A_l(q).$ Note that we are using Lie notation $A_l(q)$ for $PSL_{l+1}(q).$ Recall $d=orbdiam(G,V),$ $n=\dim(V)$ and $V=V(\lambda).$

\begin{theorem}\label{al(q)}
Let $G$ be as in Hypothesis \ref{hypothesis} with $X_l(q)\cong A_l(q).$ 
\begin{enumerate}[(i)]
    \item \label{ali1} If $\lambda$ (or $\lambda^\star)$ is  in Table \ref{ali}, then the value of $n$
 and a lower bound for $d$ is given.
\begin{table}[]
\[\begin{array}{c|c|c|c|c|c}
    &\lambda & n &d\geq&\text{extra conditions}  \\\hline
   1.& \omega_1 &l+1&=1&q=q_0\\
     2.&\omega_2 &\frac{l(l+1)}{2}&\lfloor\frac{l+1}{2}\rfloor &q=q_0\\
      3.&2\omega_1&\frac{(l+2)(l+1)}{2} & l+1 &q=q_0\,\,\&\,\, p>2\\
      4.&\omega_1+\omega_l& l^2+2l&l+1  & q=q_0\,\,\&\,\, p\not \vert l+1 \\
      5.&\omega_1+\omega_l& l^2+2l-1& \frac{l^2+2l-1}{2l+1}& q=q_0\,\,\&\,\,p\vert l+1  \\
        6.&\omega_3&{l+1 \choose 3 }&\frac{l^2}{18}& q=q_0\\
         7.&\omega_1+p^i\omega_1&(l+1)^2&l+1& q=q_0\,\,\&\,\,p^i\neq q^{\frac{1}{2}}\\
         8.&\omega_1+p^i\omega_l&(l+1)^2&\frac{(l+1)^2}{2l+1}& q=q_0\\ 
         9.&\omega_1+q_0\omega_1&(l+1)^2&(l+1)/2& q=q_0^2   \\\hline
\end{array}\]
\caption{Bounds in Theorem \ref{al(q)} part \ref{ali1}}
    \label{ali}
\end{table}
    \item \label{part2ali}
    If $l\geq 9,$ then for all $\lambda$ not in Table \ref{ali},
    $$
        d\geq \frac{l(l-1)(l-2)}{12(l+2)}.
    $$
   
For $l\leq 8$ and $\lambda$ not in Table \ref{ali}, a lower bound on $d$ is as follows.
    \begin{align}\label{smalll}
     \begin{array}{c|c|c|c|c|c|c|c|c|c|c}
            l &8&7&6&5&4&3&2&1\\\hline
            d\geq &4&4&6&4&3&3&3&3
        \end{array}
     \end{align}
 \end{enumerate}   

\end{theorem}
The lower bounds in the statement of Theorem \ref{al(q)} are greater than the lower bounds contained in Theorem \ref{lhalf}, so Theorem \ref{lhalf} holds for $X_l(q)=A_l(q)$. 

\begin{proof} 
\par 
\underline{\textbf{(i) Proof of the bounds in Table \ref{ali}}}
\par Recall that $G_0^\infty\cong A_l(q)$ and let $W$ denote the natural module of dimension $l+1$ over $\mathbb{F}_q.$ \par 
We consider the weights $\lambda$ in Table \ref{ali}. Note that $q=q_0$ except in the last entry by Lemma \ref{q0lemma} part \ref{qnode1}. \par 
1. Here $\lambda=\omega_1$ and $V(\lambda)=W,$ the natural module. By Lemma \ref{naturalmod}, $G_0$ acts transitively, so the orbital diameter is 1.\par 
2. Here $\lambda=\omega_2$ and $V(\lambda)=\bigwedge^2W,$ the alternating square of $W.$ 
Choose a basis of $W$, $\{v_1,\dots,v_n\},$ so that $\{v_i\otimes v_j\vert 1\leq i,j\leq n\}$ is a basis of $W\otimes W.$ Now we have a $G_0$-isomorphism $$\phi:W\otimes W\to M_n(q)$$ via $$\phi \colon x=\sum a_{i,j}(v_i\otimes v_j )\mapsto A,$$ where $[A]_{i,j}=a_{i,j}.$ For $g\in G_0,$ if $g\in GL(W)$ then the action of $g$ sends $A\to gAg^T$ and for a field automorphism $\sigma,$ $\sigma$ sends $A\to(a_{i,j}^{\sigma}).$ Furthermore, $x\in \bigwedge^2W$ if and only if $\phi(x)$ is skew-symmetric with zeroes on the diagonal in characterictic 2.  \par 
For $q$ odd we identify $V(\lambda)=\bigwedge^2W$ with the space of $(l+1)\times (l+1)$ skew-symmetric matrices, $$\{A\in M_{l+1}(q)\vert A^T=-A\}.$$ For $q$ even, we identify $V=\bigwedge^2W$ with the set of symmetric matrices with zeroes on the diagonal. 
Since the action of $G_0$ preserves the rank,  all elements in an orbit have the same rank. 
Furthermore we know that skew-symmetric matrices have even rank. Let $A,B\in V$ with rank $a$ and $b$, respectively. Then $rank(A+B)\leq a+b,$ so we need to add up at least $\lfloor\frac{l+1}{2}\rfloor$ rank 2 skew-symmetric matrices to get a skew-symmetric matrix of maximal rank. Hence $d\geq \lfloor\frac{l+1}{2}\rfloor.$ \par 
3. Here $\lambda=2\omega_1$ and $V(\lambda)=S^2W$ with $p\neq 2.$ For each $x\in S^2W,$ $\phi(x)$ is symmetric, so we identify $V(\lambda)$ with the space of symmetric matrices, $$\{A\in M_{l+1}(q)\vert A^T=A\}.$$ The action of $G_0$ is the same as on the skew-symmetric matrices. Again, the rank is preserved, and as all ranks are possible, our lower bound is $l+1.$\par 
4. and 5. Here $\lambda=\omega_1+\omega_l,$ $n=l^2+2l-1$ or $l^2+2l,$ and $V(\lambda)$ is the adjoint module. \par Suppose $p$ does not divide $l+1.$ Then the adjoint module can be identified with $V_{ad}=\{A\in M_{l+1}(q)\,\,\vert\,\, tr(A)=0\},$ the space of the traceless $(l+1)\times (l+1)$ matrices and $G_0$ acts by conjugation. Conjugation preserves the rank, so on orbits the rank is constant. To get a traceless matrix of rank $l+1,$ we need to add up at least $l+1$ elements of an orbit with rank 1 traceless matrices, so $d\geq l+1$ in this case.   \par Suppose $p$ divides $l+1.$ Then $n=l^2+2l-1,$ so to find the bound we can use Lemma \ref{para}, which tells us that the parabolic $P_{\lambda}=P_{1,l}$ fixes a $1$-space. As in Example \ref{paralambdaexample}, we see using
Lemma \ref{lemma2.1}(\ref{boundformula1}) that $$ d\ge\frac{l^2+2l-1}{2l+1}.$$ \par 
6. Here $\lambda=\omega_3$ and $n=\frac{l^3-l}{6}.$ 
 Lemma \ref{para} tells us that the parabolic $P_{\lambda}=P_3$ fixes a $1$-space in $V$. Then $\vert G_0\colon P_{3}\vert = \frac{(q^{l-1}-1)(q^l-1)(q^{l+1}-1)}{(q^3-1)(q^2-1)(q-1)}$ and Lemma  \ref{orbit-stabparaborel} gives the bound $\vert \mathcal{O} \vert \leq(q-1) \frac{(q^{l-1}-1)(q^l-1)(q^{l+1}-1)}{(q^3-1)(q^2-1)(q-1)}\leq \frac{q^{3l-3}}{2}.$ Now Lemma \ref{lemma2.1}(\ref{boundeq2}) gives us the result.\par 
7. Here $\lambda=\omega_1+p^i\omega_1$ and $n=(l+1)^2$ with $p^i\neq q^{\frac{1}{2}}.$   By Lemma \ref{para} the parabolic $P_1$ fixes a $1$-space. By Lemma  \ref{orbit-stabparaborel} there is an orbit $\mathcal{O}$ such that $\vert \mathcal{O} \vert \leq q^{l+1}-1.$ Hence Lemma \ref{lemma2.1}(\ref{boundformula1}) gives $d\geq l+1.$ \par 
8. Here $\lambda=\omega_1+p^i\omega_l$ and $n=(l+1)^2.$  The parabolic $P_{1,l}$ fixes a $1$-space.  By Lemma  \ref{orbit-stabparaborel} there is an orbit $\mathcal{O}$ such that $\vert \mathcal{O} \vert \leq q^{2l+1}-1.$ Hence, Lemma \ref{lemma2.1}(\ref{boundformula1}) gives $d\geq \frac{(l+1)^2}{2l+1}.$  \par
9. Here $\lambda=\omega_1+q_0\omega_1,$ $n=(l+1)^2$ and $V=V_n(q_0)$ where $q=q_0^2.$ Now Lemma \ref{q0lemma} part \ref{subfield} holds. By Lemma \ref{para} the parabolic $P_1$ fixes a $1$-space. By Lemma  \ref{orbit-stabparaborel} there is an orbit $\mathcal{O}$ such that $\vert \mathcal{O} \vert \leq \frac{ (q_0)^{2l+2}}{2},$ so the bound follows from Lemma \ref{lemma2.1}(\ref{boundeq2}). \par 
From now on assume that part (\ref{ali1}) does not hold, i.e. $\lambda$ is not in Table \ref{ali}. \par 
\underline{\textbf{Proof of the bounds in (\ref{part2ali})}} \par
By Lemma \ref{q0lemma}(\ref{qnode1}) either $q=q_0$ or $q=q_0^k$ for some $k\geq 2.$ We will prove (\ref{part2ali}) for these cases in turn.\par 
\par \underline{\textbf{Case 1: $q=q_0$}}\par 
\par 
\underline{\textit{Case 1.a: $l\geq 9$ and $n\geq{l+1 \choose 4}$ }} \par 
We know that a maximal $1$-space is fixed by a Borel subgroup, and hence there is an orbit $\mathcal{O}$ of $G_0$ with $\vert \mathcal{O} \vert \leq \frac{q^{\frac{(l^2+3l+2)}{2}}}{2}.$  Now we use Lemma \ref{lemma2.1}(\ref{boundeq2}) which gives us that $$d\geq \frac{2{l+1 \choose 4}}{l^2+3l+2}=\frac{l(l-1)(l-2)}{12(l+2)}$$ as required for conclusion (\ref{part2ali}). 
\par 
\underline{\textit{Case 1.b: $\lambda$ is $p$-restricted, $l\geq 9$ and $n<{l+1 \choose 4}$ }} \par 
By Theorem \ref{alvarobounds}, since $n<{l+1 \choose 4},$ $\lambda$ is as in Table \ref{tab:alvarodims}. Since $\lambda$ is not in Table \ref{ali}, we have $\lambda=3\omega_1,\omega_1+\omega_2,\omega_1+\omega_{l-1}$ or $2\omega_1+\omega_{l}$. Using Lemma \ref{para} and Lemma  \ref{orbit-stabparaborel}, we find upper bounds for the size of the orbit of the $1$-space fixed by the respective parabolic. Lemma \ref{lemma2.1} parts  (\ref{boundeq1})) and (\ref{boundeq2}) give us the bounds \begin{align}\label{alvaro}
     \begin{array}{c|c|c|c|c|c|c|c|c|c}
            \lambda& 3\omega_1&\omega_1+\omega_2&\omega_1+\omega_{l-1}& 2\omega_1+\omega_{l}\\\hline
            d\geq &\frac{(l+1)(l+2)(l+3)-6}{6(l+1)}& \frac{l(l+2)(l+5)-6}{6(2l+1)}&\frac{(l+1)(l^2+l+4)-2}{6l} &\frac{(l+1)(l^2+3l-2)-2}{2(2l+1)}
        \end{array}
\end{align} All of these bounds are more than $\frac{l(l-1)(l-2)}{12(l+2)},$ so the result follows in this case. \par 
\underline{\textit{Case 1.c: $\lambda$ is not $p$-restricted,  $l\geq 9$ and $n<{l+1 \choose 4}$}} \par 
By Theorems \ref{alvarobounds} and \cite[Thm 5.4.5]{kleidmanliebeck},  either
$\lambda=\mu_0+p^i\mu_1+p^j\mu_2,$ where each $\mu_i=\omega_1$ or $\omega_l,$ or $\lambda=\omega_1+p^i\omega_2$ or $\omega_1+p^i2\omega_1.$ Using Lemmas \ref{para}, \ref{lemma2.1}(\ref{boundformula1}) and Lemma \ref{orbit-stabparaborel} we get the bounds 
\begin{align}
    \begin{tabular}{c|c|c|c}\label{nonrestrictedal}
    $\lambda$ & $\lambda=\mu_0+p^i\mu_1+p^j\mu_2,$ & $\omega_1+p^i\omega_2$ & $\omega_1+p^i2\omega_1$ \\\hline
    $d\geq$ &  $\frac{(l+1)^3}{(2l+1)}$ &$\frac{l(l+1)^2}{4l+2}$ & $\frac{(l+1)(l+2)}{2}.$ 
\end{tabular}
\end{align} These bounds are all more than $\frac{l(l-1)(l-2)}{12(l+2)}.$
 
\underline{\textit{Case 1.d:  $l\leq 8$}} \par 

First assume $l=1.$ The $p$-restricted simple modules for $A_1(q)$ are $V=V(r\omega_1),$ the space of homogeneous polynomials in $x$,$y$ of degree $r.$ Then $V$ has dimension $r+1$ and basis $x^r,x^{r-1}y,\dots,y^r.$ The smallest orbit $\Delta$ of $G_0$ is the one containing $x^r,$ $$\Delta=\{(ax+by)^r\vert (a,b)\neq (0,0)\}.$$ Clearly $\vert\Delta\vert =q^2-1.$ Now we can use Lemma \ref{lemma2.1}(\ref{boundformula1}), which says $ 1+(q^2-1)+\dots+(q^2-1)^d \geq q^{n}$ and so $d\geq \frac{n}{2}=\frac{r+1}{2}.$ Hence $d\geq 3$ for $r\geq 4.$ 
Since $1+(q^2-1)+(q^2-1)^2 < q^4$, we also have $d\geq 3$ if $r=3$ Hence $d\geq 3$ in all cases. \par 
In fact, since a Borel fixes a maximal $1$-space, by Lemma \ref{orbit-stabparaborel} there is always an orbit of size at most $(q^2-1).$ Since $1+(q^2-1)+(q^2-1)^2 < q^4,$ if $n\geq 4,$ then by Lemma \ref{lemma2.1}(\ref{boundformula1}), $d\geq 3.$ The non-restricted cases with $n\leq 3$ are in Table \ref{ali}, so (\ref{part2ali}) holds for $l=1.$\par 
 
For $l=2$ or $3,$ we need to prove that $d\geq 3.$ We use the fact that a Borel fixes a maximal $1$-space, Lemma \ref{orbit-stabparaborel} and Lemma \ref{lemma2.1}(\ref{boundformula1}), which tells us that if $d\leq 2,$ then $n\leq 13$ or $18,$ respectively. By Theorem \ref{lowerlprestrict} and \cite[Thm 5.4.5]{kleidmanliebeck} all modules satisfying this are in Tables \ref{ali}, \ref{alvaro} or \ref{nonrestrictedal}.  \par 
Note that the bounds for the weights in (\ref{alvaro}) or (\ref{nonrestrictedal}) hold also for $l\leq 8$. These bounds are greater than those in (\ref{part2ali}) except when $(\lambda,l)=(\omega_1+p^i\omega_2,2).$ In this case, using that fact that a Borel subgroup fixes a maximal $1$-space we find that there is an orbit of size at most $(q^3-1)(q+1),$ so by Lemma \ref{lemma2.1}(\ref{boundformula1}) $d\geq 3,$ as required. \par 
Now assume $4\leq l \leq 8.$  By Theorem \ref{lowerlprestrict}, either $n\geq N,$ where $N$ is as in  (\ref{nagyN}), or $\lambda$ is in (\ref{alvaro}) or (\ref{nonrestrictedal}) or $l=7$ or $8$ and $(\lambda, n)=(\omega_4, {l+1 \choose 4}).$ \par 
Suppose $n \geq N.$ Then using the fact that a Borel fixes a maximal $1$-space, Lemma \ref{orbit-stabparaborel} and Lemma \ref{lemma2.1} we get the following lower bounds for $d.$ \par 
\begin{align}\label{nagyN}
    \begin{tabular}{c|c|c|c|c|c}
        $l$& 8&7&6&5&4\\\hline
        $N$&156&112&147&90&45\\\hline
        $d\geq$  &4&4 &6&4&3
    \end{tabular}
\end{align} 
Finally suppose $(\lambda, n)=(\omega_4, {l+1 \choose 4})$ and $7\leq l\leq 8.$ We can use Proposition \ref{para} to estimate the size of a small orbit and it follows that $d\geq 4$ for $7\leq l\leq 8.$ \par

\underline{\textbf{Case 2: $q=q_0^k$ for $k\geq 2$ as in Lemma \ref{q0lemma}(\ref{subfield})}}\par 
\underline{\textit{Case 2.a: $k=2$ }} \par 
Here Lemma \ref{q0lemma} gives $V(\lambda)=V(\lambda')\otimes V(\lambda')^{q_0}$ and $q=q_0^2.$ The case when $\lambda'=\omega_1$ is in Table \ref{ali} so is excluded. Hence $\dim(V(\lambda'))\geq \frac{1}{2}l(l+1)$ by Theorems \ref{alvarobounds} and \ref{lowerlprestrict} and so $n\geq \frac{l^2(l+1)^2}{4}.$ \par 
Using the fact that a Borel fixes a $1$-space, by Lemma \ref{orbit-stabparaborel}, we have $\vert \mathcal{O} \vert \leq \frac{q_0^{(l^2+3l+2)}}{2}$ and so using Lemma \ref{lemma2.1}(\ref{boundeq2}) it follows that \begin{equation*}\label{q0rlarge1}
   d\geq \frac{(l(l+1))^2}{4(l+1)(l+2)}.\end{equation*} 
This satisfies the bounds in (\ref{part2ali}) for $l\geq 4.$ \par 
Assume $l=3.$ For $n=\frac{(l(l+1))^2}{4}=36$ with $\lambda=\omega_2+q_0\omega_2,$ by Lemma \ref{para}, the parabolic $P_2$ fixes a $1$-space. Using Lemma  \ref{orbit-stabparaborel} and Lemma \ref{lemma2.1}(\ref{boundeq2}), $d\geq 4.$ For $\lambda\neq \omega_2+q_0\omega_2,$ by Theorem \ref{lowerlprestrict}, now $n\geq \frac{((l+2)(l+1))^2}{4}=100.$ Using the fact that Borel fixes a $1$-space, it follows that $d\geq \frac{(l+1)(l+2)}{4}=5$ which satisfies the bound for $l=3$ in (\ref{part2ali}). \par
Assume $l=2.$ By Theorem \ref{lowerlprestrict}, either $\lambda$ is in Table  \ref{ali} or  $n\geq 9.$ In the latter case, $d\geq 3$ as required for (\ref{part2ali}). \par
Assume $l=1.$ Now by Theorem \ref{lowerlprestrict}, when $\lambda\neq \omega_1+q_o\omega_1,$ we have $n\geq 9.$ Since a Borel fixes a maximal $1-$space,  by Lemma  \ref{orbit-stabparaborel} we have an orbit of size at most $(q_0-1)(q_0^2+1)$ so by Lemma \ref{lemma2.1}(\ref{boundformula1}) $d\geq 3$ as required. \par 
\underline{\textit{Case 2.b: $k\geq 3$ }} \par
Here we have $V(\lambda)=V(\lambda')\otimes \dots \otimes V(\lambda')^{q_0^{k-1}}.$ Now either $\lambda=\omega_1+q_0\omega_1+\dots+q_0^{k-1}\omega_1,$ (or their duals) or $\dim(V(\lambda'))\geq \frac{l(l+1)}{2}$ so $n\geq \frac{(l(l+1))^k}{2^k}.$ 
If $\lambda=\omega_1+q_0\omega_1+\dots +q_0^{k-1}\omega_1,$ then the parabolic $P_1$ fixes a $1$-space so have an orbit of size $\vert \mathcal{O} \vert \leq \frac{(q_0)^{kl+k}}{2},$ and so by Lemma \ref{lemma2.1}(\ref{boundeq2}) \begin{equation*}\label{q0w13}
    d\geq \frac{(l+1)^{k-1}}{k}.\end{equation*} This bound is greater than the one in (\ref{part2ali})  unless $(l,k)=(1,3)$ or $(1,4).$ In the latter cases, by Lemma  \ref{orbit-stabparaborel} we have an orbit of size at most $(q_0-1)(q_0^k+1)$ so by Lemma \ref{lemma2.1}(\ref{boundformula1}), $d\geq 3$ as required for (\ref{part2ali}). \par 
Finally, suppose $n\geq \frac{(l(l+1))^k}{2^k}.$ Again, using the fact that a Borel fixes a $1$-space, it follows that $\vert \mathcal{O} \vert \leq \frac{q_0^{k\frac{(l^2+3l+2)}{2}}}{2}$ and so by Lemma \ref{lemma2.1} part(\ref{boundeq1}) it follows that \begin{equation*}\label{q0rlarge2}
    d\geq \frac{2(l(l+1))^k}{k2^k(l+1)(l+2)},\end{equation*} so the bound in (\ref{part2ali}) holds. This concludes the proof of Theorem \ref{al(q)}.
\end{proof}

Now we can provide a complete classification of groups of the form $G$ as in Hypothesis \ref{hypothesis} with $\frac{G_0^\infty}{Z(G_0^\infty)}=A_l(q)$ which have orbital diameter 2, as stated in Theorem \ref{diam2}.

\begin{theorem}\label{aldiam2}
Let $G$ be as in Hypothesis \ref{hypothesis} with $\frac{G_0^\infty}{Z(G_0^\infty)}=A_l(q).$ Then $orbdiam(G,V)\leq 2$ if and only if one of the following holds.
\begin{itemize}
    \item $V$ is the natural $(l+1)$-dimensional module
    \item $(\lambda,l)=(\omega_2,4)$
     \item $(\lambda,l)=(\omega_2,3)$
      \item $(\lambda,l)=(2\omega_1,1)$ and $G_0$ contains the group $\mathbb{F}_{q_0}^*$ of scalars
     \item $(\omega_1+q_0\omega_1,1)$ and $q=q_0^2$
\end{itemize}
\end{theorem}

\begin{proof}
Assume $orbdiam(G,V)=2$. Looking at every lower bound in Theorem \ref{al(q)}, we get that either $orbdiam(G,V)\geq 3$ or $(\lambda,l)$ are as in the Table below.  \newline
\begin{tabular}{c|c|c|c|c|c|c|c|c}
   $\lambda$  & $\omega_2$& $\omega_2$&$2\omega_2$&$\omega_1+\omega_2$&$\omega_3$&$\omega_1+p^i \omega_2$&$\omega_1+q_0 \omega_1$&$\omega_1+q_0 \omega_1$\\ \hline
    $l$ & $4$&$3$&$1$&$2$&$5$&$2$&$1$&$2$ \\\hline
    extra conditions &&&&$3\vert q$&&&$q=q_0^2$&$q=q_0^2$
\end{tabular}
\par 
\underline{\textit{Case $(\lambda,l)=(\omega_2,4)$}}\newline
This was handled in Lemma \ref{smalldiamlemma}. We show here that this produces an example for orbital diameter $2$ even if $G_0$ does not contain the scalars in $GL_n(q_0).$
In this case $SL_5(q)\unlhd G_0\leq GL_{10}(q),$ and so by \cite{liebeckaffine}, $G_0$  has 2 orbits on $1$-spaces. Now $V=\bigwedge^2 W,$ where $W$ is the natural module of $SL_5(q).$ Let $v_1,\dots,v_5$ be the standard basis of $W.$ Then the two orbits of $G_0$ on the $1$-spaces of $V$ are $\langle v_1\wedge v_2\rangle^{G_0}$ and $\langle v_1\wedge v_2+v_3\wedge v_4\rangle^{G_0}.$ Since the diagonal matrices $diag(\lambda,1,\lambda^{-1},1,1)$ and $diag(\lambda,1,\lambda,1,\lambda^{-2})$ are in $SL_5(q),$ it has two orbits of non-zero vectors as well, so the diameter is 2 for any $G_0$ containing $SL_5(q)$.\par 
\underline{\textit{Case $(\lambda,l)=(\omega_2,3)$}}\newline
Now $\frac{SL_4(q)}{\langle-I\rangle}\cong \Omega^+_6(q)$ so $V$ is the natural module of $\Omega^+_6(q),$ and so by Lemma \ref{naturalmod}, $orbdiam(G,V)=2.$\par 
\underline{\textit{Case $(\lambda,l)=(2\omega_1,1)$ }}\newline
Now $PSL_2(q)\cong \Omega_3(q),$ and $V$ is the natural module of $\Omega_3(q),$ so by Lemma \ref{naturalmod}, $orbdiam(G,V)=2$ when $G_0$ contains the scalars in $GL_n(q_0).$\par 

\underline{\textit{Case $(\lambda,l)=(\omega_1+\omega_2,2)$  with $3\vert q$}}\newline
Now $(\lambda,l)=(\omega_1+\omega_2,2)$ with $3\vert q.$ As $3=l+1,$ the adjoint module can be identified with $V_{ad}=\{A\in M_{3}(q)\,\,\vert\,\, tr(A)=0\}/Z$ where $Z=\{\alpha I_{3}\vert \alpha\in \mathbb{F}_q\}.$ Let $E$ be a rank 1 traceless matrix. Now every element of the orbit of $Z+E$ will have a coset representative of rank 1 as conjugation preserves the rank. Hence showing that there exists a traceless matrix $A$ such that $rank(A+\alpha I_{3})=3$ for all $\alpha\in \mathbb{F}_q$ shows that $orbdiam(G,V)\geq 3.$ The companion matrix of an irreducible polynomial of the form $f(x)=x^3-bx-c$ with $b,c\in \mathbb{F}_q$ and $c\neq 0$ over $\mathbb{F}_q$ satisfies this property, so we want to show that such an irreducible polynomial exists. There are $q(q-1)$ polynomials of the form $x^3-\alpha x-\beta$ with $\beta \neq 0.$ This is reducible, if there is $\gamma, \delta\in \mathbb{F}_q$ such that $x^3-\alpha x-\beta=(x-\gamma)(x^2+\gamma x+\delta).$ Now $\gamma\delta=-\beta,$ which tells us that neither $\gamma$ or $\delta$ can be 0. Hence there are at most $(q-1)^2$ reducible such polynomials, so there is at least one irreducible polynomial of the form $f(x)=x^3-bx-c$ with $b,c\in \mathbb{F}_q$ and $c\neq 0.$ \par 
\underline{\textit{Case $(\lambda,l)=(\omega_3,5)$ }}\newline
In this case $V=\bigwedge^3W$ where $W$ is the natural module of $SL_6(q)\unlhd G_0.$ Since $SL_6(q)$ is transitive on the $3-$dimensional subspaces of $W$, $G_0$ has a single orbit on simple wedges, $w_1\wedge w_2\wedge w_3$. To prove that $orbdiam(G,V)\geq 3,$ it suffices to show that there are strictly fewer than $q^{20}$ distinct sums of at most two simple wedges. The number of simple wedges is $(q-1)$ times the number of 3-dimensional subspaces of $W.$ This is $$(q-1)\frac{(q^6-1)(q^5-1)(q^4-1)}{(q^3-1)(q^2-1)(q-1)}.$$ 
Now we want to count the number of sums of two simple wedges of the form $v_1\wedge v_2 \wedge v_3+w_1\wedge w_2 \wedge w_3$ with $v_i,w_i\in W.$ To do this we will first count the pairs of $3$-dimensional subspaces $A=\langle v_1,v_2,v_3\rangle$ and $B=\langle w_1,w_2,w_3\rangle.$  We have 3 cases to consider. If $dim (A\cap B)=2,$ then there are $x,y,z,k\in W$ such that $A=\langle x,y,z\rangle$ and $B=\langle x,y,k\rangle$ and so $x\wedge y\wedge z+x\wedge y\wedge k=x\wedge y\wedge (z+k),$ so $v_1\wedge v_2 \wedge v_3+w_1\wedge w_2 \wedge w_3$ is a simple wedge and we counted them already. 
For each pair $(A,B)$ such that $dim (A\cap B)=1$ or $dim (A\cap B)=0$ there are $(q-1)^2$ corresponding sums of two simple wedges. \par 
We count the number of pairs such that  $\dim (A\cap B)=1.$ We start with a $3-$dimensional subspace $A=\langle x_1,x_2,x_3\rangle.$ There are $X_a:=\frac{(q^6-1)(q^5-1)(q^4-1)}{(q^3-1)(q^2-1)(q-1)}$ choices. Choose a $1-$dimensional subspace of $A,$ which will be the intersection, call it $\langle x\rangle.$ There are $q^2+q+1$ choices for this. 
Let $\phi_1 \colon W\to \frac{W}{A}.$ 
 Then any $B$ such that $A\cap B=\langle x\rangle$ is of the form $B=Span(x,y_1,y_1),$ where $\langle \phi_1(y_1),\phi_1(y_2)\rangle $ is a $2$-dimensional subspace of $ \frac{W}{A}.$ So the number of such $B$s is $(q^2+q+1)q^4,$ so in the case when $\dim (A\cap B)=1,$ there are $\frac{1}{2}(q^6-1)(q^5-1)(q^2+1)(q^2+q+1)q^4$ sums of two simple wedges.

For $\dim (A\cap B)=0,$ by a similar argument, there are  $\frac{1}{2}(q-1)(q^3+1)(q^5-1)(q^2+1)q^9$ such sums. \par 
Adding up these quantities gives a value less than $q^{20},$ so $orbdiam(G,V)\geq 3.$ \par

\underline{\textit{Case $(\lambda,l)=(\omega_1+p^i\omega_l,2)$}}\newline
Now a Borel is the stabilizer of a maximal $1$-space, so using Lemma  \ref{orbit-stabparaborel} it folows that there is an orbit of size at most $(q^3-1)(q+1).$ As $1+(q^3-1)(q+1)+(q^3-1)^2(q+1)^2 <q^9 $ for $q\geq 2,$ this case has orbital diameter at least 3.\par 
\underline{\textit{Case $(\lambda,l)=(\omega_1+q_0\omega_1,1)$ with  $q_0^2=q$}}\newline
Now $PSL_2(q)\cong \Omega_4^-(q^{1/2}),$ and $V$ is the natural module of $\Omega_4^-(q^{1/2}),$ so by Lemma \ref{naturalmod} $orbdiam(G,V)=2.$\par 
\underline{\textit{Case $(\lambda,l)=(\omega_1+q_0\omega_1,2)$ with  $q_0^2=q$}}\newline
Define $V'$ to be the following $\mathbb{F}_{q_0}$-subspace of $M_3(q),$  $$V'=\{A\vert A^{(q_0)}=A^T\}=\langle \alpha E_{i,j} +\alpha^q E_{j,i} \colon \alpha \in \mathbb{F}_q,\,\,1\leq i,j\leq 3\rangle_{\mathbb{F}_{q_0}}.$$ 
 Here $g\in G_0$ acts on $A\in V$ as  $$A\to g^{T}Ag^\sigma$$ where $\sigma$ is the Frobenius morphism which raises matrix  entries to the power $q_0$ and  $V'$ is preserved by the action of $G_0.$ Hence we can identify $V$ with $V'.$ 
The rank of $A$ is also preserved by the $G_0$-action, so we cannot express $E_{1,1}+E_{2,2}+E_{3,3}$ as the sum of two elements of the orbit of $E_{1,1},$ so the orbital diameter is at least 3.\par

\end{proof}

\subsubsection{$G_0\triangleright\,^2\!A_l(q)$}
In this case we have that $G=V_n(q_0)G_0$ such that $\frac{G_0^\infty}{Z(G_0^\infty)}=\,^2\!A_l(q).$ Let $d=orbdiam(G,V)$. Note that $^2A_1(q)\cong A_1(q),$ so we can assume that $l\geq 2.$ Recall $\tau_0$ denotes a graph automorphism of $A_l.$
\begin{theorem}\label{2al(q)}
Let $G$ as in Hypothesis \ref{hypothesis} with  $\frac{G_0^\infty}{Z(G_0^\infty)}=\,^2\!A_l(q).$ 
\begin{enumerate}
    \item \label{1tau} If $\tau_0(\lambda)\neq \lambda,$ then $^2A_l(q)\leq A_l(q^2)\leq GL(V) $ and the lower bounds on $d$ in Theorem \ref{al(q)} hold. 
    \item \label{2tau} Suppose $\tau_0(\lambda)=\lambda.$  \begin{enumerate}[(i)]
    \item \label{2ali1} If $\lambda$ is in Table \ref{2ali}, the value of $n$ and a lower bound for $d$ are as given in the table.
\begin{table}[]
\[\begin{array}{c|c|c|c|c|c}
   & l&\lambda & n &d\geq \\\hline
 1. & 3&\omega_2&6&=2\\
 2. & \text{all}& \omega_1+\omega_l&(l+1)^2-1-\epsilon_p(l+1)&\frac{(l+1)^2-2}{2l+2}& \\\hline
\end{array}\]
    \caption{Bound in Theorem \ref{2al(q)} part \ref{2ali1}}
    \label{2ali}
    \end{table}   
    \item\label{part22ali} For $l\geq 9$ and $\lambda$ not in Table \ref{2ali}, we have $$d\geq \frac{(l-3)(l^2-l+4)}{12(l+1)}.$$ For $l\leq 8$ and $\lambda$ not in Table \ref{2ali}, we have the following bounds.
  
\begin{align}
     \begin{tabular}{c|c|c|c|c|c|c|c|c|c}
            l &8&7&6&5&4&3&2\\\hline
            $d\geq$ &26&4&12&3&5&3&3
        \end{tabular}
     \end{align}
\end{enumerate}  

\end{enumerate}

\end{theorem}
The lower bounds in the statement of Theorem \ref{2al(q)} are greater than the lower bounds contained in Theorem \ref{lhalf}, so Theorem \ref{lhalf} holds for $X_l(q)= \,^2A_l(q)$.

\begin{proof}[Proof of Theorem \ref{2al(q)}]
Part \ref{1tau} follows from Lemmas \ref{subgroup} and \ref{q0lemma}. Now we prove part \ref{2tau}, so from now on we assume that $\tau_0(\lambda)=\lambda.$ \par
We start by proving the bounds in Table \ref{2ali} in part (\ref{2ali1}). \par 
\underline{\textbf{Proof of the bounds in Table \ref{2ali}}}\par 


1. Here $\lambda=\omega_2,$ $l=3$ and $n=6.$ Now $\frac{G_0^{\infty}}{Z(G_0^{\infty})}\cong \Omega^-_6(q),$ and $V$ is the natural module of $\Omega^-_6(q)$, so the diameter is 2 by Lemma \ref{naturalmod}. \par 
2. Here $\lambda=\omega_1+\omega_l,$ so  $n=(l+1)^2-1-\epsilon_p(l+1),$ and $V=V(\lambda)$ is the adjoint module. In this case, the parabolic $P_{1,l}$ fixes a maximal $1$-space and so by Lemma  \ref{orbit-stabparaborel} there is an orbit of size at most $\vert \mathcal{O} \vert \leq \frac{q^{2l+2}}{2}.$ Now using Lemma \ref{lemma2.1}(\ref{boundeq2}) it follows that $d\geq \frac{(l+1)^2-2}{2l+2}.$\par

For the rest of the proof assume that $\lambda$ is not in Table \ref{2ali}. \par 
\underline{\textbf{Proof of the bounds in (\ref{part22ali})}}\par 
By Lemma \ref{q0lemma}, either $q=q_0$ or $q=q_0^k$ and $V=V(\lambda')\otimes V(\lambda')^{q_0}\otimes\dots\otimes V(\lambda')^{q_0^{k-1}}$ for some $p$-restricted dominant weight $\lambda'.$ 
\par \underline{\textbf{Case 1}: $q=q_0$} \par 
First we note that a Borel subgroup fixes a maximal $1$-space, so using Lemma \ref{orbit-stabparaborel} we see that there is an orbit $\mathcal{O}$ of size $\vert \mathcal{O} \vert \leq{q^{\frac{1}{2}(l^2+3l+2)}}.$\par 
\underline{\textit{Case 1.a: $\lambda$ is $p$-restricted and $l\geq 9$}}\par 

By Theorem \ref{alvarobounds}, for $l\geq 9,$  $n\geq {l+1 \choose 4}.$  Then, using the fact that a Borel fixes a $1$-space, Lemma \ref{lemma2.1}(\ref{boundeq1}) gives $d\geq{\frac{{l+1 \choose 4}-1}{\frac{1}{2}(l^2+3l+2)}}.$ This satisfies the bound in (\ref{part22ali}).
 
\underline{\textit{Case 1.b: $\lambda$ is $p$-restricted and $2\leq l\leq 8$}}\par  
 For $2\leq l\leq 8,$ by Theorem \ref{lowerlprestrict}, 
 $(l,n,\lambda)=(7,70,\omega_4),(5,20,\omega_3),$ $(3,19,\omega_2)$ or $n$ and $d$ are bounded below as follows: \par 
 \begin{align}\label{smallproofs}
     \begin{tabular}{c|c|c|c|c|c|c|c}
  $l$   &2&3& 4&5&6&7&8 \\ \hline
   $n\geq$  &19&44& 74&154&344&657&1135\\\hline
   $d\geq $ &3&5&5&8&12&19&26
   \end{tabular}.
    \end{align}\par
In the cases as $(l,n,\lambda)=(7,70,\omega_4)$ (respectively  $(5,20,\omega_3),$ $(3,19,\omega_2)$), a parabolic $P_4$ (respectively $P_3,$ $P_2$) fixes a maximal $1$-space, so there is a $G_0$-orbit on $V$ of size at most $q^{19}$ (respectively $(q^5+1)(q^3+1)(q-1)^2,$ $(q^4-1)(q^3+1)$). Now Lemma \ref{lemma2.1}(\ref{boundeq1}) gives that $d\geq 4$ (respectively 3,3).\par 
\underline{\textit{Case 1.c: $\lambda$ is not $p$-restricted}}\par 

 If $V$ is not $p$-restricted, then by Theorems \ref{alvarobounds} and \ref{lowerlprestrict}, $n\geq s^2$ where $s$ is as follows:
      \[\begin{array}{c|c|c|c}
         l & 3&5&l\neq 3\,\, \text{or} \,\, 5  \\\hline
        s  & 6&20&(l+1)^2-2.\\\hline
        \end{array}\] By the fact that a Borel fixes a $1$-space and using Lemma \ref{lemma2.1}(\ref{boundeq1}),  it follows that $d \geq \frac{s^2-1}{(l^2+3l+2)/2}.$ This satisfies the bounds in (\ref{part22ali}) . \par 
 \underline{\textbf{Case 2: $q=q_0^k$ for $k\geq 2$ as in Lemma \ref{q0lemma}(\ref{subfield})}}\par 
 Here $V(\lambda)=V(\lambda')\otimes \dots \otimes V(\lambda')^{q_0^{k-1}}$ for some $p$-restricted $\lambda'.$ Using the fact that a Borel fixes a maximal $1$-space we have that $\vert\mathcal{O}\vert \leq q_0^\frac{{kl^2+3kl+2}}{2}$ and so by Lemma \ref{lemma2.1}(\ref{boundeq1}) it follows that $d\geq \frac{2s^k-2}{{kl^2+3kl+2}}\geq\frac{2s^2-2}{{{2l^2+3(2l)+2}}}.$ This satisfies the bounds in (\ref{part22ali}) except for $l=3.$  For the case $l=3,$ $\vert\mathcal{O}\vert \leq q_0^{6k+4}$ and so using Lemma \ref{lemma2.1}(\ref{boundeq1}) it follows that $d\geq \frac{6^k-1}{6k+4}\geq \frac{35}{16}$ and so $d\geq 3.$ 
 \end{proof}

 \begin{prop}\label{2aldiam2}
Let $G$ be as in Hypothesis \ref{hypothesis} such that $\frac{G_0^\infty}{Z(G_0^\infty)}=\,\!^2A_l(q).$ Then $orbdiam(G,V)\leq 2$ if and only if one of the following holds.
\begin{itemize}
    \item $V$ is the $(l+1)$-dimensional natural module.
     \item $(\lambda,l)=(\omega_2,3).$
\end{itemize}
 \end{prop}
\begin{proof}
\par \underline{\textbf{Case 1}:  $\tau_0(\lambda)\neq \lambda$}\par 
 By Theorem \ref{aldiam2} the only candidates for $d=2$ are $\lambda=\omega_1$ and  $(\lambda,l)=(\omega_2,4)$ with $q^2=q_0$ in all cases. For $\lambda=\omega_1,$ $V$ is the natural module, so the orbital diameter is $2$ by Lemma \ref{naturalmod}.
 
In the other case $^2A_4(q)\leq A_4(q^2)\leq GL_{10}(q^2)$ and a parabolic $P_2$ fixes a maximal $1$-space, so there is an orbit of size $\vert\mathcal{O}\vert \leq (q^2-1)\frac{(q^5+1)(q^3+1)}{(q+1)}$ and so by Lemma \ref{lemma2.1}(\ref{boundformula1}), $d\geq 3.$\par  
 \underline{\textbf{Case 2}:  $\tau_0(\lambda)= \lambda$}\par  By Theorem \ref{2al(q)} the candidates for $d\leq 2$ are $(\lambda,l)=(\omega_2,3),$ $(\omega_1+\omega_3,3)$ and $(\omega_1+\omega_2,2).$ \par 
 \underline{\textit{$(\lambda,l)=(\omega_2,3)$}}\par
 In this case $d=2$ by Lemma \ref{smalldiamlemma}. \par 
 \underline{\textit{$(\lambda,l)=(\omega_1+\omega_3,3)$}}\par 
Here $n=(l+1)^2-1-\epsilon_p(l+1)=15-\epsilon_p(l+1)\geq 14.$ The parabolic subgroup $P_{1,3}$ fixes a maximal $1$-space and its orbit has size $\vert \mathcal{O} \vert \leq (q^3+1)(q^2+1)(q-1).$ Now by Lemma \ref{lemma2.1}(\ref{boundformula1}),  $d\geq 3.$ \par
  \underline{\textit{$(\lambda,l)=(\omega_1+\omega_2,2)$}} \par 
Here $n$ is either $8$ or $7.$ If $n=8,$ then a Borel fixes a maximal $1$-space, so by Lemma  \ref{orbit-stabparaborel} there is an orbit of size $\vert\mathcal{O}\vert \leq(q^3+1)(q-1),$ and so by Lemma \ref{lemma2.1}(\ref{boundformula1}), $d\geq 3.$ \par 
Now we consider the case when $(\lambda,l)=(\omega_1+\omega_2,2)$ with $3\vert q,$ and $V$ is the adjoint module. Here $V$ can be identified with $V_{ad}=\{A\in M_{3}(q^2)\,\,\vert\,\, A+\overline{A}^T=0, tr(A)=0\}/Z$ where $\overline{A}=A^{(q)}$ and $Z=\{\alpha I_{3}\vert \alpha+\overline{\alpha}=0,\alpha\in \mathbb{F}_{q^2}\}.$ Here $G_0$ acts by conjugation on $V_{ad}$. Let $E$ be a rank 1 traceless matrix in $V$. Now every element of the orbit of $Z+E$ will have a coset representative of rank 1 as conjugation preserves the rank. Hence showing that there exists a traceless matrix $A$ such that $rank(A+\alpha I_{3})=3$ for all $\alpha\in \mathbb{F}_{q^2}$ shows that $orbdiam(G,V)\geq 3.$ Hence it is sufficient to show that there is a rank 3 matrix in $V_{ad}$ with irreducible characteristic polynomial. 
Fix $a \in  \mathbb{F}_{q^2}$ such that $a \bar a = -1$, and define the $3\times 3$ matrix $$M_b=\begin{pmatrix}
0 & a & b\\
-\overline{a}& 0 & 1\\
-\overline{b}&-1&0
\end{pmatrix}$$ with characteristic polynomial $x^3+b \bar bx + a \bar b- \bar ab$. Fix $\alpha \in  \mathbb{F}_q^*$ such that $-\alpha$ is a nonsquare, and define
\[
S_\alpha = \{ \beta\in  \mathbb{F}_{q^2} :  \beta\bar \beta=\alpha\},\;\;T = \{ \beta\in  \mathbb{F}_{q^2} : \beta+ \bar  \beta= 0\}.
\]
Then $|S_\alpha| = q+1$, and $T$ is a subgroup of $ \mathbb{F}_{q^2}^+$ of size $q$. For $b \in S_\alpha$, the matrix $M_b$ has characteristic polynomial 
\[
c_b(x) = x^3+\alpha x+a \bar b- \bar ab.
\]
 We shall show that $b \in S_\alpha$  can be chosen so that $c_b(x)$ is irreducible (over $ \mathbb{F}_{q^2}$).

We first count the number of reducible cubics $x^3+\alpha x + \beta$ with $ \beta\in T$. To do this, define $\phi:T\mapsto T$  to send $x \mapsto x^3+\alpha x$ for $x \in T$. Then $\phi$ is an additive homomorphism, and ${\rm ker}(\phi)$ consists of the solutions of $x(x^2+\alpha)=0$. As we chose $-\alpha$ to be a nonsquare in $ \mathbb{F}_q$, it has two square roots in $ \mathbb{F}_{q^2}$ which we write as $\pm \gamma$; moreover $ \bar \gamma$ is also a solution, so $ \bar \gamma = -\gamma$ and so $\gamma \in T$. Thus ${\rm ker}(\phi) = \{0,\pm \gamma\}$ and so ${\rm Im}(\phi) =  \frac{q}{3}$. Thus there are $ \frac{q}{3}$  reducible cubics $x^3+\alpha x + \beta$ with $ \beta\in T$. (Note that any root in $ \mathbb{F}_{q^2}$ of such a cubic lies in $T$.)

Now we count the number of distinct cubics $c_b(x)$ for $b \in S_\alpha$. This is just the number of distinct elements 
$a \bar b- \bar ab$ for $b \in S_\alpha$. Define $\pi: S_\alpha \mapsto T$ to send $b \mapsto a \bar b- \bar ab$. For $b_1,b_2 \in S_\alpha$,  
\[
\begin{array}{ll}
\pi(b_1)=\pi(b_2) & \Leftrightarrow a( \bar b_1- \bar b_2) =  \bar a(b_1-b_2) \\
 &  \Leftrightarrow  a^2\alpha ( \frac{1}{b_1}- \frac{1}{b_2}) = b_2-b_1 \\
 &  \Leftrightarrow  a^2\alpha = b_1b_2 \\
 &  \Leftrightarrow  b_2 =  \frac{a^2\alpha}{b_1}.
\end{array}
\]
It follows that the image of $\pi$ has size at least $ \frac{1}{2}|S_\alpha| =  \frac{1}{2}(q+1)$. Since $ \frac{1}{2}(q+1)> \frac{q}{3}$, by the previous paragraph it follows that there exists $b \in S_\alpha$ such that $c_b(x)$ is irreducible (over $ \mathbb{F}_{q^2}$), as required. 

\end{proof}
\subsubsection{$G_0\triangleright C_l(q)$}
Recall that $C_1(q)\cong A_1(q).$
\begin{theorem}\label{cl(q)}
Let $G$ be as in Hypothesis \ref{hypothesis} such that $\frac{G_0^\infty}{Z(G_0^\infty)}=C_l(q)$ with $l\geq 2.$ 
\begin{enumerate}[(i)]
    \item \label{cli1} If $\lambda$ is in Table \ref{cli}, and the value of $n$
 and a lower bound for $d$ are as given in the table.
 \begin{table}[]

 \begin{center}
\begin{tabular}{c|c|c|c|c|c}
   & $\lambda$ & $n$ &$d\geq$&extra conditions  \\\hline
   1.& $\omega_1$&$2l$&$=1$&$q=q_0$\\
    2.& $\omega_2$ &$2l^2-l-1-\epsilon_p(l)$&$\frac{2l^2-l-2}{4l-2}$ & 
  $q=q_0$\\
    3.& $2 \omega_1$&$2l^2+l$&$\frac{2l+1}{2}$&$q=q_0$\\
   4.&  $\omega_1+p^i\omega_1$&$4l^2$&$2l$&$q=q_0$\\
5.&$\omega_1+q_0\omega_1$&$4l^2$&$l$& $q=q_0^2$\\
6.&$\omega_l$&$2^l$&$=2$&$l=3$ or $4,$ $q=q_0$ and $p=2$ \\\hline

\end{tabular}
\end{center}
     \caption{Bounds in Theorem \ref{cl(q)}(\ref{cli1})}
     \label{cli}
 \end{table}
    \item \label{part2cli}For $l\geq 14$ and $\lambda$ is not in Table \ref{cli}, we have
    $$
        d\geq \frac{l(4l^2-6l-10)-3}{18l-30}
    $$ For $l\leq 13$ and $\lambda$ is not in Table \ref{cli}, the lower bound is as follows.
    \begin{align*}
     \begin{tabular}{c|c|c|c|c|c|c|c|c|c|c|c|c|c|c|c}
            l &13&12&11&10&9&8&7&6&5&4&3&2\\\hline
            $d\geq$&13&12&11&10 &6&4&3&3&3&3&3&3
        \end{tabular}
     \end{align*}
\end{enumerate}

\end{theorem}
The lower bounds in the statement of Theorem \ref{cl(q)} are greater than the lower bounds contained in Theorem \ref{lhalf}, so Theorem \ref{lhalf} holds for $X_l(q)=C_l(q)$. 

\begin{proof}[Proof of Theorem \ref{cl(q)}]

\underline{\textbf{Proof of the bounds in Table \ref{cli}}}

1. Here $\lambda=\omega_1$ and $V$ is the natural module so the result follows from Lemma \ref{naturalmod}. \par
2-5. These cases are proved by Lemma \ref{orbit-stabparaborel} in the usual way. \par 
 6. Here $\lambda=\omega_l,$ $n=2^l$ and $l=3$ or $4.$ Both of these cases are in  Lemma \ref{smalldiamlemma}.


From now on, assume that $\lambda$ is not in Table \ref{cli}. \par 
\underline{\textbf{Proof of the bounds in part (\ref{part2cli})}}\par 
Recall that either $q=q_0$ or $q=q_0^k$ and $V=V(\lambda')\otimes V(\lambda')^{q_0}\otimes\dots\otimes V(\lambda')^{q_0^{k-1}}$ for some $p$-restricted dominant weight $\lambda'$ by Lemma \ref{q0lemma}. 
\par \underline{\textbf{Case 1}: $q=q_0$}\par 
 \underline{\textit{Case 1.a}: $\lambda$ is $p$-restricted and $l\geq 14$}\par 
 By Theorem \ref{alvarobounds}, if $n< 16{l \choose 4},$ then either $\lambda=\omega_3,$ $3\omega_1$ or $\omega_1+\omega_2.$ Using Lemma \ref{para} and Lemma  \ref{orbit-stabparaborel}, we find upper bounds for the orbit of a maximal $1$-space fixed by the respective parabolic. Lemma \ref{lemma2.1} parts (\ref{boundformula1}) and(\ref{boundeq1}) give us the bounds 
\begin{align*}
\begin{tabular}{c|c|c|c|c}
    $\lambda$ & $\omega_3$ &$3\omega_1$ & $\omega_1+\omega_2$ \\\hline
    $d\geq$ & $\frac{l(4l^2-6l-10)-3}{18l-30}$& $\frac{ (1 + 2 l) (2 + 2 l)}{6}$& $\frac{l(4l^2+6l-4)-3}{12l-6}$\\\hline
\end{tabular}
\end{align*}
These satisfy the bound in part (\ref{part2cli}).\par 
 Now suppose $n\geq 16{l \choose 4}.$ A Borel subgroup fixes a maximal $1$-space, so using Lemma  \ref{orbit-stabparaborel} we see that there is an orbit $\mathcal{O}$ of size $\vert \mathcal{O} \vert \leq{q^{l^2+l}},$ so we can use Lemma \ref{lemma2.1}(\ref{boundeq1}), which gives that $$d\geq \frac{16{l \choose 4}-1}{l^2+l}.$$ This satisfies the bound in part (\ref{part2cli}).\par 

 \underline{\textit{Case 1.b}: $\lambda$ is $p$-restricted and $2\leq l\leq 13$}\par 
We can see from Theorem \ref{lowerlprestrict} that either $\lambda$ is as in (\ref{smallclexceptions}) or the lower bounds for $n$ and $d$ in (\ref{mdiam}) hold. 
\begin{align}\label{smallclexceptions}
    \begin{tabular}{c|c|c|c|c|c|c}
    &1.&2.&3.&4.&5.&6.\\ \hline
    $l$ & $6$& $5$&$4$&$4$&$3$&$2$\\\hline
    $\lambda$ & $\omega_6$&$\omega_5$&$\omega_3$&$\omega_4$&$\omega_3$&$2\omega_2$\\\hline
    $n$& $64$&$32$&$48-8\epsilon_3(p)$&$\geq 41,$ &$14$&$10$\\\hline
    $p$&$2$&$2$&&&&\\\hline
\end{tabular}
\end{align} 

\begin{align}\label{mdiam}
    \begin{tabular}{c|c|c|c|c|c|c|c|c|c|c|c|c}
    $l$& $2$&$3$&$4$&$5$&$6$&$7$&$8$&$9$&$10$&$11$&$12$&$13$\\\hline
       $n\geq $  & $11$&$25$&$64$&$100$&$208$&$128$&$256$&$512$&$1000$&$1331$&$1728$&$2197$ \\\hline
       $d\geq $  &$3$&$3$&$4$&$4$&$3$&$3$&$4$&$6$&$10$&$11$&$12$&$13$ \\\hline
    \end{tabular}
\end{align}
Using the fact that a Borel fixes a maximal $1$-space we obtain the bounds in (\ref{mdiam}) using Lemma \ref{lemma2.1}(\ref{boundeq1}) for $4\leq l\leq 13$ and part (\ref{boundformula1}) for $2\leq l\leq 3.$
Now we find a lower bound for $d$ for all cases in (\ref{smallclexceptions}) in turn. 

1. Here $(l,\lambda,n,p)=(6,\omega_6,64,2).$  \par 
   By Lemma \ref{para}  the parabolic $P_6$ fixes a maximal $1$-space. Using Lemma  \ref{orbit-stabparaborel}, $\vert \mathcal{O} \vert \leq q^{27}$ and so by Lemma \ref{lemma2.1}(\ref{boundformula1}), $d\geq 3.$\par 
2. Here $(l,\lambda,n,p)=(5,\omega_5,32,2).$ \par
   Here $C_5(q)\leq D_6(q)\leq GL_{32}(q)$ so it suffices to show that for $G_0\triangleright D_6(q)$ the orbital diameter is at least 3. 
\begin{claim}\label{d6}
 Let $G=VG_0$ be a primitive affine group such that $\frac{G_0^\infty}{Z(G_0^\infty)}=D_6(q)$ with $V=V(\omega_5)=V_{32}(q).$ Then $orbdiam(G,V)\geq 3.$  
\end{claim}
\begin{proof}[Proof of Claim \ref{d6}]
First consider the algebraic group $\overline{G_0}=D_6(K),$
 where $K$ is the algebraically closed field $\mathbb{\overline{F}}_p$ acting on $\overline{V}=V_{32}(K).$ All stabilizers are listed in \cite{igusa} and \cite[Proof of Lemma 2.11]{glms}.\par
 Let $\Delta_0$ be the orbit of $\overline{G_0}$ on $\overline{V}$ of a maximal vector. We will call the elements of $\Delta_0$ pure spinors. The stabilizer of the 1-space of a pure spinor is $P_6=Q_{15}A_5T_1,$ and of a pure spinor is $P_6'.$ By \cite{igusa}, there exists $v\in V$ such that $(\overline{G_0})_v=Q_{14}C_3.$ We want to show that we cannot express $v$ as a sum of at most two pure spinors. It is sufficient to show that for all $g\in \overline{G_0},$ $P_6'\cap P_6'^g$ is not contained in $Q_{14}C_3.$ We prove this by contradiction. \par 
 Suppose that $P_6'\cap P_6'^g\leq Q_{14}C_3.$ Recall that by Lemma \ref{fullrank}, $T_6\leq P_6\cap P_6^g.$ By the second isomorphism theorem, $$T_1\cong \frac{P_6}{P_6'}\cong \frac{T_6P_6'}{P_6'}\cong \frac{T_6}{T_6\cap P_6'},$$ and so $$T_5\leq T_6\cap P_6'\leq P_6'\cap P_6^g.$$ By applying the second isomorphism theorem again, $T_4\leq P_6'^g\cap T_5\leq P_6'\cap P_6'^g.$ Since $T_4$ is not contained in $Q_{14}C_3,$ we reach a contradiction. Hence, we cannot express any elements of $\Delta$ as a sum of two pure spinors. \par 
 Now consider the finite group $G_0=\overline{G_0}^{(q)}$ acting on $V=\overline{V}^{(q)}=V_n(q).$ Choose $w\in V$ such that $w\in \Delta.$ By the above argument, there does not exists $a,b \in \Delta_0$ such that $w=a+b.$ Now it follows that there does not exists $a,b\in \Delta_0^{(q)}\cap V$ such that $w=a+b$ either. Observe that the orbit of pure spinors is preserved by all automorphisms of $D_6(q),$ so the Claim now follows.

\end{proof} \par 
3-6 (in  (\ref{smallclexceptions})). Using Lemma \ref{orbit-stabparaborel} we find upper bounds for the size of the orbit of a maximal vector, whose 1-space is stabilized by the respective parabolic. For 3. and 4. Lemma \ref{lemma2.1}(\ref{boundeq1}) gives the bounds of $d\geq 3$ and $d\geq 4,$ respectively. For 5. and 6. by Lemma \ref{lemma2.1}(\ref{boundformula1}), $d\geq 3.$
\underline{\textit{Case 1.c: $\lambda$ is not $p$-restricted}}\par 
By Theorems \ref{alvarobounds} and \ref{lowerlprestrict} we have that either $n\geq l^4$ or $\lambda=\omega_1+p^i \omega_1+p^j\omega_1,$ $\omega_1+p^i\omega_2,$  $\omega_1+p^i2\omega_1,$ or $\omega_1+p^i\omega_l$ with $3\leq l\leq 7.$ \par 
For $n\geq l^4$ we use the fact that a Borel fixes a maximal $1$-space to get that $d\geq \frac{l^4}{l^2+l}$ and so part (\ref{part2cli}) is satisfied.\par
For the other possibilities, using Lemma \ref{para} and Lemma  \ref{orbit-stabparaborel}, we find upper bounds for the orbit of the $1$-space fixed by the respective parabolic. Lemma \ref{lemma2.1} parts (\ref{boundformula1}) and(\ref{boundeq1}) give us the bounds 
\begin{align}
\begin{tabular}{c|c|c|c|c}
    $\lambda$ & $\omega_1+p^i \omega_1+p^j\omega_1$& $\omega_1+p^i\omega_2$&  $\omega_1+p^i2\omega_1$ &$\omega_1+p^i\omega_l$  \\\hline
    $d\geq$ & $4l^2$&$\frac{2l(2l^2-l-2)}{4l-2}$&$2l^2+l$&$\frac{2^{l+1}-2}{l^2-l}$\\\hline
\end{tabular}
    \end{align}
These satisfy the bound in part  (\ref{part2cli}). \par 
 \underline{\textbf{Case 2: $q=q_0^k$ for $k\geq 2$ as in Lemma \ref{q0lemma}(\ref{subfield})}}\par 
Here $V(\lambda)=V(\lambda')\otimes V(\lambda')^{q_0}\otimes \dots \otimes V(\lambda')^{q_0^{k-1}}.$ We can see from Theorems \ref{alvarobounds} and \ref{lowerlprestrict} that either $\lambda'=\omega_1$ and $n=(2l)^k$ or $\lambda'\neq \omega_1$ and $\dim V(\lambda')\geq l^2-1$ so that $n\geq (l^2-1)^k.$  Note that the second lowest dimension is usually even higher than $l^2-1,$ but we choose this value as this is a lower bound that works for all values of $l,$ in particular also for $l=3,$ where the second lowest dimension is $2^3=3^2-1=8.$ We will consider these two cases in turn.\par 
\underline{$n=(2l)^k$ and $\lambda=\omega_1+q_0\omega_1+\dots q_0^{k-1}\omega_1.$}\par 
Here $k\geq 3,$ as the $k=2$ case is in Table \ref{cli}. In this case the parabolic subgroup $P_1$ fixes a $1$-space and so by Lemma  \ref{orbit-stabparaborel}, $\vert\mathcal{O}\vert \leq \frac{q_0^{2kl}}{2}$ and so using Lemma \ref{lemma2.1}(\ref{boundeq2}) we can see that $d\geq \frac{(2l)^k}{2kl}\geq \frac{(2l)^3}{6l}.$ This satisfies the bound in part  (\ref{part2cli}). \par 
\underline{$\lambda\neq\omega_1+q_0\omega_1+\dots q_0^{k-1}\omega_1$}\par 
Now assume $\lambda'\neq \omega_1.$
Here $n> (2l)^k,$ and so $n\geq (l^2-1)^k$ by Theorems \ref{alvarobounds} and \ref{lowerlprestrict}. Since a Borel is fixing a maximal $1$-space, by Lemma \ref{orbit-stabparaborel}, $\vert\mathcal{O}\vert \leq \frac{q_0^{k(l^2+l)}}{2}$ and so by Lemma \ref{lemma2.1}(\ref{boundeq2}) $d\geq \frac{(l^2-1)^k}{k(l^2+l)}\geq \frac{(l^2-1)^2}{2l^2+2l}.$ This satisfies the bound in part  (\ref{part2cli})  for $l\geq 3.$ For $l=2$ we have $n\geq 5^k,$ and so $d\geq \frac{5^k}{6k}\geq \frac{5^2}{12}$ satisfying part (\ref{part2cli}). 
\end{proof}
We can achieve the following classification.
 \begin{prop}\label{cldiam2}
Let $G$ be as in Hypothesis \ref{hypothesis} such that $\frac{G_0^\infty}{Z(G_0^\infty)}=C_l(q).$ Then $orbdiam(G,V)\leq 2$ if and only if one of the following holds.
\begin{itemize}
    \item $V$ is the natural module.
     \item $G_0\triangleright C_3(q)$ and $V=V_{8}(q)$ with $q$ even.
     \item $G_0\triangleright C_4(q)$ and $V=V_{16}(q)$ with $q$ even.
\end{itemize}
 \end{prop}
 \begin{proof}
 By Theorem \ref{cl(q)}, the candidates for $orbdiam(G,V)\leq 2$ are $(\lambda,l)=(\omega_2,2),$ $(\omega_2,3),$  $(\omega_2,4),$ and $(\omega_3,3)$ and $(\omega_4,4)$ with $p=2,$ all of these with $q=q_0,$ and $(\omega_1+q_0\omega_1,2)$ with $q=q_0^2.$ \par
 \underline{ $(\omega_3,3)$ or $(\omega_4,4)$} These cases produce an example for a group with orbital diameter 2 by Lemma \ref{smalldiamlemma}.\par 
 \underline{$(\lambda,l)=(\omega_2,2)$.} This is the natural module for $C_2(q)\cong B_2(q)$ so the orbital diameter is $2$ by Lemma \ref{naturalmod}.\par 
 \underline{$(\lambda,l)=(\omega_2,3)$.}
Let $W$ be the natural module for $G_0^\infty=Sp_6(q)$. By \cite[page 103]{liebeckshalevbases}, $V(\lambda)$ is an irreducible 
composition factor of $\bigwedge^2W.$ Let ${e_1,e_2,e_3,f_1,f_2,f_3}$ be a standard basis of $W$ and let $J=\sum_{i=1}^3 e_i\wedge f_i$ and define a symmetric bilinear form on $\wedge^2W$ by $(v\wedge w,v'\wedge w')=B(v,v')B(w,w')-B(v,w')B(w,v')$ as in \cite[page 103]{liebeckshalevbases}. Now $(J,J)=0$ if and only if $p=3.$ For $p\neq 3$ we take $V=J^\perp$ of dimension $14$ and for $p=3$ we take $V=\frac{J^\perp}{\langle J\rangle}$ of dimension 13. \par 
For $3\nmid q,$ there is an orbit containing only simple wedges and $e_1\wedge e_2+ f_3\wedge f_2+e_3\wedge f_1 \in J^\perp$ cannot be expressed as a sum of at most 2 simple wedges, so $d\geq 3.$\par 
Now consider the case when $n=13,$ so $p=3$ and $V=\frac{J^\perp}{\langle J\rangle}.$ Consider the algebraic group $\overline{G_0}=C_3(K)$ where $K=\overline{\mathbb{F}_3}$ acting on $\overline{V}=V_{13}(K).$ Now the stabilizer of a maximal 1-space in $G_0$ is $P_2.$ Let $\Delta_0$ be the orbit of $\overline{G_0}$ containing the maximal vectors. If $v=a+b,$ where $a,b\in \Delta_0,$ then $\overline{G_0}_a\cap \overline{G_0}_b\leq \overline{G_0}_v.$ The stabilizers of $a$ and $b$ are conjugates of $P_2'.$ Without loss of generality, they are $P_2'$ and $P_2'^g$ for some $g\in \overline{G_0}.$ We see from Lemma \ref{alunalemma} what the intersections of two parabolics can be. Now we will show that if $H\leq P_2\cap P_2^g$ and $H$ is either a unipotent subgroup or $A_1A_1,$ then $H\leq P_2'\cap P_2'^g.$ By the second isomorphism theorem, $$\frac{H}{H\cap P_2'}\cong \frac{HP_2'}{P_2'}\leq \frac{P_2}{P_2'}\cong T_1.$$ Since the only unipotent or simple subgroup of $T_1$ is the identity, we can deduce that $H\cong H\cap P_2',$ so $H\leq P_2'.$ Similarly we can see that $H\leq P_2'^g.$ \par 
Note that $C_1(K)^3$ 
is a subgroup of $\overline{G_0}\cong C_3(K).$ Let $\overline{W_6}$ be the natural module of $C_3(K)$ and $\overline{W_2}$ the natural module of $C_1(K).$ Then $$\wedge^2 \overline{W}_6\downarrow C_1(K)^3=\wedge^2 (\overline{W}_2^1+\overline{W}_2^2+\overline{W}_2^3)=(\wedge^2 \overline{W_2})^3+\sum_{1\leq i\neq j\leq 3} \overline{W}_2^i \otimes \overline{W}_2^j. $$ Since $\wedge^2 \overline{W_2}$ is trivial and   $\overline{V}\downarrow C_1(K)^3=\bigwedge^2 \overline{W_2}+\sum_{1\leq i\neq j\leq 3} \overline{W}_2^i \otimes \overline{W}_2^j,$ it follows that $C_1(K)^3$ fixes a vector in $\overline{V}.$ \par 
Let $\sigma_q$ be the standard Frobenius morphism of $\overline{G_0}$ and let $\omega\in \overline{G_0}$ be the map permuting the terms in $X=C_1(K)^3,$ so for $(x_1,x_2,x_3)\in X,$ $\omega$ maps $(x_1,x_2,x_3)$ to $(x_3,x_1,x_2).$ By the Lang-Steinberg Theorem \cite[Theorem 21.7]{malletesterman}, $\overline{G_0}^{\sigma_q}\cong \overline{G_0}^{\sigma_q\omega}\cong C_3(q),$ acting on $V=\overline{V}^{(q)}=V_{13}(q).$ Also $X^{\sigma_q\omega}=\{(x,x^{(q)},x^{(q^2)})\vert x\in C_1(q^3)\}\cong C_1(q^3)\leq C_3(q),$ which fixes a vector in $V.$ \par 
The possible intersections of $P_2'\cap P_2'^g$ with $\overline{G_0}^{(q)}$ by Lemma \ref{alunalemma} contain either a unipotent subgroup of order at least $q^{5}$ or $A_1(q)^2.$  Since these are not contained in $C_1(q^3),$ we cannot express $w$ as a sum of at most two elements in $\Delta_0\cap V,$ so the orbital diameter is at least 3.

   \par 
\underline{$(\lambda,l)=(\omega_2,4)$.} Now $n=27-\epsilon_2(p)$ and the parabolic $P_2$ fixes a maximal $1$-space. Hence by Lemma  \ref{orbit-stabparaborel} we have that $\vert\mathcal{O}\vert \leq \frac{(q^8-1)(q^6-1)}{(q+1)}$ and so by Lemma \ref{lemma2.1}(\ref{boundformula1}) $d\geq 3.$ \par 
\underline{$(\lambda,l)=(\omega_1+q_0\omega_1,2)$} Using the fact that the parabolic $P_1$ fixes a maximal $1$-space and Lemma  \ref{orbit-stabparaborel}, it follows that $\vert\mathcal{O}\vert \leq (q-1)(q^2+1)(q^4+1)$ and so by Lemma \ref{lemma2.1}(\ref{boundformula1}) again $d\geq 3.$\par 
 \end{proof}
\subsubsection{$G_0\triangleright B_l(q)$}
Recall that $B_1(q)\cong A_1(q),$  $B_2(q)\cong C_2(q)$ and $B_l(2^r)\cong C_l(2^r).$
\begin{theorem}\label{bl(q)}
Let $G$ be as in Hypothesis \ref{hypothesis} such that $\frac{G_0^\infty}{Z(G_0^\infty)}=B_l(q) $ with $l\geq 3$ and $q$ odd.  \begin{enumerate}[(i)]
    \item \label{bili1} If $\lambda$ is in Table \ref{bili}, and the value of $n$
 and a lower bound for $d$ are as given in the table. \begin{table}[]
 \begin{center}
\begin{tabular}{c|c|c|c|c|c}
    &$\lambda$ & $n$ &$d\geq$ &extra conditions\\\hline
   1.& $\omega_1$&$2l+1$&$=2$&$q=q_0$\\
    2.& $\omega_2$ &$2l^2+l$&$\frac{2l^2+l}{4l-2}$ & $q=q_0$\\
    3.& $2\omega_1$&$2l^2+3l-\epsilon_p(l+1)$&$\frac{2l^2+3l-1}{2l}$& $q=q_0$\\
    4.& $\omega_1+p^i\omega_1$&$(2l+1)^2$&$\frac{(2l+1)^2}{2l}$&$q=q_0$\\
     5.& $\omega_1+q_0\omega_1$&$(2l+1)^2$&$\frac{(2l+1)^2}{4l}$&$q=q_0^2$\\
     6.&$\omega_l$&$2^l$&$=2$&$l=3$ or $4$ and $q=q_0$ \\\hline
    
\end{tabular}
\end{center}
        \caption{Bounds in Theorem \ref{bl(q)}(\ref{bili1}) }
        \label{bili}
    \end{table}
    \item \label{part2bili} For $l\geq 14$ and $\lambda$ not in Table \ref{bili}, we have
    $$
        d\geq \frac{4l^3-l-3}{18l-30}.
    $$ For $l\leq 13$  and $\lambda$ not in Table \ref{bili}, the lower bound is as follows.
    \begin{align*}
     \begin{tabular}{c|c|c|c|c|c|c|c|c|c|c|c|c|c|c|c}
            l &13&12&11&10&9&8&7&6&5&4&3&2\\\hline
            $d\geq$&13&12&11&10 &6&4&3&3&3&3&3&3
        \end{tabular}
     \end{align*}
\end{enumerate}


\end{theorem}
The lower bounds in the statement of Theorem \ref{bl(q)} are greater than the lower bounds contained in Theorem \ref{lhalf}, so Theorem \ref{lhalf} holds for $X_l(q)=B_l(q)$.
     Moreover, for $n>(2l+1)^2,$  $$orbdiam(G,V)\geq \frac{l^2}{18}.$$ 

\begin{proof}[Proof of Theorem \ref{bl(q)}]

\underline{\textbf{Proof of the bounds in Table \ref{bili}}}

1. Here $\lambda=\omega_1$ and $V$ is the natural module so the result follows from Lemma \ref{naturalmod}.\par 
2-5. These cases are proved by Lemma \ref{orbit-stabparaborel} in the usual way. \par 
6. Here $(\lambda,n)=(\omega_l,2^l)$ and $l=3$ or $4.$\par
   For $l=4,$ as discussed already for the even characteristic case $G_0\triangleright C_4(2^r)\cong B_4(2^r)$, also in odd characteristic, the orbital diameter is 2 by Lemma \ref{smalldiamlemma}.\par 
For $l=3,$ if $G_0$ contains the scalars in $GL_n(q_0),$ then $G$ is a rank 3 group by \cite{liebeckaffine} so the orbital diameter is 2 by Lemma \ref{smalldiamlemma}.

From now on we assume that $\lambda$ is not in Table \ref{bili}. \par 
\underline{\textbf{Proof of the bounds in part (\ref{part2bili})}}\par 
Recall that either $q=q_0$ or $q=q_0^k$ and $V=V(\lambda')\otimes V(\lambda')^{q_0}\otimes\dots\otimes V(\lambda')^{q_0^{k-1}}$ for some $p$-restricted dominant weight $\lambda'$ by Lemma \ref{q0lemma}. 
\par \underline{\textbf{Case 1}: $q=q_0$}\par 
A Borel subgroup fixes a maximal $1$-space, so using Lemma  \ref{orbit-stabparaborel} we see that there is an orbit $\mathcal{O}$ of size $\vert \mathcal{O} \vert \leq{q^{l^2+l}}.$ \par 
 \underline{\textit{Case 1.a}: $\lambda$ is $p$-restricted and $l\geq 14$}\par 
 By Theorem \ref{alvarobounds}, if $n\leq 16{l \choose 4},$ then either $\lambda=\omega_3,$ $3\omega_1$ or $\omega_1+\omega_2.$ Using Lemma \ref{para} and Lemma  \ref{orbit-stabparaborel}, we find upper bounds for the orbit of the $1$-space fixed by the respective parabolic. Lemma \ref{lemma2.1} parts (\ref{boundformula1}) and(\ref{boundeq1}) give us the bounds 
\begin{align*}
\begin{tabular}{c|c|c|c|c}
    $\lambda$ & $\omega_3$ &$3\omega_1$ & $\omega_1+\omega_2$ \\\hline
    $d\geq$ & $\frac{4l^3-l-3}{18l-30}$& $\frac{ (3 + l) (-1 + 2 l) (1 + 2 l)}{6l}$& $\frac{  4 l^3+ 12 l^2 - 7 l-6}{12l-6}$\\\hline
\end{tabular}
\end{align*}
These satisfy the bound in part  (\ref{part2bili}).\par 
 Suppose $n\geq 16{l \choose 4}.$ A Borel subgroup fixes a maximal $1$-space and so Lemma \ref{lemma2.1}(\ref{boundeq1}), which gives that $$d\geq \frac{16{l \choose 4}-1}{l^2+l}.$$ This satisfies the bound in part  (\ref{part2bili}).\par 
\underline{\textit{Case 1.b}: $\lambda$ is $p$-restricted and $3\leq l\leq 13$}\par 
We can see from Theorems \ref{alvarobounds} and \ref{lowerlprestrict} that either $(\lambda,l)=(\omega_6,6),$ $(\omega_5,5)$
or the lower bounds for $n$ and $d$ in (\ref{bmdiam}) hold. The bounds on $d$ in (\ref{bmdiam}) are obtained using the fact that a Borel fixes a maximal $1$-space.
\begin{align}\label{bmdiam}
    \begin{tabular}{c|c|c|c|c|c|c|c|c|c|c|c}
    $l$& $3$&$4$&$5$&$6$&$7$&$8$&$9$&$10$&$11$&$12$&$13$\\ \hline
       $n\geq $  & $27$&$64$&$100$&$208$&$128$&$256$&$512$&$1000$&$1331$&$1728$&$2197$ \\\hline
       $d\geq $  &$3$&$4$&$4$&$3$&$3$&$4$&$6$&$10$&$11$&$12$&$13$\\\hline
    \end{tabular}
\end{align}
 
\underline{$(\lambda,l)=(\omega_6,6)$} \par  The parabolic $P_l$ fixes a maximal $1$-space so using Lemma  \ref{orbit-stabparaborel} there is an orbit of size at most $\vert \mathcal{O} \vert \leq \frac{q^{23}}{2}.$ By Lemma \ref{lemma2.1}(\ref{boundeq2}) it follows that $d\geq \frac{64}{23}$ and so $d\geq 3.$  \par
\underline{$(\lambda,l)=(\omega_5,5)$} \par 
Here $B_5(q)\leq D_6(q)\leq GL_{32}(q)$ and since if $G_0\triangleright D_6(q)$ the orbital diameter is at least 3 by Claim \ref{d6},  the orbital diameter this case is also at least 3 by Lemma \ref{subgroup}. \par 
\underline{\textit{Case 1.c: $\lambda$ is not $p$-restricted}}\par 
By \cite[Thm 5.4.5]{kleidmanliebeck}, and Theorems \ref{alvarobounds} and \ref{lowerlprestrict}, we have that either $n\geq l^4$ or $\lambda=\omega_1+p^i \omega_1+p^j\omega_1,$ $\omega_1+p^i\omega_2,$  $\omega_1+p^i2\omega_1$ or for $3\leq l\leq 7,$ $\omega_1+p^i\omega_l.$   \par 
For $n\geq l^4$ we use the fact that a Borel fixes a maximal $1$-space to get that $d\geq \frac{l^4}{l^2+l}$ and so the bound in part (\ref{part2bili}) is satisfied.\par
Using Lemma \ref{para} and Lemma  \ref{orbit-stabparaborel}, we find upper bounds for the orbit of the $1$-space fixed by the respective parabolic. Lemma \ref{lemma2.1} parts (\ref{boundformula1}) and(\ref{boundeq1}) give us the bounds 
\begin{align}
\begin{tabular}{c|c|c|c|c}
    $\lambda$ & $\omega_1+p^i \omega_1+p^j\omega_1$& $\omega_1+p^i\omega_2$&  $\omega_1+p^i2\omega_1$ & 
    $\omega_1+p^i\omega_l$
    \\\hline
    $d\geq$ & $\frac{(2l+1)^3}{2l}$&$\frac{(2l+1)(2l^2+l)}{4l-2}$&$\frac{(2l+1)(2l^3+3l-1}{2l}$ &$\frac{2^{l+1}-2}{l(l-1)}$\\\hline
\end{tabular}
    \end{align}
These satisfy the bound in part  (\ref{part2bili}). \par 
 \underline{\textbf{Case 2: $q=q_0^k$ for $k\geq 2$ as in Lemma \ref{q0lemma}(\ref{subfield})}}\par 
 Here $V(\lambda)=V(\lambda')\otimes\dots\otimes V(\lambda')^{q_0^{k-1}}.$ By Theorems \ref{alvarobounds} and \ref{lowerlprestrict}, either $\lambda'=\omega_1$ or $\dim V(\lambda')\geq l^2-1.$ 
We will consider these two cases in turn.\par 
\underline{$n=(2l+1)^k$ and $\lambda=\omega_1+q_0\omega_1+\dots q_0^{k-1}\omega_1.$}\par 
Here $k\geq 3,$ as the $k=2$ case is in Table \ref{bili}. \par 
In this case the parabolic subgroup $P_1$ fixes a $1$-space and so by Lemma  \ref{orbit-stabparaborel} $\vert\mathcal{O}\vert \leq \frac{q_0^{2kl}}{2}$ and so using Lemma \ref{lemma2.1}(\ref{boundeq2}) we can see that $d\geq \frac{(2l+1)^k}{2kl}\geq \frac{(2l+1)^3}{6l}.$ This satisfies the bound in part  (\ref{part2bili}). \par 
\underline{$\lambda\neq\omega_1+q_0\omega_1+\dots q_0^{k-1}\omega_1$}\par Now assume $\lambda'\neq \omega_1.$
Here $n> (2l+1)^k,$ and so by Theorems \ref{alvarobounds} and \ref{lowerlprestrict}, $n\geq (l^2-1)^k.$ Since a Borel is fixing a maximal $1$-space, $\vert\mathcal{O}\vert \leq \frac{q_0^{k(l^2+l)}}{2}$ and so by Lemma \ref{lemma2.1}(\ref{boundeq2}) $d\geq \frac{(l^2-1)^k}{k(l^2+l)}\geq \frac{(l^2-1)^2}{2l^2+2l}.$ This satisfies the bound in part  (\ref{part2bili}) for $l\geq 3.$ 

\end{proof}
From Theorem \ref{bl(q)} and Lemma \ref{smalldiamlemma}, we can achieve the following classification.
 \begin{prop}\label{bldiam2}
Let $G$ be as in Hypothesis \ref{hypothesis} such that $\frac{G_0^\infty}{Z(G_0^\infty)}=B_l(q).$ Then $orbdiam(G,V)\leq 2$ if and only if one of the following holds.
\begin{itemize}
    \item $V$ is the natural module.
     \item $G_0\triangleright B_3(q)$ and $V=V_{8}(q).$
     \item $G_0\triangleright B_4(q)$ and $V=V_{16}(q).$
\end{itemize}
 \end{prop}

\subsubsection{$G_0\triangleright D_l(q)$}
Note that for $l\leq 3,$ $D_l(q)$ is isomorphic to other classical groups and have already been considered.
\begin{theorem}\label{dl(q)}
Let $G$ be as in Hypothesis \ref{hypothesis} such that $\frac{G_0^\infty}{Z(G_0^\infty)}=D_l(q)$ with $l\geq 4.$ 
\begin{enumerate}[(i)]
    \item \label{dli1} If $\lambda$ is in Table \ref{dli}, and  the value of $n$
 and a lower bound for $d$ are as given in the table.\begin{table}[]
         \begin{center}
\begin{tabular}{c|c|c|c|c|}
    &$\lambda$ & $n$ &$d\geq$&extra conditions \\\hline
    1.&$\omega_1$&$2l$&$=2$&$q=q_0$\\
    2.& $\omega_2$ &$2l^2-l-\gamma$&$\frac{2l^2-l-3}{4l-5}$&$\gamma=gcd(2,l)$ if $p=2$ and $\gamma=0$ otherwise and $q=q_0$ \\
     3.& $2 \omega_1$&$2l^2+l-1-\epsilon_p(l)$&$\frac{2l^2+l-3}{2l}$&$q=q_0$\\
   4.&  $\omega_1+p^i\omega_1$&$4l^2$&$\frac{4l^2-1}{2l}$&$q=q_0$\\
5.&$\omega_1+q_0\omega_1$&$4l^2$&$l$& $q=q_0^2$\\
6.&$\omega_l$&$2^l$&$=2$&$l=5$ and $q=q_0$\\\hline
\end{tabular}
\end{center}
\caption{Bounds in Theorem \ref{dl(q)}(\ref{dli1})}
        \label{dli}
    \end{table}    
    \item \label{partedli}  For $l\geq 16$ and $\lambda$ not in Table \ref{dli}, we have

        $$d \geq \frac{4l^3-6l^2-4l}{18l-39}$$
    
    For $l\leq 15$ and $\lambda$ not in Table \ref{dli}, a lower bound is as follows.
   \begin{align*}
     \begin{tabular}{c|c|c|c|c|c|c|c|c|c|c|c|c|c|c|c}
            l &15&14&13&12&11&10&9&8&7&6&5&4\\\hline
            $d\geq$&15&14&13&12 &16&10&6&4&3&3&2&3
        \end{tabular}
      \end{align*}
\end{enumerate}

\end{theorem}
The lower bounds in the statement of Theorem \ref{dl(q)} are greater than the lower bounds contained in Theorem \ref{lhalf}, so Theorem \ref{lhalf} holds for $X_l(q)=D_l(q)$.
\begin{proof}[Proof of Theorem \ref{dl(q)}]    

\underline{\textbf{Proof of the bounds in Table \ref{dli}}}

1. Here $\lambda=\omega_1$ and $V$ is the natural module so the result follows from Lemma \ref{naturalmod}. \par
2-5. These cases are proved by Lemma \ref{orbit-stabparaborel} in the usual way. \par 
6. Here $(\lambda,n,l)=(\omega_5,16,5).$
    \par $G$ is a rank three group by \cite{liebeckaffine} when $G_0$ contains the scalars in $GL_n(q_0)$. The orbital diameter is 2 by Lemma \ref{smalldiamlemma}. \par

From now on we will assume that $\lambda$ is not in Table \ref{dli}. \par 
 \underline{\textbf{Proof of the bounds in part (\ref{partedli})}} \par
Recall that either $q=q_0$ or $q=q_0^k$ with $k\geq 2$ and $V=V(\lambda')\otimes V(\lambda')^{q_0}\otimes\dots\otimes V(\lambda')^{q_0^{k-1}}$ for some $p$-restricted dominant weight $\lambda'$ by Lemma \ref{q0lemma}. 
\par \underline{\textbf{Case 1}: $q=q_0$}\par 

A Borel subgroup fixes a maximal $1$-space, so using Lemma  \ref{orbit-stabparaborel} we see that there is an orbit $\mathcal{O}$ of size $\vert \mathcal{O} \vert \leq{q^{l^2}}.$ \par 
 \underline{\textit{Case 1.a}: $\lambda$ is $p$-restricted and $l\geq 16$}\par 
 By Theorem \ref{alvarobounds}, if $n < 16{l \choose 4},$ then either  $\lambda=\omega_3,$ $3\omega_1$ or $\omega_1+\omega_2.$ Using Lemma \ref{para} and Lemma  \ref{orbit-stabparaborel}, we find upper bounds for the orbit of the $1$-space fixed by the respective parabolic. Lemma \ref{lemma2.1} parts (\ref{boundformula1}) and(\ref{boundeq1}) give us the bounds 
\begin{align*}
\begin{tabular}{c|c|c|c|c}
    $\lambda$ & $\omega_3$ &$3\omega_1$ & $\omega_1+\omega_2$ \\ \hline
    $d\geq$ & $\frac{4l^3-6l^2-4l}{18l-39}$&$\frac{4l^3+6l^3-10l-3}{6l}$&$\frac{4l^3+6l^2-16l-3}{12l-12}$\\\hline
\end{tabular}
\end{align*}
These satisfy the bound in part  (\ref{partedli}).\par 
 Suppose $n\geq 16{l \choose 4}.$ A Borel subgroup fixes a maximal $1$-space and so Lemma \ref{lemma2.1}(\ref{boundeq1}), which gives that $$d\geq \frac{16{l \choose 4}-1}{l^2}.$$ This satisfies the bound in part  (\ref{partedli}).\par

 \underline{\textit{Case 1.b}: $\lambda$ is $p$-restricted and $4\leq l\leq 15$}\par  

Now we need to consider cases that are not in Table \ref{dli} for $4\leq l\leq 15.$ By Theorem \ref{lowerlprestrict} either $\lambda=\omega_l$ or $(\lambda,l)=(\omega_l, 4\leq l\leq 11),$ $(\omega_3,5),$ $(\omega_3,6),$ $(\omega_3,7),$ $(\omega_1+\omega_3,4)$ or $n\geq l^3.$ We prove the result for these in turn. \par 
\underline{$\lambda=\omega_l$ and $4\leq l \leq 11$}\par 
Note that for $l\geq 12,$ $2^{l-1}\geq l^3$ so we only need to consider $4\leq l \leq 11.$
For $l=4$ this is the natural module, and for $l=5$ this is in Table \ref{dli} in part (\ref{dli1}). For $l=6,$ the bound is $d\geq 3$ by Claim \ref{d6}.\par 
The parabolic $P_l$ fixes a maximal $1$-space so by Lemma  \ref{orbit-stabparaborel} we have an orbit of size $\vert \mathcal{O} \vert \leq q^{\frac{l(l+1)}{2}}$ and so using Lemma \ref{lemma2.1}(\ref{boundeq1}) we get the bounds. 
\begin{align*}
\begin{tabular}{c|c|c|c|c|c}
    $l$&11 & 10&9&8&7 \\\hline
   $d\geq $& 16 &10&6&4&3 \\\hline
\end{tabular}
\end{align*}\par 
\underline{$\lambda=\omega_{3}$ and $5\leq l\leq 7$}\par  
The parabolic $P_{3}$ fixes a maximal $1$-space so by Lemma  \ref{orbit-stabparaborel} we have an orbit of size $\vert \mathcal{O} \vert\leq  q^{7l-16}$ and so using Lemma \ref{lemma2.1}(\ref{boundeq1}) we get the bounds for $d$ as in the following Table.
\begin{align*}
\begin{tabular}{c|c|c|c}
    $l$&7&6&5 \\\hline
   $d\geq $& 11&8&6\\\hline
\end{tabular}
\end{align*}\par 
\underline{$(\lambda,l)=(\omega_1+\omega_3,4)$}\par 
The parabolic $P_{1,3}$ fixes a maximal $1$-space, so it follows that there is an orbit such that $\vert\mathcal{O}\vert \leq \frac{(q^4-1)^2(q^3+1)}{(q-1)}.$ Now, since $n\geq 48,$ Lemma \ref{lemma2.1}(\ref{boundformula1}) shows that $d\geq 3.$\par 
\par 
\underline{$n\geq l^3$}\par 
In the case when $n\geq l^3$, a Borel subgroup fixes a maximal $1$-space, so using Lemma  \ref{orbit-stabparaborel} we see that there is an orbit $\mathcal{O}$ of size $\vert \mathcal{O} \vert \leq{q^{l^2}},$ so we can use Lemma \ref{lemma2.1}(\ref{boundeq1}). This gives that $d\geq \frac{l^3-1}{l^2}$ and satisfies the bound in part (\ref{partedli}).\par 
 \underline{\textit{Case 1.c}: $\lambda$ is not $p$-restricted}\par 
By Theorems \ref{alvarobounds} and \ref{lowerlprestrict} we have that either $n\geq l^4$ or $\lambda=\omega_1+p^i \omega_1+p^j\omega_1,$ $\omega_1+p^i\omega_2,$  $\omega_1+p^i2\omega_1$ or $\omega_1+p^i\omega_l$ with $5\leq l\leq 7$. \par 
For $n\geq l^4$ we use the fact that a Borel fixes a maximal $1$-space to get that $d\geq l^2$ and so the bound in part (\ref{partedli})  is satisfied.\par 
Using Lemma \ref{para} and Lemma  \ref{orbit-stabparaborel}, we find upper bounds for the orbit of the $1$-space fixed by the respective parabolic. Lemma \ref{lemma2.1} parts (\ref{boundformula1}) and(\ref{boundeq1}) give us the bounds 
\begin{align}
\begin{tabular}{c|c|c|c}
    $\lambda$ & $\omega_1+p^i \omega_1+p^j\omega_1$& $\omega_1+p^i\omega_2$&  $\omega_1+p^i2\omega_1$  \\\hline
    $d\geq$ & $4l^2$&$\frac{(2l)(2l^2-l-2)}{4l-4}$&$2l^2-l-2$\\\hline
\end{tabular}
    \end{align}
These satisfy the bound in part  (\ref{partedli}). \par 
For $\lambda=\omega_1+p^i\omega_l$ and $5\leq l\leq 7$ we use the fact that the parabolic $P_{1,l}$ fixes a maximal $1$-space. Using Lemma  \ref{orbit-stabparaborel}, we find upper bounds for the orbit of the maximal vector. Lemma \ref{lemma2.1} parts (\ref{boundeq1}) we get the following bounds. 
\[\begin{array}{c|c|c|c}
   l  & 5&6&7 \\\hline
   d\geq  & 10&17&30.\\\hline
\end{array}\]
 \underline{\textbf{Case 2: $q=q_0^k$ for $k\geq 2$ as in Lemma \ref{q0lemma}(\ref{subfield})}}\par 
Here $V(\lambda)=V(\lambda')\otimes \dots \otimes V(\lambda')^{q_0^{k-1}}.$ By Theorems \ref{alvarobounds} and \ref{lowerlprestrict}, either $\lambda'=\omega_1$ or $\dim(V(\lambda'))\geq l^2-9.$ \par 
\underline{$n=(2l)^k$ and $\lambda=\omega_1+q_0\omega_1+\dots q_0^{k-1}\omega_1.$}\par 
Here $k\geq 3,$ as the $k=2$ case is in Table \ref{dli} in part (\ref{dli1}).
In this case the parabolic subgroup $P_1$ fixes a $1$-space and so by Lemma  \ref{orbit-stabparaborel}, $\vert\mathcal{O}\vert \leq \frac{q_0^{2kl}}{2}$ and so using Lemma \ref{lemma2.1}(\ref{boundeq2}) we can see that $d\geq \frac{(2l)^k}{2kl}\geq \frac{(2l)^3}{6l}.$ This satisfies the bound in part  (\ref{partedli}). \par 
\underline{$\lambda\neq \omega_1+q_0\omega_1+\dots q_0^{k-1}\omega_1.$} \par 
Now assume $\lambda'\neq \omega_1.$
By Theorems \ref{alvarobounds} and \ref{lowerlprestrict}, $n\geq (l^2-9)^k$ in this case.  
Since a Borel is fixing a maximal $1$-space, $\vert\mathcal{O}\vert \leq \frac{q_0^{kl^2}}{2}$ and so by Lemma \ref{lemma2.1}(\ref{boundeq2}) $d\geq \frac{(l^2-9)^k}{k(l^2)}\geq \frac{(l^2-9)^2}{2l^2}$ for $l\geq 5.$ If $l=4,$ then $n\geq l^2$ and so $d\geq \frac{(l^2)^k}{k(l^2)}\geq \frac{(l^2)^2}{2l^2}$ and $d\geq 3.$ Hence the bound in part (\ref{partedli}) holds.
\end{proof}
The following classification is immediate.
 \begin{prop}\label{dldiam2}
Let $G$ be as in Hypothesis \ref{hypothesis} such that $\frac{G_0^\infty}{Z(G_0^\infty)}=D_l(q).$ Then $orbdiam(G,V)\leq 2$ if and only if one of the following holds.
\begin{itemize}
    \item $V$ is the natural module.
     \item $G_0\triangleright D_5(q)$ and $V=V_{16}(q).$
\end{itemize}
 \end{prop}
\subsubsection{$G_0\triangleright \,^2\!D_l(q)$}
We will consider the case when $l\geq 4,$ as the lower rank cases are isomorphic to other classical groups and have already been considered.\par 
\begin{theorem}\label{2dl(q)}
Let $G$ be as in Hypothesis \ref{hypothesis} such that $\frac{G_0^\infty}{Z(G_0^\infty)}=\,^2\!D_l(q)$ with $l\geq 4.$ \par
\begin{enumerate}
\item If $\tau_0(\lambda)\neq \lambda,$ then $^2D_l(q)\leq D_l(q^2)\leq GL(V) $ so the bounds from Theorem \ref{dl(q)} hold. \par
\item Suppose $\tau_0(\lambda)=\lambda.$  \begin{enumerate}[(i)]

    \item \label{makosnudli} If $\lambda$ is as in Table \ref{nudli}, and   the value of $n$
 and a lower bound for $d$ are as given in the table.
    \begin{table}[]

    \begin{center}
\begin{tabular}{c|c|c|c|c|}
 & $\lambda$ & $n$ &$d\geq$& extra conditions \\ \hline
 1.  & $\omega_1$&$2l$&$=2$&$q=q_0$\\
   2. & $\omega_2$ &$2l^2-l-\gamma$&$\frac{2l^2-l-3}{4l-4}$&$\gamma=gcd(2,l)$ if $p=2$ and $\gamma=0$ otherwise and  $q=q_0$ \\
    3.  &  $2\omega_1$&$2l^2+l-1-\epsilon_p(l)$&$\frac{2l^2-l-2}{2l+1}$&$q=q_0$\\
     4.  &  $\omega_1+p^i\omega_1$&$(2l)^2$&$\frac{4l^2-1}{2l+1}$&$q=q_0$\\
 5.  &  $\omega_1+q_0\omega_1$&$(2l)^2$&$\frac{4l^2-1}{4l+1}$&$q=q_0^2$\\\hline

\end{tabular}
\end{center}

        \caption{Bound in Theorem \ref{2dl(q)} part \ref{makosnudli}}
        \label{nudli}
    \end{table}
    \item \label{nudli2}  For $l\geq 16$ and $\lambda$ not in Table \ref{nudli}, we have
$$d\geq \frac{4l^3-6l^2-4l}{18l-39}.$$
  For $l\leq 15$  and $\lambda$ not in Table \ref{nudli}, the lower bound is as follows.
      \begin{align*}
     \begin{tabular}{c|c|c|c|c|c|c|c|c|c|c|c|c|c|c|c}
            l &15&14&13&12&11&10&9&8&7&6&5&4\\\hline
            $d\geq$&15&14&13&12 &11&10&9&8&7&6&5&4
        \end{tabular}
      \end{align*}

\end{enumerate}   
\end{enumerate}

\end{theorem}
The lower bounds in the statement of Theorem \ref{2dl(q)} are greater than the lower bounds contained in Theorem \ref{lhalf}, so Theorem \ref{lhalf} holds for $X_l(q)=^2D_l(q)$.

\begin{proof}[Proof of Theorem \ref{2dl(q)}]

\underline{\textbf{Case I}:$\tau_0(\lambda)\neq\lambda$}\par 
By Lemma \ref{q0lemma}(\ref{twistedexceptional}), $^2D_l(q)\leq D_l(q^2)\leq GL(V)$ and so the lower bounds for the orbital diameter from the case of $G_0\triangleright D_l(q)$ hold. \par  
\underline{\textbf{Case II: $\tau_0(\lambda)=\lambda$}} \par
\underline{\textbf{Proof of the bounds in Table \ref{nudli}}}\par 

1. Here $(\lambda,n)=(\omega_1,2l).$ \par In this case $V$
is the natural module so by Lemma \ref{naturalmod} the orbital diameter is 2.\par 
2-5. These cases are proved by Lemma \ref{orbit-stabparaborel} in the usual way. \par 

 \par
 From now we assume that $\lambda$ is not in Table \ref{nudli}. \par 
\underline{\textbf{Proof of the bounds in (\ref{nudli2})}}\par Recall that either $q=q_0$ or $q=q_0^k$ and $V=V(\lambda')\otimes V(\lambda')^{q_0}\otimes\dots\otimes V(\lambda')^{q_0^{k-1}}$ for some $p$-restricted dominant weight $\lambda'$ by Lemma \ref{q0lemma}. 
\par 
\underline{\textbf{Case 1}: $q=q_0$}\par 
A Borel subgroup fixes a maximal $1$-space, so using Lemma  \ref{orbit-stabparaborel} we see that there is an orbit $\mathcal{O}$ of size 
$\vert \mathcal{O} \vert\leq\frac{q^{l^2+2}}{2}.$ \par
\underline{\textit{Case 1.a: $\lambda$ is $p$-restricted and $l\geq 16$}}\par
By Theorem \ref{alvarobounds}, if $n< 16{l \choose 4},$ then either $\lambda=\omega_3,$ $3\omega_1$ or $\omega_1+\omega_2.$ Using Lemma \ref{para} and Lemma  \ref{orbit-stabparaborel}, we find upper bounds for the orbit of the $1$-space fixed by the respective parabolic. Lemma \ref{lemma2.1} parts (\ref{boundformula1}) and(\ref{boundeq1}) give us the bounds 
\begin{align*}
\begin{tabular}{c|c|c|c}
    $\lambda$ & $\omega_3$ &$3\omega_1$ & $\omega_1+\omega_2$ \\\hline
    $d\geq$ & $\frac{4l^3-6l^2-4l}{18l-39}$&$\frac{4l^3+6l^3-10l-3}{6l+3}$&$\frac{4l^3+6l^2-16l-3}{12l-12}$\\\hline
\end{tabular}
\end{align*}
These satisfy the bound in part  (\ref{nudli2}).\par 
 Suppose $n\geq 16{l \choose 4}.$ a Borel subgroup fixes a maximal $1$-space and so Lemma \ref{lemma2.1}(\ref{boundeq1}), which gives that $$d\geq \frac{16{l \choose 4}-1}{l^2+2}.$$ This satisfies the bound in part  (\ref{nudli2}).\par

\underline{\textit{Case 1.b: $\lambda$ is $p$-restricted and $4\leq l\leq 15$}}\par 

 By Theorem \ref{lowerlprestrict}, either $n\geq l^3$ or $(\lambda,l)=(\omega_3,5),$ $(\omega_3,6),$ $(\omega_3,7)$ or $(\omega_1+\omega_2,4).$ \par 
Suppose $n\geq l^3.$ a Borel subgroup fixes a maximal $1$-space, so we can use Lemma \ref{lemma2.1}(\ref{boundeq2}). This gives that $$d\geq \frac{l^3}{l^2+2}$$ and so the bound in part (\ref{nudli2}) is satisfied. \par 
Suppose $\lambda=\omega_3$ and $l=5,$ $6$ or $7.$ The parabolic $P_{3}$ fixes a $1$-space. Using Lemma  \ref{orbit-stabparaborel}, we can bound the size of an orbit and using Lemma \ref{lemma2.1}(\ref{boundeq1}) find a bound for the diameter, as in the following Table. \par 
\begin{tabular}{c|c|c|c}
$l$&7&6&5\\\hline
    $\vert \mathcal{O} \vert \leq $ &$q^{29}$&$q^{23}$& $q^{17}$ \\\hline
    $n\geq$ & $336$ &$208$&$100$\\\hline
    $d\geq $&$12$&$9$& $6$\\\hline
\end{tabular} \par 
Suppose $(\omega_1+\omega_2,4).$ Then $n=48$ or $56.$ In this case the parabolic $P_{3,4}$ fixes a $1$-space so there is an orbit $\vert\mathcal{O}\vert \leq q^{12}$ and so using Lemma \ref{lemma2.1}(\ref{boundeq1}) it follows that $d\geq 4.$ 
\par 
\underline{\textit{Case 1.c: $\lambda$ is not $p$-restricted}}\par 
By Theorems \ref{alvarobounds} and \ref{lowerlprestrict} we have that either $n\geq l^4$ or $\lambda=\omega_1+p^i \omega_1+p^j\omega_1,$ $\omega_1+p^i\omega_2,$  $\omega_1+p^i2\omega_1$. \par 
For $n\geq l^4$ we use the fact that a Borel fixes a maximal $1$-space to get that $d\geq \frac{l^4-1}{l^2+2}$ and so the bound in part (\ref{nudli2}) is satisfied.\par 
Using Lemma \ref{para} and Lemma  \ref{orbit-stabparaborel}, we find upper bounds for the orbit of the $1$-space fixed by the respective parabolic. Lemma \ref{lemma2.1} parts (\ref{boundformula1}) and(\ref{boundeq1}) give us the bounds 
\begin{align}
\begin{tabular}{c|c|c|c}
    $\lambda$ & $\omega_1+p^i \omega_1+p^j\omega_1$& $\omega_1+p^i\omega_2$&  $\omega_1+p^i2\omega_1$  \\\hline
    $d\geq$ & $\frac{(2l)^3}{2l+1}$&$\frac{(2l)(2l^2-l-2)}{4l-4}$&$\frac{(2l)(2l^2-l-2)}{2l+1}$\\ \hline
\end{tabular}
    \end{align}
These satisfy the bound in part  (\ref{nudli2}). \par

 \underline{\textbf{Case 2: $q=q_0^k$ for $k\geq 2$ as in Lemma \ref{q0lemma}(\ref{subfield})}}\par 

Here $V(\lambda)=V(\lambda')\otimes \dots \otimes V(\lambda')^{q_0^{k-1}}.$ By Theorems \ref{alvarobounds} and \ref{lowerlprestrict}, either $\lambda'=\omega_1$ or  $\dim(V(\lambda'))\geq l^2.$  Note that this number, as for the $C_l(q)$ case comes from examining the lowest dimensions for each $l$ using Theorems \ref{alvarobounds} and \ref{lowerlprestrict}.\par 

\underline{$n=(2l)^k$ and $\lambda=\omega_1+q_0\omega_1+\dots q_0^{k-1}\omega_1.$}\par 
Here $k\geq 3,$ as the $k=2$ case is in Table \ref{nudli} in part (\ref{makosnudli}).
In this case the parabolic subgroup $P_1$ fixes a $1$-space and so by Lemma  \ref{orbit-stabparaborel} $\vert\mathcal{O}\vert \leq {q_0^{2kl+1}}$ and so using Lemma \ref{lemma2.1}(\ref{boundeq1}) we can see that $d\geq \frac{(2l)^k-1}{2kl+1}\geq \frac{(2l)^3-1}{6l+1}.$ This satisfies the bound in part  (\ref{nudli2}). \par 
\underline{$\lambda\neq \omega_1+q_0\omega_1+\dots q_0^{k-1}\omega_1.$}
Now assume $\lambda'\neq \omega_1.$
By Theorems \ref{alvarobounds} and \ref{lowerlprestrict}, $n\geq (l^2)^k$ in this case.  
Since a Borel is fixing a maximal $1$-space so $\vert\mathcal{O}\vert \leq {q_0^{kl^2+1}}$ and so by Lemma \ref{lemma2.1}(\ref{boundeq1}) $d\geq \frac{(l^2)^k-1}{k(l^2)+1}\geq \frac{l^4-1}{2l^2+1}$ satisfying part (\ref{nudli2}). 
\end{proof}
The following classification is immediate.
 \begin{prop}\label{2dldiam2}
Let $G$ be as in Hypothesis \ref{hypothesis} such that $\frac{G_0^\infty}{Z(G_0^\infty)}=^2D_l(q).$ If $orbdiam(G,V)\leq 2$  then one of the following holds.
\begin{itemize}
    \item $V$ is the natural module.
     \item  $G_0\triangleright \,^2\!D_5(q)$ and $V=V_{16}(q^2).$
\end{itemize}
 \end{prop}
 \underline{Note} We have not been able to determine the orbital diameter in the second case.
\subsection{Exceptional Stabilizers}
In this section we will prove Theorem \ref{lhalf} for the case when $\frac{G_0^\infty}{Z(G_0^\infty)}$ is a simple group of exceptional Lie type. In fact we prove the following stronger result. 
\begin{theorem}\label{exceptionalcases}
Let $G$ be as in Hypothesis \ref{hypothesis} such that $\frac{G_0^\infty}{Z(G_0^\infty)}\cong X_l(q)$ an exceptional group. Denote the minimal module of $G_0$ by $V_{min}$  and the adjoint module of $G_0$ by $V_{ad}.$  A lower bound for $d$ is as in Table \ref{tab:exceptional}. 
\begin{table}[]

  \[ \begin{array}{c|c|c|c}
         & V_{min}& V_{ad}& \text{rest} \\\hline
       E_8  & 4&4&29\\\hline
      E_7  & 3&4&13\\\hline
      E_6  & 3&4&8\\\hline
      ^2E_6  & 3&4&8\\\hline
      F_4  & 3&3&7\\\hline
      ^2F_4  & 3&3&7\\\hline
      G_2  & 2-\epsilon_p(2)&3&3\\\hline
      ^2G_2  & 3&3&3\\\hline
       ^2B_2  & 2&3&3\\\hline
        ^3D_4  & 2&3&3\\\hline
    \end{array}\] 
    \caption{Lower bounds for the orbital diameter with exceptional stabilizer}
    \label{tab:exceptional}
\end{table}
\end{theorem}
Now we prove this result for each exceptional group in turn. Recall that in all of these cases, $q_0$ is as in Lemma \ref{q0lemma}. \par 
\subsubsection{$G_0\triangleright E_8(q)$}
\par \underline{\textbf{Case 1}: $q=q_0$}\par 
A Borel subgroup fixes a maximal $1$-space, so using Lemma  \ref{orbit-stabparaborel}  we see that there is an orbit $\mathcal{O}$ of size $\vert \mathcal{O} \vert \leq\frac{q^{129}}{2},$ so we can use Lemma \ref{lemma2.1}(\ref{boundeq2}). \par 
\underline{\textit{Case 1.a: $\lambda$ is $p$-restricted}}\par 
By Theorem \ref{lowerlprestrict}, either $n=248$ and $\lambda=\omega_8$ or $n\geq 3626.$ \par 
For $n\geq 3626$ we use the fact that a Borel subgroup fixes a maximal $1$-space and Lemma \ref{lemma2.1}(\ref{boundeq2}) gives us that $d\geq 29.$ \par 
In the case when $(\lambda,n)=(\omega_8,248),$ $V$ is simultaneously the adjoint and the minimal module. By Lemma \ref{para}, the parabolic $P_8$ fixes a $1$-space. By Lemma  \ref{orbit-stabparaborel}, $E_8(q)$ has an orbit of size at most $q^{65},$ so using Lemma \ref{lemma2.1}(\ref{boundeq1}) it follows that $d\geq 4.$ \par
\underline{Case 1.b: $\lambda$ is not $p$-restricted} \par 
Using  Theorem \ref{lowerlprestrict}, in this case $n\geq 248^2$ and so using the fact that a Borel fixes a maximal $1$-space, we can conclude using Lemmas \ref{orbit-stabparaborel} and \ref{lemma2.1} that $d\geq 477.$ \par 
 \underline{\textbf{Case 2: $q=q_0^k$ for $k\geq 2$ as in Lemma \ref{q0lemma}(\ref{subfield})}}\par
We use the fact that a Borel fixes a maximal $1$-space so it follows that $\vert\mathcal{O}\vert \leq q_0^{120k+8}$ and so by Lemma \ref{lemma2.1}(\ref{boundeq1}) it follows that $d\geq \frac{248^k-1}{120k+8}\geq\frac{248^2-1}{248} \geq 248.$

\subsubsection{$G_0\triangleright E_7(q)$}

\par \underline{\textbf{Case 1: $q=q_0$}}\par 
A Borel subgroup fixes a maximal $1$-space, so using Lemma  \ref{orbit-stabparaborel} we see that there is an orbit $\mathcal{O}$ of size $\vert \mathcal{O} \vert \leq\frac{q^{71}}{2},$ so we can use Lemma \ref{lemma2.1}(\ref{boundeq2}). \par 
\underline{\textit{Case 1.a: $\lambda$ is $p$–restricted}} \par
By Theorem \ref{lowerlprestrict}, either  $n=56,$ $133-\epsilon_p(2),$ or $n\geq 856.$ \par 
In the case when $n\geq 856,$ we use the fact that a Borel fixes a maximal $1$-space and by Lemma \ref{lemma2.1}(\ref{boundeq2}), $d\geq 13.$\par 
The case when $n=133-\epsilon_p(2),$ and $\lambda=\omega_1,$ $V$ is the adjoint module. By Lemma \ref{para}, it follows that the parabolic $P_1$ fixes a maximal $1$-space so using Lemma  \ref{orbit-stabparaborel}, $\vert \mathcal{O} \vert \leq q^{35}$ and so by Lemma \ref{lemma2.1}(\ref{boundeq1}), $d\geq 4.$\par 
In the case of $n=56,$ $V$ is the minimal module.  The stabilizers of the vectors are classified in \cite[Lemma 4.3]{e7}.
Consider the algebraic group $\overline{G_0}=E_7(K)$ where $K=\overline{\mathbb{F}_p}$ acting on $\overline{V}=V_{56}(K).$ Now the stabilizer of a maximal 1-space in $G_0$ is $P_7.$ Let $\Delta_0$ be the orbit of $\overline{G_0}$ containing the maximal vectors. If $v=a+b,$ where $a,b\in \Delta_0,$ then $\overline{G_0}_a\cap \overline{G_0}_b\leq \overline{G_0}_v.$ The stabilizers of $a$ and $b$ are conjugates of $P_7'.$ Without loss of generality, they are $P_7'$ and $P_7'^g$ for some $g\in \overline{G_0}.$ 


Now consider the finite group $G_0=\overline{G_0}^{(q)}=E_7(q)$ acting on $V=\overline{V}^{(q)}=V_{56}(q).$ By \cite[Lemma 4.3]{e7},  $^2E_6(q).2,$ stabilizes a vector $w\in V.$  By Lemma \ref{alunalemma}, the possible intersections of $P_7'\cap P_7'^g$ with $\overline{G_0}^{(q)}$ either contain a unipotent subgroup of order at least $q^{42},$ or contain either $D_5(q)$ or $E_6(q).$ 
Since none of these is contained in $^2E_6(q).2,$ we cannot express $w$ as a sum of at most two elements in $\Delta_0\cap V,$ so the orbital diameter is at least 3. \par 
\underline{Case 1.b: $\lambda$ is not $p$-restricted }\par 
In this case $n\geq 56^2$ and so using that fact that a Borel fixes a maximal $1$-space, $d\geq 45.$ \par 
 \underline{\textbf{Case 2: $q=q_0^k$ for $k\geq 2$ as in Lemma \ref{q0lemma}(\ref{subfield})}}\par
Using the fact that a Borel fixes a maximal $1$-space we have that there is an orbit $\vert\mathcal{O}\vert \leq q_0^{63k+7}$ and so by Lemma \ref{lemma2.1}(\ref{boundeq1}) it follows that $d\geq \frac{56^k}{63k+7}\geq\frac{56^2}{133} \geq 24.$

\subsubsection{$G_0\triangleright E_6(q)$}
\par \underline{\textbf{Case 1: $q=q_0$}}\par 
A Borel subgroup fixes a maximal $1$-space, so using Lemma  \ref{orbit-stabparaborel} we see that there is an orbit $\mathcal{O}$ of size $\vert \mathcal{O} \vert \leq\frac{q^{43}}{2},$ so we can use Lemma \ref{lemma2.1}(\ref{boundeq2}).
\par 
\underline{\textit{Case 1.a: $\lambda$ is $p$-restricted}}\par 
By Theorem \ref{lowerlprestrict}, either $n=27,$ $78-\epsilon_p(3)$ or $n\geq 324.$\par 
In the case when $n\geq 324,$ we use the fact a Borel fixes a maximal $1$-space and by Lemma \ref{lemma2.1}(\ref{boundeq2}), $d\geq 8.$ \par 
In the case when $n=78-\epsilon_p(3),$ $\lambda=\omega_2$ and $V$ is the adjoint module. By Lemma \ref{para}, it follows that the parabolic $P_2$ fixes a maximal $1$-space so using Lemma  \ref{orbit-stabparaborel}, $\vert \mathcal{O} \vert \leq q^{23}$ and so by Lemma \ref{lemma2.1}(\ref{boundeq1}), $d\geq 4.$ \par 
In the case when $n=27,$ $V$ is the minimal module. Here $G_0$ has 3 non-zero orbits on $V$ by \cite[Remark on page 468]{cohencooperstein} and hence $G$ has 3 non-diagonal orbitals, and so the orbital diameter is bounded above by 3. By \cite[Thm 1.1]{distancetransitive}, one of the orbital graphs is distance-transitive when $G_0$ contains the scalars in $GL_n(q_0)$, and so the orbital diameter is exactly 3.\par
\underline{\textit{Case 1.b: $\lambda$ is not $p$-restricted}} \par 
Hence $n\geq 27^2$ and so using the fact that a Borel fixes a maximal $1$-space, $d\geq 17.$ \par 
 \underline{\textbf{Case 2: $q=q_0^k$ for $k\geq 2$ as in Lemma \ref{q0lemma}(\ref{subfield})}}\par
Using the fact that a Borel fixes a maximal $1$-space we have that there is an orbit $\vert\mathcal{O}\vert \leq q_0^{36k+6}$ and so by Lemma \ref{lemma2.1}(\ref{boundeq1}) it follows that $d\geq \frac{27^k-1}{36k+6}\geq\frac{27^2-1}{78} \geq 10.$
\subsubsection{$G_0\triangleright\,^2\!E_6(q)$}
\par \underline{\textbf{Case 1}: $\tau_0(\lambda)\neq\lambda$}\par 
Now by Lemma \ref{q0lemma}, $^2E_6(q)\leq E_6(q^2)\leq GL(V)$ and so the lower bounds for the orbital diameter from the case of $G_0\triangleright E_6(q)$ hold. \par 
\underline{\textbf{Case 2}:$\tau_0(\lambda)= \lambda$ and $q=q_0$}\par
\underline{\textit{Case 2.a: $\lambda$ is $p$-restricted}}\par 
By Theorem \ref{lowerlprestrict}, either $n= 78-\epsilon_p(3)$ or $n\geq 572.$ \par 
The case when $n=78-\epsilon_p(3),$  $\lambda=\omega_2,$ the parabolic subgroup $P_2$ fixes a $1$-space. By Lemma  \ref{orbit-stabparaborel} there is an orbit of size at most $\vert \mathcal{O} \vert \leq \frac{q^{23}}{2},$ so by Lemma \ref{lemma2.1}(\ref{boundeq2}), $d\geq 4.$ \par 
For $n\geq 572$  we can use the fact that a Borel subgroup fixes a maximal $1$-space, so using Lemma  \ref{orbit-stabparaborel} we can see that there is an orbit of size $\vert\mathcal{O}\vert \leq {q^{42}},$ so we can use Lemma \ref{lemma2.1}(\ref{boundeq1}) to get that $d\geq 14.$\par
\underline{\textit{Case 2.b: $\lambda$ is not $p$-restricted}}
\par Hence $n\geq 77^2,$ and since a Borel fixes the $1$-space and by Lemma \ref{lemma2.1}(\ref{boundeq1}), $d\geq 142.$ 
\par 
 \underline{\textbf{Case 3: $\tau_0(\lambda)= \lambda$ and  $q=q_0^k$ as in Lemma \ref{q0lemma}(\ref{subfield})}}\par 
Using the fact that a Borel fixes a maximal $1$-space we have that there is an orbit $\vert\mathcal{O}\vert \leq q_0^{36k+6}$ and so by Lemma \ref{lemma2.1}(\ref{boundeq1}) it follows that $d\geq \frac{77^k-1}{36k+6}\geq\frac{77^2-1}{78} \geq 76.$


\subsubsection{$G_0\triangleright\, F_4(q)$}
 \par
\par \underline{\textbf{Case 1: $q=q_0$}}\par 

A Borel subgroup fixes a maximal $1$-space, so using Lemma  \ref{orbit-stabparaborel} we see that there is an orbit $\mathcal{O}$ of size $\vert \mathcal{O} \vert \leq\frac{q^{29}}{2},$ so we can use Lemma \ref{lemma2.1}(\ref{boundeq2}) to find lower bounds on the orbital diameter. 
\par 
\underline{\textit{Case 1.a: $\lambda$ is $p$-restricted}}\par
By Theorem \ref{lowerlprestrict}, either $n=26-\epsilon_p(3),$ $52$ or $n\geq 196.$ \par 
For $n=26-\epsilon_p(3),$ consider the algebraic group $\overline{G_0}=F_4(K)$ where $K=\overline{\mathbb{F}_p}$ acting on $\overline{V}=V_{26-\epsilon_p(3)}(K).$ Now the stabilizer of a maximal 1-space in $G_0$ is $P_4=Q_{15}B_3T_1.$ Let $\Delta_0$ be the orbit of $\overline{G_0}$ containing the maximal vectors. If $v=a+b,$ where $a,b\in \Delta_0,$ then $\overline{G_0}_a\cap \overline{G_0}_b\leq \overline{G_0}_v.$ The stabilizers of $a$ and $b$ are conjugates of $P_4'.$ Without loss of generality, they are $P_4'$ and $P_4'^g$ for some $g\in \overline{G_0}.$ We see from \cite[Lemma 5.38]{aluna2} what the intersections of two parabolics can be. 
Now consider the finite group $G_0=\overline{G_0}^{(q)}=F_4(q)$ acting on $V=\overline{V}^{(q)}=V_{26-\epsilon_p(3)}(q).$  By \cite[Table 2]{cohencooperstein}, $^3D_4(q).3$ stabilizes a vector $w\in V.$ By \cite[Lemma 5.38]{aluna2}, the possible intersections of $P_4'\cap P_4'^g$ with $\overline{G_0}^{(q)}$  either contain a unipotent subgroup of order at least $q^{13}$ or contain $B_3(q).$ Since neither of these is contained in $^3D_4(q).3,$ we cannot express $w$ as a sum of at most two elements in $\Delta_0\cap V,$ so the orbital diameter is at least 3. \par 

 The case when  $n=52,$ and $\lambda=\omega_1,$ by Lemma \ref{para}, the parabolic $P_1$ fixes a maximal $1$-space so using Lemma  \ref{orbit-stabparaborel}, $\vert \mathcal{O} \vert \leq q^{17}$ and so by Lemma \ref{lemma2.1}(\ref{boundeq1}), $d\geq 3.$\par  
 For $n\geq 196,$ we can use the fact that a Borel fixes a maximal $1$-space and so by Lemma \ref{lemma2.1}(\ref{boundeq2}) we deduce that $d\geq 7.$ 
\par
\underline{\textit{Case 1.b: $\lambda$ is not $p$-restricted}} \par 
 
Hence $n\geq 25^2$ and using the fact that a Borel fixes a maximal $1$-space,  $d\geq 22.$ \par 
 \par 
 \underline{\textbf{Case 2: $q=q_0^k$ for $k\geq 2$ as in Lemma \ref{q0lemma}(\ref{subfield})}}\par
 Using the fact that a Borel fixes a maximal $1$-space we have that there is an orbit $\vert\mathcal{O}\vert \leq q_0^{24k+4}$ and so by Lemma \ref{lemma2.1}(\ref{boundeq1}) it follows that $d\geq \frac{25^k-1}{24k+4}\geq\frac{25^2-1}{52} \geq 12.$

\subsubsection{$G_0\triangleright\,^2\!F_4(q)$}
By Lemma \ref{q0lemma}(\ref{twistedexceptional}), for each $V$ we have that $$^2F_4(q)\leq F_4(q)\leq GL(V)$$ and so all bounds from the case $G_0 \triangleright F_4(q)$ hold and the result follows.

\subsubsection{$G_0\triangleright\, G_2(q)$}
\par \underline{\textbf{Case 1: $q=q_0$}}\par 
 A Borel subgroup fixes a maximal $1$-space, so using Lemma  \ref{orbit-stabparaborel} we see that there is an orbit $\mathcal{O}$ of size $\vert \mathcal{O} \vert \leq{q^{8}},$ so we can use Lemma \ref{lemma2.1}(\ref{boundeq1}) to find lower bounds on the orbital diameter. 
 \par 
\underline{\textit{Case 1.a: $\lambda$ is $p$-restricted}}\par
By Theorem \ref{lowerlprestrict}, either $n=7-\epsilon_p(2),$ $14$ or $n\geq 26.$ \par 
For $\lambda=\omega_2$ and $p=2$, $G_0$ acts transitively on $V$, so the orbital diameter is 1 when $G_0$ contains the scalars of $GL_n(q_0)$. For $p$ odd, $G$ has orbital diameter 2 by Lemma \ref{smalldiamlemma}. \par
For $n=14,$  the parabolic subgroup $P_2$ fixes a maximal $1$-space, so using Lemma  \ref{orbit-stabparaborel} we can see that there is an orbit of size at most $q^6-1.$ Now by Lemma \ref{lemma2.1}(\ref{boundformula1}), $d\geq 3.$ \par 
For $n\geq 26$ using the fact that a Borel fixes a maximal $1$-space, $d\geq 3.$ 
\par
\underline{\textit{Case 1.b: $\lambda$ is not $p$-restricted}} \par 
 
Hence $n\geq 6^2$ and using the fact that a Borel fixes a maximal $1$-space,  $d\geq 5.$ \par 
 \par 
 \underline{\textbf{Case 2: $q=q_0^k$ for $k\geq 2$ as in Lemma \ref{q0lemma}(\ref{subfield})}}\par
Using the fact that a Borel fixes a maximal $1$-space we have that there is an orbit $\vert\mathcal{O}\vert \leq {q_0^{6k+2}}$ and so by Lemma \ref{lemma2.1} part \ref{boundformula1} it follows that $d\geq \frac{6^k-1}{6k+2}\geq\frac{6^2-1}{14} \geq 3.$
\subsubsection{$G_0\triangleright\,^2\!G_2(q)$}
By Lemma \ref{q0lemma}(\ref{twistedexceptional}), for each $V,$ we have that $$^2\!G_2(q)\leq G_2(q)\leq GL(V)$$ and so all bounds from the case $G_0 \triangleright G_2(q)$ hold. As $p=3$ in this case, the representation for $n=6$ does not exist. \par 
For the case when $n=7$, we show that $d\geq 3.$ The orbits of $^2G_2(q)$ are described in \cite{reewilson}. Let $q=3^{2m+1}.$ Sticking to the notation in \cite{reewilson}, the module $V$ has basis $\{e_{-3},e_{-2},e_{-1},e_0,e_1,e_2,e_3,\}$ where $e_i$ is the row vector that has all zeros except in position $i+4,$ where it has a 1. We want to show that we cannot express every element in the vector space as a sum of two elements in the orbit of $e_{-3},$ call this orbit $\mathcal{O}.$ From the proof of \cite[Lemma 3]{reewilson} we know that the orbit $\mathcal{O}$ consists of the $q-1$ scalar multiples of $e_{-3}$ and the $(q-1)q^3$ of images of $e_{3}$ under the action of the stabilizer of $\langle e_{-3} \rangle.$ The generators  of the stabilizer of $\langle e_{-3} \rangle$ are given in \cite[Lemma 1]{reewilson} and the proof of \cite[Lemma 2]{reewilson} and they are acting by right multiplication on the vectors. These generators are diagonal matrices of the form $$D_i=diag(a_i^{-1},b_i^{-1},\lambda_i^{-1},1,\lambda_i,b_i,a_i)$$ where $\lambda_i\in \mathbb{F}_q^\star,$ $a_i=\lambda_i^{3^{m+1}+2}$ and $b_i=\lambda_i^{3^{m+1}+1},$\tiny $$A=\begin{pmatrix}
1 &  &&&&&& \\
1 & 1 & &&&&&\\
-1&-1&1&&&&\\
0&0&1&1&&&\\
1&0&1&-1&1&&\\
-1&0&1&-1&1&1&\\
-1&-1&-1&0&0&-1&1
\end{pmatrix},\,\,B=\begin{pmatrix}
1 &  &&&&&& \\
0 & 1 & &&&&&\\
-1&0&1&&&&\\
0&1&0&1&&&\\
1&0&0&0&1&&\\
0&1&0&-1&0&1&\\
1&0&-1&0&1&0&1
\end{pmatrix}, \,\text{and}\,\, C=\begin{pmatrix}
1 &  &&&&&& \\
0 & 1 & &&&&&\\
0&0&1&&&&\\
1&0&0&1&&&\\
0&-1&0&0&1&&\\
1&0&1&0&0&1&\\
1&-1&0&-1&0&0&1
\end{pmatrix}.$$ \normalsize The images of $\langle e_3 \rangle$ are represented by the last row of the matrices. We can see that the orbit of $e_3$ under the stabilizer contains vectors of the form $$e_3A^{D_1}C^{D_3}B^{D_2}=\lambda(x_{-3}e_{-3}+x_{-2}e_{-2}+x_{-1}e_{-1}+x_{0}e_{0}+x_1e_1+x_2e_2+e_3),$$ where $$x_{-3}=a_1^{-2}+a_1^{-1}\lambda_1^{-1}\lambda_2a_2^{-1}-b_1a_1^{-1}(b_3^{-1}a_3^{-1}-b_3^{-1}\lambda_3^{-1}\lambda_2a_2^{-1})-a_3^{-2}+a_2^{-2},$$ $$x_{-2}=-a_1b_1^{-1}+b_1a_1^{-1}b_2^{-2}-a_3^{-1}b_3^{-1}-a_3^{-1}b_2^{-1},$$ $$x_{-1}=-a_1^{-1}\lambda_1^{-1}-a_2^{-1}\lambda_2^{-1},$$ $$x_0=-b_2^{-1}b_1a_1^{-1}-a_3^{-1},$$ $$x_1=a_2^{-1}\lambda_2$$ and $$x_2=-b_1a_1^{-1}.$$  \par 
This is the whole orbit, as the number of distinct vectors of this form is exactly the size of $\mathcal{O}.$
Now we will show that there is a vector in $\langle e_{-2}\rangle$ that we cannot express as a sum of two elements in the orbit $\mathcal{O}.$ This will imply that the orbital diameter is at least 3, since $\langle e_{-2}\rangle \cap \mathcal{ O}= \emptyset.$ Suppose $v_1,v_2\in \mathcal{O}$ such that $k_1 v_1+k_2 v_2=e_{-2},$ so $v_1+k_1^{-1}k_2 v_2=k_1^{-1}e_{-2}.$ Let $k_1^{-1}k_2=\lambda.$ Without loss of generality, the $v_i$s are of the form $x_{-3,i}e_{-3}+x_{-2,i}e_{-2}+x_{-1,i}e_{-1}+x_{0,i}e_{0}+x_{1,i}e_1+x_{2,i}e_2+e_3,$ as if any of them were in $\langle e_{-3}\rangle$ then we would arrive to a contradiction immediately. Consider $v= v_1+\lambda v_2.$ Now the coefficients of $e_3,$ $e_2,$ $e_1,$ $e_0,$ $e_{-1}$ and $e_{-3}$ in $v$ are all $0.$ The coefficient of $e_{3}$ in $v$ is $1+\lambda=0,$ so we conclude that $\lambda=-1.$ \par 
Using this we see that the coefficient of $e_{2}$ in $v$ is $x_{2,1}-x_{2,2}=0,$ so $b_{1,1}a_{1,1}^{-1}=b_{1,2}a_{1,2}^{-1},$ so $\lambda_{1,1}^{-1}=\lambda_{1,2}^{-1}$ and so $\lambda_{1,1}=\lambda_{1,2},$ $a_{1,1}=a_{1,2}$ and $b_{1,1}=b_{1,2}.$ \par 
The coefficient of $e_1$ in $v$ is $x_{1,1}-x_{1,2}=0,$ so $a_{2,1}^{-1}\lambda_{2,1}-a_{2,2}^{-1}\lambda_{2,2}=0,$ so $\lambda_{2,1}^{-3^{n+1}-1}=\lambda_{2,2}^{-3^{n+1}-1},$ and since $\gcd(3^{n+1}+1,3^{2n+1}-1)=2,$ $\lambda_{2,1}=\pm \lambda_{2,2},$ $a_{2,1}=\pm a_{2,2}$ and $b_{2,1}=b_{2,2},$ as $b_{2,i}$ is an even power of $\lambda_{2,i}.$ \par 
The coefficient of $e_0$ in $v$ is $x_{0,1}-x_{0,2}=0,$ so $-b_{2,1}^{-1}b_{1,1}a_{1,1}^{-1}-a_{3,1}^{-1}+b_{2,2}^{-1}b_{1,2}a_{1,2}^{-1}+a_{3,2}^{-1}=0.$ Now using the facts that $a_{1,1}=a_{1,2},$ $b_{1,1}=b_{1,2}$ and $b_{2,1}=b_{2,2},$ this shows that $a_{3,1}^{-1}=a_{3,2}^{-1},$ so $a_{3,1}=a_{3,2},$ and since $\gcd(3^{n+1}+2,3^{2n+1}-1)=1,$ also $\lambda_{3,1}=\pm \lambda_{3,2}$ and $b_{3,1}=b_{3,2}.$ \par 
Now the coefficient of $e_{-2}$ in $v$ is $x_{-2,1}-x_{-2,2}=-a_{1,1}b_{1,1}^{-1}+b_{1,1}a_{1,1}^{-1}b_{2,1}^{-2}-a_{3,1}^{-1}b_{3,1}^{-1}-a_{3,1}^{-1}b_{2,1}^{-1}+
a_{1,2}b_{1,2}^{-1}+b_{1,2}a_{1,2}^{-1}b_{2,2}^{-2}-a_{3,2}^{-1}b_{3,2}^{-1}-a_{3,2}^{-1}b_{2,2}^{-1}=0$ as  $b_{1,1}=b_{1,2},$ $b_{2,1}=b_{2,2},$ $b_{3,1}=b_{3,2},$ $a_{1,1}=a_{1,2}$ and $a_{3,1}=a_{3,2}.$ This shows that we cannot express $k_1^{-1}e_{-2}$ as a sum of two elements in $\mathcal{O}$ and so the orbital diameter is at least 3.
\subsubsection{$G_0\triangleright\,^2\!B_2(q)$}
By Lemma \ref{q0lemma}(\ref{twistedexceptional}), for each $V$ we have that $$^2\!B_2(q)\leq B_2(q)\leq GL(V)$$ and so all bounds from the case $G_0 \triangleright B_2(q)$ hold. This shows that  either $d\geq 3$ or $V=V_4(q).$ In the latter case $d=2,$ because $G$ has rank 3 by \cite{liebeckaffine} when $G_0$ contains the scalar matrices of $GL_n(q_0)$ and so the orbital diameter is 2 by Lemma \ref{smalldiamlemma}.
\subsubsection{$G_0\triangleright\, ^3D_4(q)$}\par
Recall that here $\tau_0$ is the graph automorphism of $D_4$ of order $3.$
\par \underline{\textbf{Case 1}: $\tau_0(\lambda)\neq \lambda$}\par
Now $^3D_4(q)\leq D_4(q^3)\leq GL(V)$ and so the lower bounds for the orbital diameter from the case of $G_0\triangleright D_4(q)$ hold. \par 
\underline{\textbf{Case 2}: $\tau_0(\lambda)= \lambda$ and $q=q_0$}\par
\underline{\textit{Case 2.a: $\lambda$ is $p$-restricted}}\par 
By Theorem \ref{lowerlprestrict}, either $n=28-2\epsilon_p(2)$ or $n\geq 195.$ \par 
In the case when $n=28-2\epsilon_p(2)$ and $\lambda=\omega_2,$ the parabolic subgroup $P_2$ fixes a $1$-space. By Lemma  \ref{orbit-stabparaborel} there is an orbit of size at most $\vert \mathcal{O} \vert \leq \frac{q^{11}}{2},$ so by Lemma \ref{lemma2.1}(\ref{boundeq2}), $d\geq 3.$ \par 
For $n\geq 195$ we can use the fact that a Borel subgroup fixes a maximal $1$-space, so using Lemma  \ref{orbit-stabparaborel} we can see that there is an orbit of size $\vert\mathcal{O}\vert \leq\frac{q^{17}}{2},$ so we can use Lemma \ref{lemma2.1}(\ref{boundeq2}) to get that $d\geq 12.$\par
\underline{\textit{Case 2.a: $\lambda$ is not $p$-restricted}}\par
Here $n\geq 26^2,$ and since a Borel fixes a maximal $1$-space, by Lemma \ref{lemma2.1}(\ref{boundeq2}), $d\geq 40.$\par 
\par 
 \underline{\textbf{Case 3: $\tau_0(\lambda)= \lambda$ and  $q=q_0^k$ as in Lemma \ref{q0lemma}(\ref{subfield})}}\par
 Using the fact that a Borel fixes a maximal $1$-space we have that there is an orbit $\vert\mathcal{O}\vert \leq \frac{q_0^{12k+3}}{2}$ and so by Lemma \ref{lemma2.1}(\ref{boundeq2}) it follows that $d\geq \frac{26^k}{12k+3}\geq\frac{26^2}{27} \geq 25.$

\section{Alternating Stabilizer}

In this section we prove our results for the case when $G_s\cong A_n$.
The groups considered in this section all satisfy the following hypothesis.
\begin{hypothesis}\label{hypothesisalternating}
 Let $G=VG_0$ be a primitive affine group such that $G_s:=\frac{G_0^\infty}{Z(G_0^\infty)}=A_r,$ where $A_r$ is an  alternating group. Suppose that $V$ is an absolutely irreducible $\mathbb{F}_{q_0}G_0^\infty-$module in characteristic $p.$ Assume that $V$ cannot be realised over a proper subfield of $\mathbb{F}_{q_0}.$ Let $n=\dim V$ and $d=orbdiam(G,V).$ 
\end{hypothesis}
\subsection{The Fully Deleted Permutation  Module}
 

 Here we prove Proposition \ref{permalt}. Recall that the fully permutation  module for $G=A_r$ or $S_r$  is ${W}/{W\cap T},$ where $W:=\{(a_1,\dots,a_r):\sum a_i=0\}\leq \mathbb{F}_{q_0}^r$ and $T:=Span(1,\dots,1).$

\begin{proof} [Proof of Proposition \ref{permalt}]
Note that it is sufficient to prove the result for $G_0=\mathbb{F}_{q_0}^\star S_r,$ and then Proposition \ref{permalt} will follow from Lemma \ref{subgroup}. \par 

For $a\in \mathbb{F}_{q_0}^\star,$ write $\underline{a}=(a,\dots,a).$
Consider the orbital $\Delta:=\{\underline{0},(1,-1,0,\dots,0)\}^{G}.$ Denote the distance in the corresponding orbital graph between two elements $v,w\in V_n(q_0)$ by $d(v,w).$  We have two cases to consider.\par
\textbf{\underline{Case 1}} Assume $p\nmid r$ and so $T$ is not contained in $W$ and $n=r-1.$ \par
Denote the number of zeros of $v\in V_n(q_0)$ by $z(v).$ \par 
\underline{Claim 1.1} Let $m\leq \frac{r}{2}$. For all $v\in V_n(q_0)$ such that $d(\underline{0},v)\leq m$ we have $z(v)\geq r-2m.$ \par This is clear, as the neighbours of a vector $h$ are of the form $h\pm (1,-1,0,\dots,0)^{G_0}$ so they have a maximum of 2 extra non-zero entries. \par 
\underline{Claim 1.2} There exists an element $w\in V_n(q_0)$ such that $d(\underline{0},w)\geq \frac{r-1}{2}.$\par 
This element is $$w=\begin{cases} (1,-1,\dots,1,-1),\,\, r\,\,\mbox{ even} \\ (1,-1,\dots,1,-1,0),\,\, r\,\,\mbox{ odd}
\end{cases}$$ 
Now $z(w)\leq r-2(\frac{r-1}{2})$ so by the contrapositive of Claim 1.1, $d(\underline{0},w)\geq \frac{r-1}{2}.$
Hence we get $\frac{r-1}{2}$ as a lower bound on the orbital diameter, as required for Proposition \ref{permalt}. \par 
\textbf{\underline{Case 2}} Assume $p \vert r$ and so $T$ is contained in $W$ and $V=\frac{W}{T}.$ \par 
\underline{Claim 2.1} Let $m\leq \frac{r}{2}$. Every coset of distance at most $m$ away from $T$ has a coset representative $v$ such that $z(v)\geq r-2m.$ \par 
The proof is analogous to the proof of Claim 1.1.\par 
\underline{Claim 2.2} There exists a coset $w+T$ such that $d(T,w+T)\geq \frac{r-2}{4}.$ \par 
To prove this, we need to find a coset $w+T$ such that for all $u\in w+T$, $z(u)\leq r-\frac{r-1}{2}.$ Note that $u$ is of the form $w+\underline{a}$ where $a\in \mathbb{F}_{q_0}^\star.$ The elements $$w=\begin{cases} (1,-1,\dots,1,-1)+T,\,\,r \,\, \mbox{even,}\,\,p\,\, \mbox{odd} \\ (1,-1,\dots,1,-1,0)+T, \,\,r \,\, \mbox{odd,}\,\,p\,\, \mbox{odd}\\
(1,0,\dots,1,0)+T,\,\, 4\vert r,\,\,p\,\,\mbox{ even} \\
(1,0,\dots,1,0,0,0)+T,\,\, 4\nmid r,\,\,r \,\, \mbox{even,}\,\,p\,\, \mbox{even}
\end{cases}$$ satisfy this, and so they are at distance at least $\frac{r}{2},$ $\frac{r-1}{2},$ $\frac{r}{2},$ or $\frac{r-2}{2},$ away from $T$, respectively. This is lowest for $\frac{r-2}{4},$ so the lower bound for the orbital diameter holds. 

\end{proof}
We use this to prove Corollary \ref{classification}, which is a classification of such groups with orbital diameter at most 2.
\begin{proof}[Proof of Corollary \ref{classification}]

For $r\geq 5$ and $r\neq 6$ the automorphism group of $A_r$ is $S_r,$ so  the orbital diameter is minimal for $G_0=\mathbb{F}_{q_0}^\star S_r.$ Note that for $r=6,$ the fully deleted permutation module only exists for $G_0^\infty = S_6$ or $A_6$, so again the orbital diameter is minimal for $G_0=\mathbb{F}_{q_0}^\star S_r.$\par  

\textbf{\underline{Case 1}} $p\vert r.$\par 
Now Proposition \ref{permalt} gives us that $r-2\leq 4d,$ and so if $d\leq2$ then $r\leq 10.$ Hence we have the following possibilities: \[\begin{array}{c|c|c|c|c|c|c|c|c}
   r  &  10&10&9&8&7&6&6&5\\ \hline
   p  & 2&5&3&2&7&2&3&5\\
\end{array}\]
To determine whether these indeed have orbital diameter 2 or 1 we use computational methods as described in Section 2.  We find that $orbdiam(G,V)=1$ if and only if $r=6$ and $q_0=2$ and $orbdiam(G,V)=2$ if and only if $q_0=2$ and $r=8$ or $10$ or $q_0=3$ and $r=6$ or $p=5,$ $r=5$ and $4\times A_5\leq  G_0.$ \par 

 \par 
\textbf{\underline{Case 2}}  $p\nmid r.$\par Now Proposition \ref{permalt} tells us that there is no possibility for the orbital diameter $1$ case. For the orbital diameter $2$ case the only possibility is $r=5$ and $p\neq 5.$ 
Using computation we can show that for the case of $(r,q_0)=(5,2)$ the orbital diameter is $2$ for any $G_0$ with $G_s\cong A_5,$ and that for $(r,q_0)= (5,3),$ and $(5,4)$ the orbital diameter is at least 3.\par 
For the case $r=5$ and $p\geq 7,$ we use another method. It suffices to show that $orbdiam(V\mathbb{F}_p^\star S_r,V)\geq 3.$ As in the proof of Proposition \ref{permalt}, consider the orbit $\Delta=(1,-1,0,0,0)^{G_0}=\{k(1,-1,0,0,0)^\sigma \vert k\in \mathbb{F}_{q_0}^\star,\,\, \sigma\in S_r\}.$ 
Then $(1,1,1,-3,0)\in W$ cannot be expressed as a sum of one or two vectors in $\Delta.$ Hence in these cases, the orbital diameter is at least 3, and the result follows. \par 
\end{proof}

\subsection{The Bounds on the Orbital Diameter} 

We now have all the information we need to prove the remaining bounds on the orbital diameter for the case when $G_s\cong A_r$.

\begin{proof}[Proof of Theorem \ref{fdtheoremalt}]
Since $V_n(q_0)$ is not the fully deleted permutation  module,  \cite[Thm 2.2]{tiep} shows that for $r\geq 15,$ $n\geq \frac{r(r-5)}{2}\geq \frac{1}{3}r^2$ which is equivalent to $r\leq \sqrt{3n}.$ From Lemma \ref{eq3} we have $n\leq 1+dr\log_2(r)\leq 2dr\log_2(r).$ Putting these together we have that $n\leq \sqrt{3}d\sqrt{n}\log_2(3n) \leq2\sqrt{3}d\sqrt{n}\log_2(n)$ so $\sqrt{n}\leq 2\sqrt{3}d\log_2(n).$ Now let $\delta=\frac{\epsilon}{2(2+\epsilon)}.$ For sufficiently large $n,$ we have $ 2\sqrt{3}\log_2(n)\leq n^\delta$ and so $\sqrt{n}\leq dn^\delta$ which gives the desired result $n\leq d^{2+\epsilon}.$ 
\end{proof}
 
\begin{proof}[Proof of Theorem \ref{lowerdiamalt}]
Assume $V$ is not the fully deleted permutation module.
Lemma \ref{eq3} gives that $ n\leq 1+d\log_2(\vert Aut(G_s)\vert)$.We get the bound \begin{equation}\label{eqbound1}
    d\geq \frac{n-1}{\log_2(\vert Aut(G_s)\vert)}.
\end{equation}


  By \cite[Thm 2.2]{tiep} for $r\geq 15,$ we have  $n\geq r(r-5)/2.$ Hence (\ref{eqbound1}) gives $$d\geq \frac{r^2-5r-2}{2r\log_2(r)}\geq \frac{r-6}{2\log_2(r)}.$$
  
  For $r\leq 14,$ the orbital diameter 1 cases are given by Theorem \ref{diam1grps}, so $d\geq 2$ for $r\geq 9.$

\end{proof}

\subsection{Alternating Affine Groups with Orbital Diameter 2}
In this section we provide a classification of groups of the form $V_n(q_0)G_0$ where $\frac{G_0^\infty}{Z(G_0^\infty)}\cong A_r$ an alternating group, and $orbdiam(V_n(q_0)G_0)\leq 2.$ 

\begin{proof}[Proof of Theorem \ref{if2}]
Assume $orbdiam(G,V)\leq 2$ and $V$ is not the fully deleted permutation module. By \cite[Theorem 2.2]{tiep} we have  $n\geq \frac{r(r-5)}{2}$ for $r\geq 15$. First note that the bound from Theorem \ref{lowerdiamalt} gives that $r\leq 23.$ Using Lemma \ref{lemma2.1}(\ref{boundformula1})  we can improve this, as if $orbdiam(G,V)\leq 2,$ then $1+(q_0-1)r!+((q_0-1)r!)^2\geq q_0^{\frac{r(r-5)}{2} }.$ By Lemma \ref{log2.2}, $G$ can only have orbital diameter $2$ if $$log_2(1+r!+(r!)^2)\geq \frac{r(r-5)}{2},$$ hence we can conclude that $r\leq 16.$ Now to prove the theorem we will consider the alternating groups $A_r$ with $5\leq r \leq 16$ in turn. We use Lemma \ref{lemma2.1} to bound the dimension $n$ and \cite{hissmalle} which lists all irreducible representations of dimension up to 250. \par 
\begin{itemize} 
    \item ${A_{16}}, A_{15}, A_{14}.$ \par  Lemma \ref{lemma2.1}(\ref{boundformula1})  and \cite{hissmalle} gives that the only representation possible has dimension 64 over the field of $2$ elements. Here $A_{14}\leq A_{15}\leq A_{16}\leq \Omega_{14}^+(2)\leq GL_{64}(2),$ see \cite[p.187, 195]{kleidmanliebeck}, and so this is the restriction of the spin representation of $\Omega_{14}^+(2)$ to $A_{r}.$ Since for $ \Omega_{14}^+(2)$ the orbital diameter is at least $3$ by Theorem \ref{diam2}, the same holds for $G_s=A_r$ here.
    \item ${A_{13}}.$\par  By Lemma  \ref{lemma2.1} and \cite{hissmalle} there are three possible cases: $(n,q_0)=(64,2),$ $(32,4)$ and $(32,3).$   \par The 64-dimensional case is again excluded by Lemma \ref{subgroup} since  $A_{13}\leq A_{16}\leq GL_{64}(2)$. \par 
    In the case of $(n,q_0)=(32,4),$ we see using Magma \cite{magma} that an element of order 13 stabilizes a vector so there is an orbit of size at most $\frac{\vert G_0 \vert }{13}.$ Lemma \ref{lemma2.1}(\ref{boundformula1}) shows that the orbital diameter is greater than or equal to 3. \par 
    Now consider the case $(n,q_0)=(32,3),$ which exists only for $G_0^{\infty}\cong 2.A_{13}.$ This the restriction of the spin representation of $D_6(3)$ to $A_{13}.$ Since the case when $G_s\cong P\Omega_{12}^+(3)$ has orbital diameter greater than 3 by Theorem \ref{diam2}, the same holds for $G_s\cong A_{13}$ by Lemma \ref{subgroup}.  \par 
\item $A_{12}.$ \par By Lemma  \ref{lemma2.1} and \cite{hissmalle} there are three possible cases; $(n,q_0)=(44,2),$ $(16,4)$ and $(16,3).$\par  Suppose $(n,q_0)=(44,2).$ This representation is an irreducible composition factor of the exterior square of the fully deleted permutation module. A subgroup $A_8$ fixes a 3-dimensional subspace pointwise in the fully deleted permutation module. This implies, that $A_8$ fixes at least a 1-dimensional space pointwise in the 44-dimension module as well, and so $G_0$ has an orbit on $V$ of size $\leq \frac{\vert S_{12}\vert}{\vert A_8\vert}.$ We can apply Lemma \ref{lemma2.1}(\ref{boundformula1}) to show that $orbdiam(G,V)\geq 3$. \par Suppose $(n,q_0)=(16,3).$ This only exists for for $G_0^{\infty}\cong 2.A_{12}.$ We can construct the representation of $G_0$ in Magma and compute all orbits of $G_0$ on $V$. By our computations, there is a vector $v\in V$ that cannot be expressed as a sum of at most two elements from the orbit of size  at most 60480, hence $orbdiam(G,V)\geq 3.$ \par 
The case when $n=(16,4)$ is in part (2) of Theorem \ref{if2}. 
\item $A_{11}.$ \par  By Lemma  \ref{lemma2.1} and \cite{hissmalle} there are five possible cases; $(n,q_0)=(44,2),$ $(16,4)$ and for $G_0^{\infty}\cong 2.A_{11},$ $(n,q_0)=(16,3), $ $(16,5),$ and $(16,11).$ \par 
For $(n,q_0)=(44,2),$ $A_{11}\leq A_{12}\leq GL_{44}(2)$ so by Lemma \ref{subgroup} the orbital diameter here is also at least 3. \par 
For $(n,q_0)=(16,3),$  $2.A_{11}\leq 2.A_{12}\leq GL_{16}(3),$ so by Lemma \ref{subgroup} the orbital diameter at least 3. \par 
The cases where $(n,q_0)=(16,4)$ or $(16,5)$ are in part (2) of Theorem \ref{if2}. \par 
For $(n,q_0)=(16,11),$ we see from the Brauer character table \cite{bratlas}, that an element $h$ of order 7 takes value $2,$ hence fixes a vector,
and so $G_0$ has an orbit of size $\leq \frac{10\vert S_{11}\vert }{7}$, so $orbdiam(G,V)\geq 3$ by Lemma \ref{lemma2.1}(\ref{boundformula1}).  \par 
\item $A_{10}.$  \par By Lemma  \ref{lemma2.1} and \cite{hissmalle} these are the possibilities: $(n,q_0)=(26,2)$ and $(16,2)$ for $G_0^{\infty}\cong A_{10}$ and $(16,3),$ $(16,7)$ and $(8,5)$ for $G_0^{\infty}\cong 2.A_{10}.$ \par 
In the case where $(n,q_0)=(26,2),$ $V$ is an irreducible composition factor of the exterior square of the fully  deleted permutation  module. We can show that $A_6,$ (and if $S_{10}\leq G_0$ then $S_6,$) fixes a 3-space pointwise in the fully  deleted permutation  module. This means that $A_6,$ respectively $S_6,$ fixes at least one vector in the 26-dimensional module as well. Now using Lemma \ref{lemma2.1}(\ref{boundformula1})  we can exclude this case.    \par 
For $(n,q_0)=(16,2),$ the representation can be constructed in Magma and we can show that  the orbital diameter is at least 3, because we cannot express every vector in $V$ as a sum of two elements in the orbit of size $945.$ \par 
For $(n,q_0)=(16,3),$ we can see from the character tables in \cite{bratlas} that this is the restriction of the irreducible $16$-dimensional module of $A_{12}$ to $A_{10.}$ Now $2.A_{10}\leq 2.A_{12}\leq GL_{16}(3),$ so by Lemma \ref{subgroup} the orbital diameter is at least 3. \par 
For $(n,q_0)=(16,7)$ then we can see from the character table in \cite{bratlas} that there is an element of order 8 that fixes a vector. Hence this is excluded using Lemma \ref{lemma2.1}(\ref{boundformula1}). 
\par For $(n,q_0)=(8,5)$ with $G_0^{\infty}\cong 2.A_{10}$ we can compute the orbits with GAP. We check that we can cannot express all elements of $V$ as a sum of at most two vectors in an orbit of size $2400$ and so the orbital diameter is at least 3. 

\item $A_{9}.$ \par  By Lemma  \ref{lemma2.1} and \cite{hissmalle} the possibilities are: $(n,q_0)=(26,2),$ $(21,3), (20,2),$ $(8,2)$ and $(8, \text{odd}) $ for $G_0^{\infty}\cong 2.A_9.$ \par 
For $(n,q_0)=(26,2)$ we have that $A_9\leq A_{10}\leq GL_{26}(2)$ so the orbital diameter is at least 3.\par 
For $(n,q_0)=(21,3),$ $V$ is the irreducible wedge square of the fully  deleted permutation  module. A subgroup $A_6$ fixes a 2-space in the fully  deleted permutation  module pointwise, hence a vector in the wedge square, so we can exclude this case also using Lemma \ref{lemma2.1}(\ref{boundformula1}). \par Consider $(n,q_0)=(20,2).$ We see using GAP that this has an orbit of size 360 and so it is excluded by Lemma \ref{lemma2.1}(\ref{boundformula1})  as well.\par 
For $(n,q_0)=(8,2),$ $G$ is a rank 3 group, so $orbdiam(G,V)=2$ as in part (2) of the theorem.\par 
Consider the $8$-dimensional representation for $2.A_{9}$. We can see from the Brauer character table that this representation is the restriction of the $8$-dimensional representation of $2.\Omega_8^+(2).$ By \cite{ATLAS} we know that $Sp_6(2)$ is a maximal subgroup of $\Omega_8^+(2),$ and that $2.\Omega_8^+(2)$ has a subgroup of the form  $2\times Sp_6(2),$ and so $Sp_6(2)\leq 2.\Omega_8^+(2).$ We can see from \cite{ATLAS} and \cite{bratlas} that the restriction of the 8-dimensional representation in question to $Sp_6(2)$ has a $7$-dimensional composition factor and since the representation is self-dual, $Sp_6(2)$ fixes a vector. Hence $G_0$ has an orbit of size at most $(q_0-1)\frac{\vert 2.\Omega_8^+(2)\vert }{\vert Sp_6(2)\vert } \leq 240(q_0-1)$ and so this case is excluded by Lemma \ref{lemma2.1}(\ref{boundformula1}) for $q_0\geq 7.$ For $q_0=5$ and $3$ we construct the representation in GAP, which tells us that the diameter is not 2. \par 
\item $A_8.$\par Since $A_8\cong SL_4(2),$ the case when $p=2$ has already been considered in Theorem \ref{aldiam2}. 
\par  By Lemma  \ref{lemma2.1} and \cite{hissmalle} the remaining possibilities are: $(n,q_0)=(13,3), (13,5)$ for $G_0^{\infty}\cong A_8,$ and $(8, \text{odd}) $ for $G_0^{\infty}\cong 2.A_8.$ \par 
For $(n,q_0)=(13,3),$ the Brauer character in \cite{bratlas} shows that the orbital diameter is at least 3. \par 
For $(n,q_0)=(13,5)$ we can show using MAGMA that an element of order 15 in $A_8$ fixes a 1-space and so we can exclude this case too by Lemma \ref{lemma2.1}(\ref{boundformula1}.) \par  The irreducible $8-$dimensional representation of $2.A_8$ is the restriction of the $8-$dimensional representation of $2.\Omega_8^+(2),$ just like in the case of $2.A_9,$ so this case inherits the lower bound of 3 by Lemma \ref{subgroup}. \par 
\item $A_7.$ \par
The case where $n\leq 9$ is in conclusion (1) of the theorem, so assume $n\geq 10.$\par 

By Lemma  \ref{lemma2.1} and \cite{hissmalle} the possibilities with $n\geq 10$ are $(n,q_0)=(20,2),$ $(15,3),$ $(14,2),$ $(13,3),$  and $(10,7).$ 
For $(n,q_0)=(20,2)$  or $(10,7)$ we can see from the Brauer character tables in \cite{bratlas} that an order 7 or 5 element fixes a 1-space, respectively, so we exclude these using Lemma \ref{lemma2.1}(\ref{boundformula1}). \par 
The representation with $(n,q_0)=(15,3)$ is the exterior square of the fully  deleted permutation  module. A subgroup $A_4$ fixes a 2-space in the fully  deleted permutation  module pointwise, and so a vector in the wedge square as well. Hence we exclude this also by Lemma \ref{lemma2.1}(\ref{boundformula1}). \par 
The representation with $(n,q_0)=(14,2)$ is an irreducible composition factor of the wedge square of the fully  deleted permutation  module. Similarly as before, a subgroup $A_5$ fixes a 2-space in the fully  deleted permutation  module pointwise, so it fixes a non-zero vector in the 14-dimensional composition factor in question. Now we can exclude this case too. \par 
For $(n,q_0)=(13,3)$ we construct the representation in GAP and show that $G_0$ has an orbit of size at most $70,$ so we can exclude this case using Lemma \ref{lemma2.1}(\ref{boundformula1}) .
\item $A_6.$ \par  By Lemma \ref{lemma2.1}(\ref{boundformula1})  and \cite{hissmalle} the only possibility for $n\geq 10$ is $(n,p)=(10,5).$ We can see from the Brauer ATLAS \cite{bratlas} that an element of order $3$ fixes a 1-space so we can exclude this using Lemma \ref{lemma2.1}(\ref{boundformula1}).  \par \item $A_5.$ \par There are no possibilities with $n\geq 10.$
\end{itemize}
\end{proof}

The following example list some groups with small orbital diameter.
\begin{example}\label{partialconverse}
     Let $G$ be as in Hypothesis \ref{hypothesisalternating}. If one of the following holds and $G_0$ contains the scalars $\mathbb{F}_{q_0}^\star$ in $GL_n(q_0)$, then $orbdiam(G,V)\leq 2.$
     \begin{enumerate}
             \item $(r,n,q_0)$ are as in Corollary \ref{classification}.
             \item $n=2.$
             \item $(r,n,q_0)$ are as in Theorem \ref{diam1grps}, so $G$ is $2$-homogeneous.
             \item $n,$ $q_0$ and $G_0$ are as in the following cases. 
             \[\begin{array}{c|c|c|c|c|c|c|c}
           G_0\triangleright  &  A_9 & A_8 & 2.A_7 & 3.A_7 & 3.A_6 & A_6 & A_5 \\\hline
     n&           8 & 4 & 4 & 3 & 3 & 3 & 3 \\\hline
  q_0&  2 & 2 & 7 & 25 & 4 & 9 & 9 \\\hline

\end{array}\]
     \end{enumerate}
\end{example}

\begin{proof}
Parts 1-3 are clear.\par Now consider part (4).  In the case where $(n,q_0)=(8,2)$ and $A_9\triangleright G_0$ produces a rank three group by \cite{liebeckaffine}, so $d=2.$\par The case where $(n,q_0)=(4,2)$ and $A_8\triangleright G_0,$ $V$ is the natural module of $SL_4(2)\cong A_8$ so the orbital diameter is $1.$
\par The remaining cases were proved by computations in GAP and MAGMA. 

\end{proof}

\section{Lie type Stabilizer In Cross Characteristic}

In this section we prove Theorems \ref{fdtheoremlie} and \ref{lowerdiamlie}.
The groups considered in this chapter all satisfy the following hypothesis.
\begin{hypothesis}\label{hypothesiscross}
 Let $G=VG_0$ be a primitive affine group such that $G_0^\infty\cong X_l(r),$ a quasisimple group of Lie type. Suppose that $V$ is an absolutely irreducible $\mathbb{F}_{q_0}G_0^\infty-$module in characteristic $p$ such that $(r,p)=1.$ Also let $n$ be the dimension of $V$  and assume that $V$ cannot be realised over a proper subfield of $\mathbb{F}_{q_0}.$ 
\end{hypothesis}

   

\begin{proof}[Proof of Theorem \ref{fdtheoremlie}]
 Let $\delta_{r'}(G_s)$ denote the minimal dimension of a non-trivial irreducible representation of any covering group of $G_s$ in characteristic not equal to $r,$ so $\delta_{r'}(G_s)\leq n$ in this case. The values of $\delta_{r'}(G_s)$  are in \cite{tiep}\cite{landazuri}. For all $G_s$ as in Hypothesis \ref{hypothesiscross}, we have that $\delta_{r'}(G_s)$ is at least $r^{l/5},$ so we have $n\geq r^{l/5}.$ This is equivalent to $5\log_2(n)\geq l \log_2(r).$ 
 Also the orders of automorphism groups of simple groups of Lie type, are all less than $r^{4l^2+l+2}$ so we have that $\log_2(\vert Aut(G_s)\vert)\leq (4l^2+l+2)\log_2(r).$ Hence by Lemma \ref{eq3} we get the following inequality:   $$n\leq (4l^2+l+2)d\log_2(r)\leq 8dl^2\log_2^2(r) .$$ As $5\log_2(n)\geq l \log_2(r),$   $$n\leq 200d\log_2^2(n).$$  Let $\delta=\frac{\epsilon-1}{\epsilon}.$ For large enough $n,$ this gives $$n\leq dn^\delta$$ which is equivalent to $$n\leq d^{1+\epsilon}.$$

\end{proof}

\begin{proof}[Proof of Theorem \ref{lowerdiamlie}] Lemma \ref{eq3} gives that $n\leq 1+d\log_2(\vert Aut(G_s)\vert)$ and by assumption we have that $\delta_{r'}(G_s)\leq n$. Putting these together we get that $$\delta_{r'}(G_s)\leq 1+d\log_2(\vert Aut(G_s)\vert)$$ and so we get the bound \begin{equation}\label{eqbound}
    d\geq \frac{\delta_{r'}(G_s)-1}{\log_2(\vert Aut(G_s)\vert)}.
\end{equation}


\textbf{1. and 2.} Assume $G_s$ is an exceptional group of Lie type with Lie rank $l$. By \cite{tiep} and \cite{landazuri}, for each family of exceptional groups, the bound in (\ref{eqbound}) for $d$ is larger than $\frac{r^l}{l\log_2(r)}$ for $ G_s\neq X_l(r)\cong ^2B_2(r),$ $^2G_2(r)$ or $^3D_4(r)$ and  larger than $ \frac{r^{l-1}}{(l-1)\log_2(r)}$ for $ G_s\cong X_l(r)\cong ^2B_2(r),$ $^2G_2(r)$ or $^3D_4(r),$ unless  $G_s$ is $^3D_4(2),$ $^2F_4(2)',$ $^2B_2(8),$ $^2B_2(32),$ $G_2(3),$ $G_2(4),$ $G_2(5),$ $G_2(7)$ or $F_4(2).$ The bounds on $d$ in parts 1 and 2 for $G_s\cong G_2(5),$  $G_2(7)$ or $^2B_2(32)$ follow from  (\ref{eqbound}). The lower bound for the remaining exceptions is $2$ as none of them produce examples of $2$-homogeneous affine groups. \par 
\textbf{3.} Assume $G_s$ is a classical group with Lie rank $l$. By \cite{tiep} and \cite{landazuri}, for each family of classical groups, the bound in (\ref{eqbound}) for $d$ is larger than $\frac{r^{l-1}-3}{(l+1)^3\log_2(r)}.$ \par  

\end{proof}

\section{Sporadic Stabilizer }

In this section we prove Theorem \ref{sporadic}.
 
 \begin{proof}[Proof of Theorem \ref{sporadic}]
     (i) and (ii) All possible irreducible representations of dimension less than 250 of these groups are in \cite{hissmalle}. The lower bound of $3$ on $d$ and the values of $N$ follow from Lemma \ref{lemma2.1}(\ref{boundeq1}). \par 
     (ii) These groups give rise to rank 3 affine groups by \cite{liebeckaffine}, so the result follows by Lemma \ref{rankbound}. 
 \end{proof}

\section*{Acknowledgements}
This paper forms part of the PhD thesis of the author under the supervision of Professor Martin Liebeck at Imperial College London, and the author wishes to thank
him for his direction and support throughout. We also thank Aluna Rizzoli for useful discussions and for Lemmas \ref{fullrank} and \ref{alunalemma}. The author was supported by the UK Engineering and Physical Sciences Research Council.

\medskip

\printbibliography

\end{document}